\documentclass[12pt]{article}

\usepackage{amssymb}

\usepackage{amsmath}

\usepackage[mathcal]{euscript}

\usepackage{amsthm}

\newcommand{\XX}{\mathbb{X}}
\newcommand{\YY}{\mathbb{Y}}
\newcommand{\ZZ}{\mathbb{Z}}
\newcommand{\e}{\varepsilon}
\newcommand{\EE}{\mathsf{E}}
\newcommand{\PP}{\mathsf{P}}

\begin{document}

\begin{center}
{\bf \large Asymptotic Expansions \\ for Stationary Distributions  of  \vspace{1mm} \\ Perturbed  Semi-Markov Processes}
\end{center} 
\vspace{2mm}

\begin{center}
{\large Dmitrii Silvestrov\footnote{Department of Mathematics, Stockholm University, SE-106 81 Stockholm, Sweden. \\ 
Email address: silvestrov@math.su.se} 
and Sergei Silvestrov\footnote{Division of Applied Mathematics, School of Education, Culture and Communication, M{\"a}lardalen  University, SE-721 23 V{\"a}ster{\aa}s, Sweden. \\ 
Email address: sergei.silvestrov@mdh.se}}
\end{center}
\vspace{2mm}

Abstract:
New algorithms for computing of asymptotic expansions for stationary distributions  of nonlinearly perturbed semi-Markov processes are presented. 
The algorithms are based on special techniques of sequential phase space reduction, which can be applied to processes  with asymptotically coupled and uncoupled finite phase spaces.  \\
 


{\bf 1. Introduction}  \\

In this paper, we present new algorithms  for  construction asymptotic expansions for stationary distributions of nonlinearly perturbed semi-Markov processes with a finite phase space.  

We consider models, where  the phase space of embedded Markov chains for pre-limiting perturbed semi-Markov processes is one class of communicative states, while the phase space for the limiting embedded Markov chain can consist of one or several closed classes of communicative states and, possibly,  a class of transient states. 

The initial perturbation conditions are formulated  in the forms of Taylor asymptotic expansions for transition probabilities of the corresponding embedded Markov chains and Laurent asymptotic  expansions for expectations of sojourn times for perturbed semi-Markov processes. Two forms of these expansions are considered, with remainders given  without or with explicit upper bounds.

The algorithms are based  on special time-space screening procedures for sequential phase space reduction and algorithms for  re-calculation of asymptotic expansions and upper bounds for remainders, which constitute perturbation conditions for the  semi-Markov processes with reduced phase spaces. 

The final asymptotic expansions for stationary distributions of nonlinearly perturbed semi-Markov processes are given in the form of Taylor asymptotic expansions with remainders given without or with explicit upper bounds.

The model of perturbed Markov chains  and semi-Markov processes, in particular, in the most difficult case of so-called singularly perturbed  Markov chains and semi-Markov processes  with absorption and asymptotically uncoupled phase spaces, attracted attention of researchers in the mid of the 20th century.  

The first  works related to asymptotical problems for the above models are Meshalkin (1958),  Simon, and Ando (1961), Hanen (1963a, b, c, d), Seneta (1967, 1968a, b), Schweitzer (1968), and Korolyuk (1969). 

The methods used for construction of asymptotic expansions for stationary distributions and related functionals such as moments of hitting times can be split in three groups.  

The first and the most widely used methods are based on analysis of  generalized  matrix and operator inverses of resolvent type for transition matrices and operators  for singularly perturbed Markov chains and semi-Markov processes. Mainly models with linear, polynomial and analytic perturbations have been objects of studies. We refer here to works by  Schweitzer (1968), Turbin (1972),  Poli\v s\v cuk and  Turbin (1973), Koroljuk, Brodi and Turbin (1974), Pervozvanski\u\i \, and  Smirnov (1974), Courtois  and  Louchard (1976),  Korolyuk and Turbin  (1976, 1978), Courtois  (1977), Latouche and  Louchard  (1978), Kokotovi\'{c},  Phillips and  Javid (1980), Korolyuk, Penev and  Turbin (1981),  Phillips and  Kokotovi\'{c} (1981), Delebecque (1983),  Abadov (1984), Silvestrov and  Abadov (1984, 1991, 1993), Kartashov (1985d, 1988, 1996b), Haviv (1986),  Korolyuk (1989),   Stewart and  Sun (1990),  Haviv, Ritov and  Rothblum (1992), Haviv and  Ritov (1993),  Schweitzer and Stewart (1993),   Stewart (1998, 2001), Yin and  Zhang (1998, 2003, 2005, 2013), Avrachenkov (1999, 2000), Avrachenkov and  Lasserre (1999), Korolyuk, V.S. and Korolyuk, V.V. (1999),  Yin, G., Zhang,  Yang and Yin, K. (2001), Avrachenkov and  Haviv (2003, 2004), Craven (2003), Bini, Latouche and Meini (2005), Korolyuk and  Limnios (2005)  and   Avrachenkov,  Filar and Howlett  (2013).

Alternatively, the methods based on regenerative properties of Markov chains and semi-Markov processes, in particular,  relations which link stationary probabilities and expectations of return times have been used for getting approximations for expectations of hitting times and stationary distributions in works by Grassman, Taksar and Heyman (1985),  Hassin and  Haviv (1992) and Hunter (2005). Also, the above mentioned relations and  methods  based on asymptotic expansions for nonlinearly perturbed regenerative processes developed in works by Silvestrov (1995, 2007, 2010),  Englund and  Silvestrov (1997), Gyllenberg and Silvestrov (1998, 1999a, 2000a, 2008), Englund (2000, 2001), Ni, Silvestrov and Malyarenko (2008), Ni (2010a, b, 2011, 2012, 2014), Silvestrov and Petersson (2013) and Petersson (2013a, b, 2014) have been  used for getting asymptotic expansions for stationary and quasi-stationary distributions for nonlinearly perturbed Markov chains and semi-Markov processes with absorption.  

Aggregation/disaggregation methods based on various modification of Gauss elimination method and space screening procedures for perturbed Markov chains have been employed for approximation of stationary distributions for Markov chains in works by Coderch, Willsky, Sastry and Casta\~{n}on (1983), Delebecque (1983), Ga\u\i tsgori and  Pervozvanski\u\i \, (1983), Chatelin and  Miranker (1984),  Courtois and Semal  (1984b), Seneta (1984), Cao and  Stewart  (1985), Vantilborgh (1985), Feinberg and  Chiu (1987),   Haviv (1987, 1992, 1999), Sumita and  Reiders (1988), Meyer (1989), Schweitzer (1991), Stewart and Zhang (1991), Stewart (1993a), Kim and  Smith (1995) and Avrachenkov,  Filar and Howlett  (2013), Silvestrov, D, and Silvestrov S. (2015).

We would like to mention that the present paper  contains also  a more extended bibliography of works in the area supplemented by short bibliographical remarks given in the last section of the paper.  

In the present paper, we combine methods based on stochastic aggregation/disaggregation approach with  methods  based on asymptotic expansions for perturbed regenerative processes applied to perturbed  semi-Markov processes.

In the above mentioned works  based on stochastic aggregation/disaggrega\-tion approach, space screening procedures for discrete time Markov chains are  used. A Markov chain with a reduced phase space  is constructed from the initial one as the sequence of its states  at sequential moment of hitting into the reduced phase space. Times between sequential hitting of a reduced phase space are ignored. Such screening procedure preserves ratios of hitting frequencies for states from the reduced phase space and, thus, the ratios of stationary probabilities are the same for the  initial and the reduced  Markov chains. This implies that the stationary probabilities for the reduced Markov chain coincide with the corresponding stationary probabilities for the initial Markov chain up to the change of the corresponding normalizing factors.

We use another more complex type of time-space screening procedures for semi-Markov processes. In this case, a semi-Markov process with a reduced phase space  is constructed from the initial one as the sequence of its states  at sequential moment of hitting into the reduced phase space and times between sequential jumps  of the  reduced semi-Markov process are  times between sequential hitting of the reduced space by the initial  semi-Markov process.  Such screening procedure preserves transition times between states from the reduced phase space, i.e., these times and, thus, their expectations are the same for the initial and the reduced semi-Markov processes. 

We also formulate perturbation conditions in terms of asymptotic expansions for transition characteristics of perturbed semi-Markov processes. The remainders in these expansions and, thus, the transition characteristics of perturbed semi-Markov processes can be non-analytical functions of perturbation parameters that makes difference with the results for models with linear, polynomial and analytical perturbations. 

We employ the methods of asymptotic analysis for nonlinearly perturbed regenerative processes developed  in works by   Silvestrov (1995, 2007, 2010) and  Gyllenberg and Silvestrov (1998, 1999a, 2000a, 2008). However, we use techniques of more general Laurent asymptotic expansions instead of Taylor  asymptotic expansions used in the  above mentioned works and combine these methods with the aggregation/disaggregation approach instead of using the approach based on generalized  matrix  inverses. This permits us consider perturbed semi-Markov processes with an arbitrary communication structure of the phase space for the  limiting semi-Markov process, including the general case, where this phase space may consist from one or several closed classes of communicative states and possibly a class of transient states.  

Another new element is that we consider asymptotic expansions with remainders given not only in the form $o(\cdot)$, but, also,  with explicit upper bounds.  

It should be mentioned that the semi-Markov setting is an adequate and  necessary element of the  method proposed in the paper. Even in the case, where the initial process is a discrete or continuous time Markov chain,  the time-space screening procedure of phase space reduction results in a semi-Markov process,  since times between sequential hitting of the reduced space by the initial process have distributions  which can differ of geometrical or exponential ones. 

Also, the use of Laurent  asymptotic expansions for expectations of sojourn times of perturbed semi-Markov processes is also a necessary element of the method. Indeed, even in the case, where expectations of sojourn times for all states  of the initial semi-Markov process are asymptotically bounded and represented by Taylor asymptotic expansions, the   exclusion of an asymptotically absorbing state from the initial  phase space  can cause appearance of states with asymptotically unbounded  expectations of sojourn times represented by Laurent  asymptotic expansions,  for the reduced semi-Markov processes.

The method proposed in the paper can be considered as a stochastic analogue of the Gauss elimination method. It is based on the procedure of sequential exclusion of states  from the phase space of a perturbed semi-Markov process accompanied by re-calculation of  asymptotic expansions penetrating perturbation conditions  for semi-Markov processes with reduced phase spaces. The corresponding  algorithms are based on some kind of  ``operational calculus''  for Laurent  asymptotic 
expansions with remainders given in two forms, without or with explicit upper bounds. 

The corresponding computational algorithms have an universal character. As was mentioned above, they can be applied to perturbed semi-Markov processes with an arbitrary communicative structure and are computationally effective due to recurrent character of computational procedures.

In conclusion, we would like to point out that, by our opinion, the results presented in the paper have a good potential for continuation of studies (asymptotic expansions for high order power and exponential moments for hitting times, aggregated time-space screening procedures, asymptotic expansions for quasi-stationary distributions, etc.). We comment some prospective directions for future studies in the end of the paper. 

The paper includes 8 sections. In Section 2, we present so-called operational rules for  Laurent asymptotic expansions. In Section 3, we  formulate basic perturbation conditions for Markov chains and semi-Markov processes. In Section 4, we give basic formulas for stationary distributions  for semi-Markov processes, in particular, formulas connecting  stationary distributions with expectations of return times. In Section 5, we present an one-step procedure of phase space reduction for semi-Markov processes. In Section 6, we present algorithms for re-calculation of asymptotic expansions for transition characteristics of perturbed semi-Markov processes with  a reduced phase space. In Section 7, we present algorithms of sequential reduction of phase space for semi-Markov processes and construction of 
asymptotic expansions for stationary distributions. In Section 8, we present some directions for future studies and short bibliographical 
remarks concerned works in the area. 

We would like to conclude the introduction with the remark that the present paper is a slightly improved version of the research report Silvestrov, D. and Silvestrov S. (2015). \\

{\bf 2. Laurent asymptotic expansions} \\

In this section, we present so-called operational rules for  Laurent asymptotic expansions. We consider  the corresponding results as possibly known, except,  some of explicit formulas for remainders, in particular, those related to product, reciprocal and quotient rules. \vspace{1mm}

{\bf 2.1. Definition of Laurent asymptotic expansions}.  Let $A(\e)$ be a real-valued function 
defined on an interval $(0, \e_0]$, for some $0 < \e_0 \leq 1$, and given on this interval
by a  Laurent asymptotic expansion,
\begin{align}\label{expap}
A(\e) = a_{h_A}\e^{h_A} + \cdots + a_{k_A}\e^{k_A} + o_A(\e^{k_A}),
\end{align}
where {\bf (a)} $- \infty < h_A \leq k_A < \infty$ are integers, {\bf
(b)} the coefficients $a_{h_A}, \ldots, a_{k_A}$ are real numbers, {\bf (c)} 
function $o_A(\e^{k_A})/\e^{k_A} \rightarrow 0$ as $\e \rightarrow 0$. 

We refer to such  Laurent asymptotic expansion as a $(h_A, k_A)$-expansion.

We say that  $(h_A, k_A)$-expansion $A(\e)$ is pivotal if it is known that $a_{h_A} \neq 0$.

We also say that $(h_A, k_A)$-expansion $A(\e)$ is a $(h_A, k_A, \delta_A, G_A,  \e_A)$-expansi\-on if  its remainder 
$o_A(\e^{k_A})$ satisfies the following inequalities  {\bf (d)} $|o_A(\e^{k_A})|$ $\leq G_A \e^{k_A + \delta_A}$,  for $0 < \e \leq \e_A$, where {\bf (e)}  $0 < \delta_A \leq 1, 0 < G_A < \infty$ and  $0 < \e_A \leq \e_0$.

In what follows, $[a]$ is the integer part of a real number $a$. 

Also, the indicator of relation $A = B$ is denoted as $I(A = B)$. It equals to $1$, if $A = B$, or $0$, if $A \neq B$.   

It is useful to note that there is no sense to consider, it seems,  a more general case of upper bounds for the remainder $o_A(\e^{k_A})$, with parameter $\delta_A > 1$. Indeed, let us define $k'_A =  k_A +  [\delta_A] - I(\delta_A = [\delta_A])$ and $\delta'_A = \delta_A -  [\delta_A] + I(\delta_A = [\delta_A])
\in (0, 1]$. 

The $(h_A, k_A, \delta_A, G_A,  \e_A)$-expansion (\ref{expap}) can be re-written in the equivalent form of the  
$(h_A, k'_A, \delta'_A, G_A, \e_A)$-expansion,
\begin{align}\label{expapas}
A(\e) = a_{h_A}\e^{h_A} + \cdots + a_{k_A}\e^{k_A} + 0 \e^{k_A +1} + \cdots + 0 \e^{k'_A} +  o'_A(\e^{k'_A}),
\end{align}
with the remainder term $o'_A(\e^{k'_A}) = o_A(\e^{k_A})$,  which satisfies inequalities $|o'_A(\e^{k'_A})| = |o_A(\e^{k_A})| \leq 
 G_A \e^{k_A + \delta_A} = G_A \e^{k'_A + \delta'_A}$, for $0 < \e \leq \e_A$.
 
Relation (\ref{expapas}) implies that the  asymptotic expansion $A(\e)$ can be represented in different forms. In such cases, we consider a more informative form with larger parameters $h_A$ and $k_A$. As far as parameters  $\delta_A, G_A$ and $\e_A$ are concerned, we consider as a more informative form, first, with larger value of parameter $\delta_A$, second, with smaller values of parameter $G_A$  and, third, with the larger values of parameter $\e_A$.   

In what follows, $a \vee b = \max(a, b)$ and $a \wedge b = \min(a, b)$,  for real numbers $a$ and $b$.
 
It is useful to note that formula (\ref{expap}) uniquely define part of coefficients $a_{h_A}, \ldots, a_{k_A}$. \vspace{1mm}

{\bf Lemma 1.} {\em If function $A(\e) = 
a'_{h'_A}\e^{h'_A} + \cdots + a'_{k'_A} \e^{k'}_A + o'_A(\e^{k'_A}) =  a''_{h''_A}\e^{h''_A} + \cdots + a''_{k''_A} \e^{k''_A} + o''_A(\e^{k''_A}), \e \in (0, \e_0]$ can be represented as, respectively, $(h'_A, k'_A)$- and $(h''_A, k''_A)$-expansion,  then the asymptotic expansion for function $A(\e)$ can be represented in the following the most informative form $A(\e) = a_{h_A}\e^{h_A} + \cdots + a_{k_A}\e^{k_A} + o_A(\e^{k_A}), \e \in (0, \e_0]$ of $(h_A, k_A)$-expansion, with parameters $h_A = h'_A \vee h''_A, k_A = k'_A \vee k''_A$, and  coefficients $a_{h_A}, \ldots, a_{k_A}$ and remainder $o_A(\e^{k_A})$ given by the following relations{\rm :}
\begin{itemize}
\item[{\rm \bf (i)}] $a'_l, a''_l = 0$, for $l < h_A$.

\item[{\rm \bf (ii)}] $a_l = a'_l = a''_l$, for $h_A  \leq l \leq  \tilde{k}_{A} = k'_A \wedge k''_A$. 

\item[{\rm \bf (iii)}] $a_l = a''_l$,  for $\tilde{k}_{A}  = k'_A < l \leq k_A$ if $k'_A < k''_A$.

\item[{\rm \bf (iv)}]  $a_l = a'_l$, for $\tilde{k}_{A} = k''_A  < l \leq k_A$ if $k''_A < k'_A$. 

\item[{\rm \bf (v)}]  The remainder term $o_A(\e^{k_A})$ is given by the following relation,
\begin{equation}\label{reequar}
o_A(\e^{k_A}) = \left\{
\begin{array}{ll}
o''_A(\e^{k''_A})  & \ \text{if} \   k'_A < k''_A, \\

o'_A(\e^{k'_A}) = o''_A(\e^{k''_A}) & \ \text{if} \   k'_A = k''_A, \\

o'_A(\e^{k'_A})  & \ \text{if} \   k'_A > k''_A.

\end{array}
\right.
\end{equation}  
\end{itemize}

The latter asymptotical expansion is pivotal if and only if $a_{h_A}
 = a'_{h_A} = a''_{h_A} \neq  0$. 
} \vspace{1mm}

It is useful to make some additional remarks.

The case $\tilde{k}_A < h_A$ is possible. In this case, the set of integers $l$ such that  $h_A  \leq l \leq \tilde{k}_A$ is empty. This can happen if $k'_A < h''_A$ or $k''_A < h'_A$. In the first case, all coefficients $a'_l = 0, l = h'_A, \ldots, k'_A$ while $h_A = h''_A, k = k''_A, a''_l = a_l, l = h_A, \ldots, k_A$. In the second case, all coefficients $a''_l = 0, l = h''_A, \ldots, k''_A$ while $h_A = h'_A, k_A = k'_A, a_l = a'_l, l = h_A, \ldots, k_A$.

If $k'_A = k''_A$ then $h_A \leq \tilde{k}_A = k_A$ and the set of integers $l$ such that $\tilde{k}_A < l \leq k_A$ is empty. In this case, all coefficients $a_l = a'_l = a''_l, l = h_A, \ldots, k_A$.

If $a'_{h_A} \neq 0$ then $h_A = h'_A$ and  $a_{h_A} = a'_{h'_A} \neq 0$. If  $a''_{h''_A} \neq 0$ then $h_A = h''_A$ and  $a_{h_A} = a''_{h''_A} \neq 0$. If $a'_{h'_A},  a''_{h''_A} \neq 0$ then $h_A = h'_A = h''_A$ and  $a_{h_A} = a'_{h'_A} = a''_{h''_A}  \neq 0$.

The following proposition supplements Lemma 1. \vspace{1mm}

{\bf Lemma 2.} {\em If $A(\e) = 
a'_{h'_A}\e^{h'_A} + \cdots + a'_{k'_A} \e^{k'_A} + o'_A(\e^{k'_A}) =  a''_{h''_A}\e^{h''_A} + \cdots + a''_{k''_A} \e^{k''_A} + o''_A(\e^{k''_A}), \e \in (0, \e_0]$ can be represented as, respectively,  $(h'_A, k'_A,  \delta'_A, G'_A$, $\e'_A)$- and $(h''_A, k''_A, \delta''_A, G''_A, \e''_A)$-expansion, then{\rm :}
\begin{itemize}
\item[{\rm \bf (i)}] The asymptotic expansion $A(\e) = a_{h_A}\e^{h_A} + \cdots + a_{k_A}\e^k + o_A(\e^{k_A}), \e \in (0, \e_0]$  given in Lemma 1 is an  $(h_A, k_A, \delta_A$, $G_A, \e_A)$-expansion with parameters $G_A, \delta_A$ and $\e_A$ which can be chosen in the following way consistent with the priority order described above{\rm :}
{\footnotesize
\begin{equation}\label{equar}
(\delta_A, G_A, \e_A) = \left\{
\begin{array}{ll}

(\delta''_A, G''_A,  \e''_A)  & \ \text{if} \   k'_A < k''_A, \vspace{1mm} \\

(\delta''_A, G''_A,  \e''_A)  & \ \text{if} \  k'_A = k''_A, \delta'_A < \delta''_A,  \vspace{1mm} \\

(\delta'_A = \delta''_A, G'_A \wedge G''_A,  \e'_A \wedge \e''_A) & \ \text{if} \   k'_A = k''_A,  \delta'_A = \delta''_A, \vspace{1mm} \\

(\delta'_A,  G'_A,  \e'_A)  & \ \text{if} \ k'_A = k''_A, \delta'_A > \delta''_A, \vspace{1mm} \\

(\delta'_A,  G'_A,  \e'_A)   & \ \text{if} \   k' _A > k''_A.

\end{array}
\right.
\end{equation} 
}

\item[{\rm \bf (ii)}]  The asymptotic expansion  $A(\e)$ can also be represented in the form $A(\e) = a_{\tilde{h}_A}\e^{\tilde{h}_A} + \cdots + a_{\tilde{k}_{A}} \e^{\tilde{k}_A} 
+ \tilde{o}_A(\e^{\tilde{k}_A})$ of an  $(\tilde{h}_A, \tilde{k}_A, \tilde{\delta}_A, 
\tilde{G}_A, \tilde{\e}_A)$-expansion, with parameters $\tilde{h}_A = h_A,  \tilde{k}_A = k'_A \wedge k''_A$  and parameters $ \tilde{\delta}_A, 
\tilde{G}_A, \tilde{\e}_A$ given by the following formulas, 
{\footnotesize
\begin{equation*}
\tilde{\delta}_A = \left\{
\begin{array}{ll}
\delta'_A &  \ \text{if} \ k'_A < k''_A , \vspace{1mm} \\
\delta'_A \wedge \delta''_A  & \ \text{if} \  k'_A = k''_A,  \makebox[66mm]{} \\
\delta''_A   & \ \text{if} \   k'_A > k''_A,
\end{array}
\right.
\end{equation*}
\begin{equation*}
\tilde{G}_A = \left\{
\begin{array}{ll}
G'_A \wedge ( \sum_{k'_A < l \leq k''_A} |a''_l | \tilde{\e}_A^{l - k'_A - \delta'_A}  + G''_A \tilde{\e}_A^{k''_A +\delta''_A - k'_A - \delta'_A})   & \ \text{if} \   k'_A < k''_A, \vspace{2mm} \\

G'_A \tilde{\e}_A^{\, \delta'_A  - \tilde{\delta}_A} \wedge   G''_A \tilde{\e}_A^{\, \delta''_A  - \tilde{\delta}_A} & \ \text{if} \  k'_A = k''_A,   \vspace{2mm}\\

G''_A \wedge ( \sum_{k''_A < l \leq k'_A} |a'_l | \tilde{\e}_A^{l - k''_A - \delta''_A}  + G'_A \tilde{\e}_A^{k'_A +\delta'_A - k''_A - \delta''_A})   & \ \text{if} \   k'_A> k''_A. \vspace{2mm}

\end{array}
\right.
\end{equation*} 
\begin{equation}
\tilde{\e}_A = \e'_A \wedge \e''_A. \makebox[92mm]{}
\end{equation} 

}

\item[{\rm \bf (iii)}] The remainders   $o'_A(\e^{k'_A}), o''_A(\e^{k''_A}), o_A(\e^{k_A})$ and $\tilde{o}_A(\e^{k_A})$ are connected by the following relations{\rm :}
\begin{equation*}
\tilde{o}_A(\e^{\tilde{k}_A}) = o_A(\e^{k_A}) + \sum_{\tilde{k}_A < l \leq k_A} a_l \e^l  \makebox[40mm]{}
\end{equation*} 
\begin{equation}
 = \left\{
\begin{array}{ll}
o'_A(\e^{k'_A})  &  \ \text{if} \ k'_A < k''_A, \vspace{1mm} \\
o'_A(\e^{k'_A})  = o''_A(\e^{k''_A})  & \ \text{if} \ k'_A = k''_A, \vspace{1mm} \\
 o''_A(\e^{k''_A}) & \ \text{if} \ k'_A > k''_A. 
\end{array}
\right.
\end{equation} 

\end{itemize} 
}
\vspace{1mm}

{\bf 2.2. Operational rules for  Laurent asymptotic expansions}.
Let us consider four Laurent asymptotic  expansions, $A(\e) = a_{h_A}\e^{h_A} +
\cdots + a_{k_A}\e^{k_A} + o_A(\e^{k_A})$, $B(\e) = b_{h_B}\e^{h_B} +
\cdots + b_{k_B}\e^{k_B} + o_B(\e^{k_B})$, $C(\e) = c_{h_C}\e^{h_C} +
\cdots + c_{k_C}\e^{k_C} + o_C(\e^{k_C})$, and $D(\e) =
d_{h_D}\e^{h_D} + \cdots + d_{k_D}\e^{k_D} + o_D(\e^{k_D})$ defined for $0 < \e \leq \e_0$, for some $0 < \e_0 \leq 1$.

The following lemma presents   ``operational'' rules for Laurent asymptotic  expansions. \vspace{1mm}

{\bf Lemma 3.} {\em The above asymptotic expansions
have the following operational rules for computing coefficients{\rm :}
\begin{itemize}
\item[{\rm \bf (i)}] If $A(\e), \e \in (0, \e_0]$ is a $(h_A, k_A)$-expansion
  and $c$ is a constant, then $C(\e) = cA(\e), \e \in (0, \e_0]$ is a $(h_{C},
  k_{C})$-expansion with parameters $h_{C} = h_A, k_{C} = k_A$ and  coefficients,  
  \begin{equation}
  c_{h_C +r} = c a_{h_C + r},  r = 0,  \ldots,   k_C - h_C.
  \end{equation}  
  
  This expansion is pivotal if and only if $c_{h_{C}} = c a_{h_A}  \neq 0$.

\item[{\rm \bf (ii)}] If $A(\e), \e \in (0, \e_0]$ is a $(h_A, k_A)$-expansion
  and $B(\e), \e \in (0, \e_0]$ is a $(h_B, k_B)$-expansion, then $C(\e) =
  A(\e)+ B(\e), \e \in (0, \e_0]$ is a $(h_C, k_C)$-expansion with parameters $h_C = h_A
  \wedge h_B, k_C = k_A \wedge k_B$, and  coefficients,
  \begin{equation}
  c_{h_C +r} = a_{h_C +r} + b_{h_C +r}, r = 0, \ldots,
  k_C - h_C,
  \end{equation} 
  where $a_{h_C +r} = 0$ for $0 \leq r < h_A - h_C$ and
  $b_{h_C +r} = 0$ for $0 \leq r < h_B - h_C$. 
  
  This expansion 
  is pivotal if and only if $c_{h_C} = a_{h_C} + b_{h_C} \neq 0$.

\item[{\rm \bf (iii)}] If $A(\e), \e \in (0, \e_0]$ is a $(h_A,
  k_A)$-expansion and $B(\e), \e \in (0, \e_0]$ is a $(h_B, k_B)$-expansion,
  then $C(\e) = A(\e) \cdot B(\e), \e \in (0, \e_0]$ is a $(h_C,
  k_C)$-expansion with parameters $h_C = h_A + h_B, k_C = (h_A + k_B) \wedge (h_B + k_A)$, and  coefficients,
  \begin{equation}
  c_{h_C + r} =
  \sum_{0 \leq i \leq r} a_{h_A + i} b_{h_B + r - i}, r = 0, \ldots,
  k_C - h_C.
  \end{equation} 
  
  This expansion is pivotal if and only if $c_{h_C}
  = a_{h_A} b_{h_B} \neq 0$.

\item[{\rm \bf (iv)}] If $B(\e), \e \in (0, \e_0]$ is a pivotal $(h_B, k_B)$-expansion, 
then there exists $0 < \e'_0 \leq \e_0$ such that $B(\e) \neq 0, \e \in (0, \e'_0]$,  and $C(\e) = 1/B(\e), \e \in (0, \e'_0]$ is a pivotal 
$(h_C, k_C)$-expansion with parameters
  $h_C = - h_B$, $k_C = k_B  - 2h_B$ and coefficients,
 \begin{equation}
c_{h_C} = b_{h_B}^{-1}, \  c_{h_C + r} = - b_{h_B}^{-1} \sum_{1 \leq i \leq r}
  b_{h_B + i} c_{h_C + r - i}, r = 1, \ldots, k_C - h_C.
  \end{equation}

  \item[{\rm \bf (v)}] If $A(\e), \e \in (0, \e_0]$ is a $(h_A, k_A)$-expansion
  $B(\e), \e \in (0, \e_0]$  is a pivotal $(h_B, k_B)$-expansion, then, there exists 
  $0 < \e'_0 \leq \e_0$ such that $B(\e) \neq 0, \e \in (0, \e'_0]$, and $D(\e) = 
  A(\e)/B(\e), \e \in (0, \e'_0]$ is a $(h_D, k_D)$-expansion with parameters $h_D = h_A - h_B,
  k_D = (k_A - h_B)  \wedge (h_A + k_B - 2 h_B)$, and coefficients,
  \begin{equation}\label{first}
  d_{h_D + r} = \sum_{0 \leq i \leq r} c_{h_C +i} a_{h_A + r - i}, r = 0, \ldots,
  k_D - h_D, 
  \end{equation}
  where $c_{h_C +j} , j = 0, \ldots, k_C - h_C$\, are
  coefficients of the $(h_C, k_C)$-expansion $C(\e) =
 1/B(\e)$ given in the above proposition {\rm \bf (iv)}, or by formulas,
 \begin{equation}\label{firsta}
 d_{h_D + r} = b_{h_B}^{-1} (a_{h_A + r} - \sum_{1
    \leq i \leq r} b_{h_B + i} d_{h_D + r - i}), r = 0, \ldots, k_D - h_D.
  \end{equation}

This expansion is
pivotal if and only if $d_{h_D} = a_{h_A}c_{h_C} = a_{h_A}/ b_{h_B} \neq 0$.
  
\end{itemize} 
} 

The following proposition presents   ``operational''  rules for computing parameters of upper bounds for remainders of  Laurent asymptotic expansions.  

\vspace{1mm}

{\bf Lemma 4.} {\em The above asymptotic expansions
have the following operational rules for computing remainders{\rm :}
\begin{itemize}
\item[{\rm \bf (i)}] If $A(\e), \e \in (0, \e_0]$ is a $(h_A, k_A,  \delta_A,  G_A,\e_A)$-expansion
  and $d$ is a constant, then $C(\e) = cA(\e), \e \in (0, \e_0]$ is a $(h_C, k_C, \delta_C, G_C,  \e_C)$-expansion with parameters $h_C = h_A, k_C = k_A$, coefficients $c_r, r = h_C,
  \ldots, k_C$  given in proposition  {\bf (i)} of Lemma 3, and parameters $\delta_C, G_C,  \e_C$ given by the following formulas,  
\begin{equation}
\delta_C = \delta_A, \  G_C = |c|G_A,  \ \e_C =  \e_A.
\end{equation}

\item[{\rm \bf (ii)}] If $A(\e), \e \in (0, \e_0]$ is a $(h_A, k_A, \delta_A,  G_A, \e_A)$-expansion
  and $B(\e), \e \in (0, \e_0]$ is a $(h_B, k_B$, $\delta_B, G_B,  \e_B)$-expansion, then $C(\e) =
  A(\e)+ B(\e)$, $\e \in (0, \e_0]$ is a $(h_C, k_C,  \delta_C, G_C$, $\e_C)$-expansion with parameters $h_C = h_A \wedge h_B,  k_C = k_A \wedge k_B$, coefficients $c_{r}, r = h_C, \ldots, k_C$ given in proposition  {\bf (ii)} of Lemma 3, and parameters $\delta_C, G_C,  \e_C$ given by the following formulas, 
\begin{equation*}
\delta_C = \left\{
\begin{array}{ll}
\delta_A & \ \text{if} \ k_C = k_A < k_B, \\
\delta_A \wedge \delta_B & \ \text{if} \ k_C = k_A = k_B,   \makebox[21mm]{}  \\
\delta_B & \ \text{if} \ k_C = k_B <  k_A,
\end{array} 
\right. 
\end{equation*} 
\begin{equation*}
\geq \delta_A \wedge \delta_B,   \makebox[57mm]{}  
\end{equation*} 
\begin{align*}
G_C = & \ G_A \e_C^{k_A + \delta_A - k_C - \delta_C}  + \sum_{k_C < i \leq k_A} |a_i|\e_C^{i - k_C - \delta_C}   \nonumber \\
&  + \ G_B \e_C^{k_B + \delta_B - k_C - \delta_C}  +  \sum_{k_C < j \leq k_B} |b_j|\e_C^{j - k_C - \delta_C},  
\end{align*}  
\begin{equation}\label{profac}
\e_C = \e_A \wedge \e_B.  \makebox[63mm]{} 
\end{equation}
 
\item[{\rm \bf (iii)}]  If $A(\e), \e \in (0, \e_0]$ is a $(h_A, k_A, \delta_A, G_A, \e_A)$-expansion
  and $B(\e), \e \in (0, \e_0]$ is a $(h_B, k_B$, $\delta_B, G_B, \e_B)$-expansion, 
  then $C(\e) = A(\e) \cdot B(\e)$, $\e \in (0, \e_0]$ is a $(h_C, k_C,  \delta_C$, $G_C, \e_C)$-expansion with parameters $h_C = h_A + h_B, k_C =  (h_A + k_B) \wedge$ $(h_B +k_A)$, coefficients $c_{r}, r = h_C, \ldots, k_C$ given in proposition  {\bf (iii)} of Lemma 3, and parameters $\delta_C, G_C, \e_C$ given by the following formulas,
\begin{equation*}
\delta_C = \left\{
\begin{array}{ll}
\delta_A & \ \text{if} \ k_C = h_B+ k_A < h_A + k_B, \\
\delta_A \wedge \delta_B & \ \text{if} \ k_C = h_B+ k_A = h_A + k_B,   \\
\delta_B & \ \text{if} \ k_C = h_A+ k_B < h_B + k_A,
\end{array}
\right. 
\end{equation*}
\begin{equation*}
\geq \delta_A \wedge \delta_B,   \makebox[54mm]{}  
\end{equation*} 
\begin{align*}
G_C = & \sum_{k_C < i + j, h_A \leq i \leq k_A, h_B \leq j \leq  k_B} |a_i| |b_j| \e_C^{i + j -  k_C - \delta_C} \makebox[5mm]{} \nonumber \\
& + G_A \sum_{h_B \leq j \leq k_B} |b_j| \e_C^{j + k_A + \delta_A  - k_C  - \delta_C}   \nonumber \\
& + G_B  \sum_{h_A \leq i \leq k_A} |a_i|  \e_C^{i + k_B  + \delta_B - k_C  - \delta_C} \nonumber \\
& + G_A  G_B  \e_C^{k_A+ k_B  + \delta_A  +\delta_B - k_C  - \delta_C}, 
\end{align*}
\begin{equation}\label{profa}
\e_C = \e_A \wedge \e_B. \makebox[61mm]{} 
\end{equation} 

\item[{\rm \bf (iv)}] If $B(\e), \e \in (0, \e_0]$ is a pivotal $(h_B, k_B, \delta_B, G_B, \e_B)$-expansion, then, there exist 
$\e_C \leq  \e'_0 \leq \e_0$ such that  $B(\e) \neq 0, \e \in (0, \e'_0]$, and $C(\e) = 1/B(\e), \e \in (0, \e'_0]$ is a pivotal $(h_C, k_C, \delta_C, G_C, \e_C)$-expansion with parameters $h_C = - h_B$, $k_C = k_B  - 2h_B$, coefficients $c_{r}, r = h_C, \ldots, k_C$ given in proposition  {\bf (iv)} of Lemma 3, and parameters $\delta_C, G_C, \e_C$ given by the following formulas,
\begin{equation*}
\delta_C = \delta_B, \makebox[94mm]{} 
\end{equation*}
\begin{align*}
G_C = & \ (\frac{|b_{h_B}|}{2})^{-1} \big( \sum_{k_B - h_B < i + j, h_B \leq i \leq k_B, h_C \leq j \leq  k_C} |b_i| |c_j| \e_C^{i + j - k_B + h_B - \delta_B} 
\nonumber \\
& \ + G_B \sum_{h_C \leq j \leq k_C} |c_j| \e_C^{j + h_B} \big), 
\end{align*}
\begin{equation}\label{naswe}
\e_C = \e_B \wedge \left\{
\begin{array}{ll}
 \frac{|b_{h_B}|}{2} ( \sum_{h_B < i \leq k_B} |b_i| \e_B^{i - h_B - 1} \vspace{2mm} \\
+  \ G_B \e_B^{k_B + \delta_B - h_B - 1})^{-1}  & \ \text{if} \ h_B < k_B,  \makebox[5mm]{}  \vspace{2mm} \\
( \frac{|b_{h_B}|}{2G_B})^{\frac{1}{\delta_B}}  & \ \text{if} \ h_B = k_B.
\end{array}
\right.
\end{equation}

  \item[{\rm \bf (v)}] If $A(\e), \e \in (0, \e_0]$ is a $(h_A, k_A,  \delta_A, G_A, \e_A)$-expansion,
  $B(\e), \e \in (0, \e_0]$ is a pivotal $(h_B, k_B, \delta_B, G_B, \e_B)$-expansion, then, there exist 
  $\e_D \leq  \e'_0 \leq \e_0$ such that  $B(\e) \neq 0, \e \in (0, \e'_0]$, and  $D(\e) = 
  A(\e)/B(\e)$ is a $(h_D, k_D$,  $\delta_D, G_D$, $\e_D)$-expansion with parameters $h_D = h_A + h_C = h_A - h_B,
  k_D = (k_A + h_C) \wedge (h_A + k_C) = (k_A - h_B)  \wedge (h_A + k_B - 2 h_B)$,  coefficients $d_{r}, r = h_D, \ldots, k_D$ given in proposition  {\bf (v)} of Lemma 3, 
  and parameters $\delta_D, G_D, \e_D$ given by the following formulas,
\begin{equation*}
\delta_D = \left\{
\begin{array}{ll}
\delta_A & \ \text{if} \ k_D = h_C+ k_A < h_A + k_C \\
\delta_A \wedge \delta_C & \ \text{if} \ k_D = h_C+ k_A = h_A + k_C,    \\
\delta_C & \ \text{if} \ k_D = h_A+ k_C < h_C + k_A,
\end{array}
\right. 
\end{equation*}
\begin{equation*}
\geq \delta_A \wedge \delta_C = \delta_A \wedge \delta_B,   \makebox[35mm]{}  
\end{equation*} 
\begin{align*}
G_D = &  \sum_{k_D < i + j, h_A \leq i \leq k_A, h_C \leq j \leq  k_C} |a_i| |c_j| \e_D^{i + j -  k_D - \delta_D} \makebox[4mm]{}  \nonumber \\
& + G_A \sum_{h_C \leq j \leq k_C} |c_j| \e_D^{j + k_A + \delta_A  - k_D  - \delta_D}   \nonumber \\
& + G_C  \sum_{h_A \leq i \leq k_A} |a_i|  \e_D^{i + k_C  + \delta_C - k_D  - \delta_D}  \nonumber \\
& +  G_A  G_C \e_D^{k_A + k_C  + \delta_A + \delta_C - k_D  - \delta_D},   
\end{align*}
\begin{equation}\label{hoputr}
\e_D = \e_A \wedge \e_C, \makebox[61mm]{} 
\end{equation}
where coefficients $c_{r}, r = h_C, \ldots, k_C$ and parameters $h_C, k_C, \delta_C, G_C, \e_C$ are given for the  $(h_C, k_C, \delta_C$, $G_C, \e_C)$-expansion  of function 
$C(\e) = 1/B(\e)$  in  proposition {\bf (iv)}, or by formulas,
\begin{equation*}
\delta_D = \left\{
\begin{array}{ll}
\delta_A & \ \text{if} \ k_D = k_A - h_B < h_A + k_B  - 2 h_B,  \\
\delta_A \wedge \delta_B & \ \text{if} \ k_D = k_A - h_B = h_A + k_B  - 2 h_B,  \makebox[20mm]{}   \\
\delta_B & \ \text{if} \ k_D = h_A + k_B  - 2 h_B < k_A - h_B ,
\end{array}
\right. 
\end{equation*} 
\begin{equation*}
\geq \delta_A \wedge \delta_B,   \makebox[87mm]{}  
\end{equation*} 
\begin{align*}
G_D = & \, (\frac{|b_{h_B}|}{2})^{-1} \big( \sum_{k_A \wedge (h_A + k_B - h_B) < i \leq k_A} |a_i| \e_D^{i  - h_B - k_D  - \delta_D} \nonumber \\
& +  \sum_{k_A \wedge (h_A + k_B - h_B) < i + j, h_A \leq i \leq k_A, h_D \leq j \leq  k_D} |a_i| |d_j|  
\e_D^{i + j  - k_D - h_B - \delta_D}  \makebox[8mm]{}  \nonumber \\
\ \\
& + G_A \e_D^{k_A + \delta_A  - h_B - k_D - \delta_D}   + G_B \sum_{h_D \leq j \leq k_D} |d_j| \e_D^{j + k_B + \delta_B - h_B  - k_D - \delta_D} \big),
\end{align*}
\begin{equation}\label{hoputrk}
\e_D = \e_A \wedge  \e_B \wedge \left\{
\begin{array}{ll}
 \frac{|b_{h_B}|}{2} ( \sum_{h_B < i \leq k_B} |b_i| \e_B^{i - h_B - 1} \vspace{2mm} \\
+  \ G_B \e_B^{k_B + \delta_B - h_B - 1})^{-1}  & \ \text{if} \ h_B < k_B,   \vspace{2mm} \makebox[6mm]{}  \\
( \frac{|b_{h_B}|}{2G_B})^{\frac{1}{\delta_B}}  & \ \text{if} \ h_B = k_B.
\end{array}
\right.
\end{equation}
  
\end{itemize} 
} 

In what follows,  the following two lemmas, which present recurrent operational rules for computing coefficients and remainders for multiple summations and multiplications of Laurent asymptotic  expansions, will also be used. These lemmas are direct corollaries of Lemmas 3 and 4.

Let  $A_m(\e) =  a_{h_{A_{m}}, m}\e^{h_{A_{m}}} +
\cdots + a_{k_{A_{m}}, m}\e^{k_{A_{m}}} + o (\e^{k_{A_{m}}}), \e \in (0, \e_0]$ be a $(h_{A_m}, k_{A_m})$-expansion, for  $m = 1, \ldots, N$, \, 
$B_n(\e) = A_1(\e) + \cdots + A_n(\e), \e \in (0, \e_0]$, and $C_n(\e) = A_1(\e) \times \cdots \times A_n(\e), \e \in (0, \e_0]$, for $n = 1, \ldots, N$.  

The following two lemmas follow, respectively,  from Lemmas 3 and 4 and recurrent relations $B_n(\e) = B_{n-1}(\e) + A_n(\e), \e \in (0, \e_0], n = 2, \ldots, N$ and  
$C_n(\e) = C_{n-1}(\e) \cdot A_n(\e), \e \in (0, \e_0], n = 2, \ldots, N$, which hold for any $N \geq 2$.

\vspace{1mm}  

{\bf Lemma 5.} {\em The above asymptotic expansions
have the following operational rules for computing coefficients{\rm :}

\begin{itemize}
\item[{\rm \bf (i)}] If $A_m(\e), \e \in (0, \e_0]$ is a $(h_{A_m}, k_{A_m})$-expansion for $m = 1, \ldots, N$ where $N \geq 2$, then $B_n(\e) = b_{h_{B_n}, n} \e^{h_{B_n}} + \cdots + b_{k_{B_n}, n }\e^{k_{B_n}} + o (\e^{k_{B_n}}), \e \in (0, \e_0]$ is  a $(h_{B_n}, k_{B_n})$-expansion, for $n = 1, \ldots, N$, with $h_{B_1} = h_{A_1}, k_{B_1} = k_{A_1}$ and    $h_{B_n} = \min(h_{A_1}, \ldots$, $h_{A_n}) = h_{B_{n-1}} \wedge  h_{A_n}, k_{B_n} =  \min(k_{A_1}$, $\ldots, k_{A_n})  = k_{B_{n-1}} 
\wedge  k_{A_n}, n = 2, \ldots, N$ and the coefficients given by formulas
 $b_{h_{B_1} + l, 1} = a_{h_{A_1} + l, 1},  l = 0, \ldots, k_{B_1}  - h_{B_1} = k_{A_1}  - h_{A_1}$ and, for $l = 0, \ldots, k_{B_n}  - h_{B_n}, n = 2, \ldots, N$, by formulas, 
\begin{align}\label{joty}
b_{h_{B_n} + l, n} =  a_{h_{B_{n}} + l, 1} + \cdots + a_{h_{B_{n}} + l, n},
\end{align}
or
\begin{align} 
b_{h_{B_n} + l, n} =  b_{h_{B_{n-1}} + l, n-1} +  a_{h_{B_{n}} + l, n}, 
\end{align}
where  $b_{h_{B_{n-1}} + l, n-1} = 0, l = 0, \ldots$, $h_{B_{n-1}} - h_{B_n}$ and 
$a_{h_{B_{n}} + l, m} = 0, l = 0, \ldots h_{A_m} - h_{B_n},  m = 1, \ldots, n$. 

Expansions $B_n(\e), n = 1, \ldots, N$ are pivotal if and only if $b_{h_{B_n}, n}  = a_{h_{A_{1}}, 1} + \cdots +  a_{h_{A_{n}}, n} \neq 0,  n = 1, \ldots, N$. 

\item[{\rm \bf (ii)}] If $A_m(\e), \e \in (0, \e_0]$ is a $(h_{A_m}, k_{A_m})$-expansion for $m = 1, \ldots, N$ where $N \geq 2$, then $C_n(\e) = c_{h_{C_n}, n} 
\e^{h_{C_n}} + \cdots + c_{k_{C_n}, n }\e^{k_{C_n}} + o (\e^{k_{C_n}}), \e \in (0, \e_0]$ is  a $(h_{C_n}, k_{C_n})$-expansion, for $n = 1, \ldots, N$, with
$h_{C_1} = h_{A_1}, k_{C_1} = k_{A_1}$ and  $h_{C_n} =  h_{A_1} + \cdots + h_{A_n} = h_{C_{n-1}} +  h_{A_n}, k_{C_n} = \min(k_{A_l} + \sum_{1 \leq r \leq n, r \neq l} h_{A_r}, l = 1, \ldots, n) = (h_{C_{n-1}} + k_{A_n}) \wedge (k_{C_{n-1}} + h_{A_n}), n = 2, \ldots, N$  and 
coefficients given by formulas,
$c_{h_{C_1} + l, 1} = a_{h_{A_1} + l, 1}, l = 0, \ldots, k_{C_1}  - h_{C_1} = k_{A_1}  - h_{A_1}$ and,  for $l = 0, \ldots, k_{C_n}  - h_{C_n}, 
n = 2, \ldots, N$, by formulas,
\begin{align}\label{jotyk}
c_{h_{C_n} + l, n} =  \sum_{l_1 + \cdots +l_n = l, 0 \leq l_i \leq k_{A_i}  - h_{A_i}, i = 1, \ldots, n}  \, \prod_{1 \leq i \leq n} a_{h_{A_{i}} + l_i, i}, 
\end{align}
or
\begin{align}
c_{h_{C_n} + l, n} =  \sum_{0 \leq l' \leq l} c_{h_{C_{n-1}} + l', n-1}  a_{h_{A_{n}} + l - l', n}.
\end{align}

Expansions $C_n(\e), n = 1, \ldots, N$ are pivotal if and only if $c_{h_{C_n}, n}  = a_{h_{A_{1}}, 1} \times \cdots  \times a_{h_{A_{n}}, n} \neq 0, n = 1, \ldots, N$.

\item[{\rm \bf (iii)}] Asymptotic expansions for functions $B_n(\e) = A_1(\e) + \cdots + A_n(\e), n = 1, \ldots, N$ and $C_n(\e) = A_1(\e) \times \cdots \times A_n(\e), n = 1, \ldots, N$ are invariant with respect to any permutation, respectively, of summation and multiplication order in the above formulas {\rm (\ref{joty})} and {\rm (\ref{jotyk})}. 
\end{itemize} 
}

{\bf Lemma 6.} {\em The above asymptotic expansions
have the following operational rules for computing remainders{\rm :}

\begin{itemize}
\item[{\rm \bf (i)}] If $A_m(\e), \e \in (0, \e_0]$ is a $(h_{A_m},  k_{A_m}, \delta_{A_m}, G_{A_m}, \e_{A_m})$-expansion for $m = 1, \ldots, N$ where $N \geq 2$, then $B_n(\e), \e \in (0, \e_0]$ is  a $(h_{B_n}, k_{B_n}, \delta_{B_n}, G_{b_n}$, $\e_{B_n})$-expansion, for $n = 1, \ldots, N$, with parameters $h_{B_1} = h_{A_1}, k_{B_1} = k_{A_1}$ and    $h_{B_n} = \min(h_{A_1}, \ldots, h_{A_n}) = h_{B_{n-1}} \wedge  h_{A_n}, k_{B_n}  = \min(k_{A_1}$, $\ldots, k_{A_n}) = k_{B_{n-1}} 
\wedge  k_{A_n}, n = 2, \ldots, N$, coefficients
$b_{h_{B_n} + l, n}, l = 0, \ldots,   k_{B_n}$ $- h_{B_n}, n = 1, \ldots, N$ given in proposition {\bf (i)} of Lemma 5 and parameters $G_{B_n}, \delta_{B_n}$, $\e_{B_n}, n = 1, \ldots, N$ given by formulas $\delta_{B_1} = \delta_{A_1} \geq \delta^*_N = \min_{ 1 \leq m \leq n} \delta_{A_m}, G_{B_1} = G_{A_1}, \e_{B_1} = \e_{A_1}$ and, for $n = 2, \ldots, N$, by formulas, 
\begin{equation*}
\delta_{B_n} =    \min_{m \in {\mathbb K}_n} \delta_{A_m} \geq \delta^*_N,  \makebox[74mm]{}  
\end{equation*} 
where 
$$
{\mathbb K}_n = \{m:  1 \leq m \leq n, \,  k_m = \min(k_1, \ldots, k_n) \}, \makebox[30mm]{}  
$$
\begin{align*}
G_{B_n} = &  \sum_{1 \leq i \leq n} \big( G_{A_i} \e_{B_{n}}^{k_{A_i} + \delta_{A_i} - k_{B_{n}}  - \delta_{B_{n}} }  
+  \sum_{k_{B_{n}}  < j \leq k_{A_i}} |a_{A_{i},j}| \e_{B_{n}}^{j - k_{B_{n}} - \delta_{B_{n}}} \big), 
\end{align*}
\begin{equation}\label{dvadat}
\e_{B_n} = \min(\e_{A_1}, \ldots, \e_{A_n}),  \makebox[61mm]{}  
\end{equation} 
or by alternative recurrent formulas,
\begin{equation*}
\delta_{B_n} =    \min_{m \in {\mathbb K}_n} \delta_{A_m} = \left\{
\begin{array}{ll}
\delta_{B_{n-1}} & \ \text{if} \  k_{B_{n}} = k_{B_{n-1}}  < k_{A_n}, \\
\delta_{B_{n-1}}  \wedge \delta_{A_n} & \ \text{if} \  k_{B_{n}} = k_{B_{n-1}}  = k_{A_n},  \makebox[12mm]{}  \\
\delta_{A_n}  & \ \text{if} \ k_{B_{n}} =  k_{A_n} < k_{B_{n-1}},
\end{array}
\right. 
\end{equation*}
\begin{equation*}
\geq \delta^*_N,  \makebox[91mm]{}  
\end{equation*} 
\begin{align*}
G_{B_n} = & \, G_{B_{n-1}} \e_{B_n}^{k_{B_{n - 1}} + \delta_{B_{n-1}} - k_{B_{n}} - \delta_{B_{n}}}  + \sum_{k_{B_{n}} < i \leq k_{B_{n-1}}} |b_{B_{n-1}, i}| \e_{B_{n}}^{i - k_{B_{n}} - \delta_{B_{n}}}  \nonumber \\
&  + \ G_{A_n} \e_{B_{n}}^{k_{A_n} + \delta_{A_n} - k_{B_{n}}  - \delta_{B_{n}} }  +  \sum_{k_{B_{n}}  < j \leq k_{A_n}} |a_{A_{n},j}| \e_{B_{n}}^{j - k_{B_{n}} - \delta_{B_{n}}}, 
\end{align*}
\vspace{1mm}
\begin{equation}\label{dvad}
\e_{B_n}  = \e_{B_{n-1}} \wedge \e_{A_n},  \makebox[77mm]{}  
\end{equation} 

\item[{\rm \bf (ii)}] If $A_m(\e), \e \in (0, \e_0]$ is a $(h_{A_m}, k_{A_m}, \delta_{A_m}, G_{A_m}, \e_{A_m})$-expansion for $m = 1, \ldots, N$ where $N \geq 2$, then $C_n(\e), \e \in (0, \e_0]$ is  a $(h_{C_n}, k_{C_n},  \delta_{C_n}, G_{C_n}$, $\e_{C_N})$-expansion, for $n = 1, \ldots, N$, with parameters
$h_{C_1} = h_{A_1}, k_{C_1} = k_{A_1}$ and  $h_{C_n} = h_{C_{n-1}} +  h_{A_n} =  h_{A_1} + \cdots + h_{A_n}, k_{C_n} = (h_{C_{n-1}} + k_{A_n}) \wedge (k_{C_{n-1}} + h_{A_n}) =  
\min_{1 \leq l \leq n}(k_{A_l} + \sum_{1 \leq r \leq n, r \neq l} h_{A_r}), n = 2, \ldots, N$, coefficients
$c_{h_{C_n} + l, n}, l = 0, \ldots,   k_{C_n} - h_{C_n}, n = 1, \ldots, N$ given in proposition {\bf (ii)} of Lemma 5 and parameters $\delta_{C_n}, G_{C_n}$, $\e_{C_n}, n = 1, \ldots, N$ given by formulas $\delta_{C_1} = \delta_{A_1} \geq \delta^*_N = \min_{ 1 \leq m \leq n} \delta_{A_m}, G_{C_1} = G_{A_1}, \e_{C_1} = \e_{A_1}$ and, for $n = 2, \ldots, N$, by formulas, 
\begin{equation*}
\delta_{C_n} = \min_{m \in {\mathbb L}_n} \delta_{A_m} \geq  \delta^*_N, \makebox[73mm]{}  
\end{equation*} 
where
\begin{align*}
{\mathbb L}_n & = \{m:  1 \leq m \leq n, \,  (k_{A_m} + \sum_{1 \leq r \leq n, r \neq m} h_{A_r})  \nonumber \\
& = \min_{1 \leq l \leq n} (k_{A_l} + \sum_{1 \leq r \leq n, r \neq l} h_{A_r}) \}, \makebox[53mm]{}  
\end{align*}
\begin{align*}
G_{C_n} = & \sum_{k_{C_n} < l_1 + \cdots + l_n, h_{A_i} \leq l_i \leq k_{A_{i}}, i = 1, \ldots, n} 
\, \prod_{1 \leq i \leq n} |a_{A_{i}, l_i}| \e_{C_n}^{l_1 + \cdots + l_n -  k_{C_n} - \delta_{C_n}}  \nonumber \\
& + \sum_{1 \leq j \leq n} \, \prod_{1 \leq i \leq n, i \neq j}  \, \big( \sum_{h_{A_i} \leq l \leq k_{A_i}} |a_{A_{i}, l}| \e_{C_{n}}^{l}  \nonumber \\ & + G_{A_i} \e_{C_{n}}^{k_{A_{i}} + \delta_{A_{i}}} \big) G_{A_{j}}  \e_{C_{n}}^{k_{A_j} + \delta_{A_j} - k_{C_n}  - \delta_{C_n}},  
\end{align*}
\begin{equation}\label{dvadatasa}
\e_{C_n} = \min_{1 \leq i \leq n} \e_{A_{i}}. \makebox[77mm]{} 
\end{equation} 
or by alternative recurrent formulas,
\begin{equation*}
\delta_{C_n} = \left\{
\begin{array}{ll}
\delta_{C_{n-1}} & \ \text{if} \ k_{C_n} = h_{A_n}+ k_{C_{n-1}} < h_{C_{n-1}} + k_{A_n}, \\
\delta_{A_n} \wedge \delta_{C_{n-1}} & \ \text{if} \ k_{C_n} = h_B+ k_A = h_A + k_B,  \makebox[21mm]{}   \\
\delta_{A_n} & \ \text{if} \ k_{C_n} = h_{C_{n-1}}+ k_{A_n} < h_{A_n} + k_{C_{n-1}},
\end{array}
\right. 
\end{equation*}
\begin{equation*}
\geq \delta^*_N,  \makebox[91mm]{}  
\end{equation*} 
\begin{align*}
G_{C_n} = & \sum_{k_{C_n} < i + j, h_{C_{n-1}} \leq i \leq k_{C_{n-1}}, h_{A_n} \leq j \leq  k_{A_n}} |c_{C_{n-1}, i}| |a_{A_{n}, j}| \e_{C_n}^{i + j -  k_{C_n} - \delta_{C_n}} \makebox[5mm]{} \nonumber \\
& + G_{C_{n-1}} \sum_{h_{A_n} \leq j \leq k_{A_n}} |a_{A_{n}, j}| \e_{C_{n}}^{j + k_{C_{n-1}} + \delta_{C_{n-1}}   - k_{C_n}  - \delta_{C_n}}   \nonumber \\
& + G_{A_N}  \sum_{h_{C_{n-1}} \leq i \leq k_{C_{n-1}} } |c_{C_{n-1}, i}|  \e_{C_{n}}^{i + k_{A_n}  + \delta_{A_n} - k_{C_n}  - \delta_{C_n}} \nonumber \\
& + G_{A_N} G_{C_{n-1}} \e_{C_{n}}^{k_{A_n}  +k_{C_{n - 1}} + \delta_{A_n} + \delta_{C_{n - 1}}  - k_{C_n}  - \delta_{C_n}},  
\end{align*}
\begin{equation}\label{doyt}
\e_{C_n} = \e_{C_{n-1}} \wedge \e_{A_n}. \makebox[77mm]{} 
\end{equation} 

\item[{\rm \bf (iii)}] Parameters $\delta_{C_n}, G_{C_n}, \e_{C_n}, n = 1, \ldots, N$ in upper bounds for remainders in the asymptotic expansions for functions $B_n(\e) = A_1(\e) + \cdots + A_n(\e), n = 1, \ldots, N$ and $C_n(\e) = A_1(\e) \times \cdots \times A_n(\e), n = 1, \ldots, N$ are invariant with respect to any permutation, respectively, of summation and multiplication order in the above formulas {\rm (\ref{dvadat})} and 
{\rm (\ref{dvadatasa})}. 
\end{itemize} 
}

It should be noted that formulas (\ref{dvadat}) and (\ref{dvadatasa})  give, in general, the values, which are less or equal than the values for these constants given in alternative formulas, respectively,  (\ref{dvad}) and (\ref{doyt}). \vspace{1mm}

{\bf 2.3. Proofs of Lemmas 1--6}. The formulas given in Lemmas 1 and 2 are quite obvious. The same relate to formulas  and in propositions {\bf (i) -- (ii)} (the multiplication by a constant and summation rules) of Lemmas 3 and 4. They can be obtained by simple accumulation of  coefficients for different powers of $\e$ and terms accumulated in the corresponding remainders, as well obvious upper bounds for absolute values of sums of terms accumulated in the corresponding remainders.  Lemmas 5 and 6 are corollaries of Lemmas 3 and 4.

Let us, therefore, give short proofs of propositions {\bf (iii)} -- {\bf (v)} of Lemmas 3 and 4. 

Multiplication of asymptotic expansions $A(\e)$ and  $B(\e)$ penetrating proposition {\bf (iii)} of Lemmas 3 and accumulation of coefficients  for powers  $\e^l$ for $l = h_C, \ldots, k_C$ yields the following relation, 
{\small
\begin{align}\label{boply}
C(\e) & = A(\e) B(\e) \nonumber \\
& = (a_{h_A}\e^{h_A} +
\cdots + a_{k_A}\e^{k_A} + o_A(\e^{k_A})) (b_{h_B}\e^{h_B} +
\cdots + b_{k_B}\e^{k_B} + o_B(\e^{k_B}) ) \nonumber \\
&  = \sum_{h_C \leq l \leq k_C} \, \sum_{i + j = l, h_A \leq i \leq k_A, h_B \leq j \leq k_B } a_i b_j \e^l \nonumber \\
& \ \ + \sum_{k_C < i + j, h_A \leq i \leq k_A, h_B \leq j \leq k_B } a_i b_j \e^{i + j} \nonumber \\
&  \ \ + \sum_{h_B \leq j \leq k_B}  b_j \e^{j} o_A(\e^{k_A}) +  \sum_{h_A \leq i \leq k_A}  a_i \e^{i} o_B(\e^{k_B}) + o_A(\e^{k_A}) o_B(\e^{k_B}) \nonumber \\
&  = \sum_{h_C \leq l \leq k_C} c_{l} \e^l + o_C(\e^{k_C}),
\end{align}}
where
\begin{align}\label{hopret}
 o_C(\e^{k_C}) & = \sum_{k_C < i + j, h_A \leq i \leq k_A, h_B \leq j \leq k_B } a_i b_j \e^{i + j}  + \sum_{h_B \leq j \leq k_B}  b_j \e^{j} o_A(\e^{k_A}) \nonumber \\
& \ \ +  \sum_{h_A \leq i \leq k_A}  a_i \e^{i} o_B(\e^{k_B}) + o_A(\e^{k_A}) o_B(\e^{k_B}). 
\end{align}  

Obviously,
\begin{equation}
\frac{o_C(\e^{k_C})}{\e^{ k_C} }  \to 0 \ {\rm as} \ 0 < \e \to 0.
\end{equation}  

It should be noted that the accumulation of coefficients for powers $\e^l$ can be made in (\ref{boply}) only up to the maximal value $l = k_C = (h_A + k_B) \wedge (h_B + k_A)$,  because of  the presence in the  expression for the remainder $o_C(\e^{k_C})$ terms  $b_{h_B} \e^{h_B} o_A(\e^{k_A})$ and 
$a_{h_A} \e^{h_A} o_B(\e^{k_B})$.

Also, relation  (\ref{hopret}) readily implies relation (\ref{profa}), which  determines parameters $\delta_C, G_C,  \e_C$ in proposition {\bf (iii)} of Lemma 4.

The  assumptions of proposition  {\bf (iv)} in  Lemma 3 imply that the following relation holds,
\begin{equation}
\e^{-h_B} B(\e) \to  b_{h_B} \neq 0 \ {\rm as} \ 0 < \e \to 0.
\end{equation}  

This relation implies that there exists $0 < \e'_0 \leq \e_0$ such that $B(\e) \neq 0$ for $\e \in (0, \e'_0]$, and, 
thus, function $C(\e) =  \frac{1}{B(\e)}$ is well defined for $\e \in (0, \e'_0]$.

The  assumptions of proposition  {\bf (iv)} of Lemmas 3 also imply that,  
\begin{align}
\e^{h_B} C(\e) & = \frac{1}{b_{h_B} + b_{h_B +1 } \e
\cdots + b_{k_B}\e^{k_B - h_B} + o_B(\e^{k_B}) \e^{- h_B}}  \nonumber \\ 
& \to \frac{1}{b_{h_B}}  =  c_{h_C} \ {\rm as} \ 0 < \e \to 0.
\end{align}

This relation means that function  $\e^{h_B} C(\e)$ can be represented in the form  $\e^{h_B} C(\e) = c_{h_C} + o(1)$, where $c_{h_C} = b_{h_B}^{-1}$,  or, equivalently, that 
the following representation holds, 
\begin{equation}\label{bopate}
C(\e) =  c_{h_C} \e^{- h_B} + C_1(\e), \e \in (0, \e'_0],
\end{equation}  
where
\begin{equation}\label{boparate}
\frac{C_1(\e)}{\e^{- h_B}} \to 0  \ {\rm as} \ 0 < \e \to 0.
\end{equation} 

Relations (\ref{bopate}) and (\ref{boparate}) prove proposition  {\bf (iv)} of Lemmas 3 for the case, where $h _B = k_B$
that is equivalent to the relation $h_C = - h_B = k_C = k_B - 2 h_B$.

Note that, in the case $h _B = k_B$, the asymptotic expansion (\ref{bopate}) for function $C(\e)$ can not  be extended. Indeed,
\begin{align}\label{utrew}
\e^{ h_B  - 1} C_1(\e)  & = \e^{h_B -1} ( C(\e) -  c_{h_C} \e^{- h_B}) \nonumber \\
& = - \frac{c_{h_C}}{b_{h_B} + o_B(\e^{h_B})\e^{-h_B}} \frac{o_B(\e^{h_B})\e^{- h_B}}{ \e}
\end{align}

The term  $\frac{o_B(\e^{h_B})\e^{- h_B}}{ \e}$ on the right hand side in (\ref{utrew})  has  an uncertain asymptotic behaviour as $0 < \e \to 0$.

Let us now assume that $h _B + 1 = k_B$  that is 
equivalent to the relation $h_C  = - h_B  = k_C -1  = k_B - 2 h_B -1$.

In this case, the  assumptions of proposition  {\bf (iv)} of Lemma 3 and relations (\ref{bopate}),  (\ref{boparate}) and (\ref{utrew}) imply that
\begin{align}
\e^{h_B  -1} C_1(\e)  & = \e^{ h_B -1} ( C(\e) -  c_{h_C} \e^{- h_B})  \nonumber \\
&  = \frac{- b_{h_B +1} c_{h_C} - o_B(\e^{h_B +1}) \e^{- h_B -1} c_{h_C} }{b_{h_B} + b_{h_B + 1} \e + o_B(\e^{h_B+ 1}) \e^{- h_B}} \nonumber \\
& \to  \frac{- b_{h_B +1} c_{h_C}}{b_{h_B}}  = c_{h_C +1} \ {\rm as} \ 0 < \e \to 0.  
\end{align}  

This relation means that function  $\e^{h_B  -1} C_1(\e)$ can be represented in the form  $\e^{h_B  -1} C_1(\e) = c_{h_C +1} + o(1)$, where $c_{h_C +1} =
 b_{h_B}^{-1} b_{h_B +1} c_{h_C}$,   or, equivalently, that 
the following representation holds, 
\begin{equation}\label{boparko}
C(\e) =  c_{h_C} \e^{- h_B} + c_{h_C +1} \e^{- h_B +1} + C_2(\e), \e \in (0, \e'_0],
\end{equation}  
where
\begin{equation}\label{boparako}
\frac{C_2(\e)}{ \e^{- h_B +1}} \to 0  \ {\rm as} \ 0 < \e \to 0.
\end{equation} 

Relations (\ref{boparko}) and (\ref{boparako}) yields proposition  {\bf (iv)} of Lemmas 3 for the case, where 
$h _B + 1 = k_B$.

Note that, in the case $h _B + 1 = k_B$, the asymptotic expansion (\ref{boparko}) for function $C(\e)$ can not  be extended. Indeed, 
\begin{align}\label{bjdsa}
\e^{h_B  - 2} C_2(\e)  & = \e^{ h_B -2} ( C(\e) -  c_{h_C} \e^{- h_B} -  c_{h_C +1} \e^{- h_B +1})  \nonumber \\
&  = - \frac{c_{h_C}}{b_{h_B} + b_{h_B + 1} \e + o_B(\e^{h_B+ 1}) \e^{- h_B}}
\frac{o_B(\e^{h_B +1}) \e^{- h_B -1} }{\e}. 
\end{align}  

The term  $\frac{o_B(\e^{h_B +1})\e^{- h_B -1}}{ \e}$ on the right hand side in (\ref{bjdsa}) has  an uncertain asymptotic behaviour as $0 < \e \to 0$.

Repeating the above arguments, we can prove that function $C(\e)$ can be represented in the form of $(h_C, k_C)$-expansion, with parameters $h_C, k_C$ and coefficients $c_{h_C}, \ldots, c_{k_C}$ given in proposition {\bf (iii)} of Lemma 3, for the  general case, where  $h _B + n = k_B$, or, equivalently, $h_C  = - h_B  = k_C - n  = k_B - 2 h_B - n$, for any $n = 0, 1, \ldots$.

The $(h_C, k_C)$-expansion for function $C(\e) = \frac{1}{B(\e)}$ can be rewritten in the equivalent form of the following relation, 
\begin{equation}\label{nopty}
1 = (b_{h_B}e^{h_B} + \cdots +  b_{k_B} \e^{k_B} + o_B(\e^{k_B}))(c_{h_C} + \cdots +  c_{h_C} \e^{k_C} + o_C(\e^{k_C})),
\end{equation}

Proposition {\bf (iii)} of Lemma 3,  applied to the product on the right hand side in (\ref{nopty}), permits to represent this product in the form of $(h, k)$-expansion 
with parameters $h =  h_B + h_C = h_B - h_B = 0$ and $k = (h_B + k_C) \wedge (k_B +  h_C) = (k_B - h_B) \wedge (k_B - 2h_B + h_B) = k_B - h_B$.

By canceling coefficient for $\e^l$ on the left and right hand sides in  (\ref{nopty}), for $l = 0, \ldots, k_B - h_B$ and then solving equation  (\ref{nopty}) with respect to  the remainder $o_C(\e^{k_C})$ permits to find the following formula for this remainder, 
\begin{align}
o_C(\e^{k_C}) &  = - \frac{  \sum_{k_B - h_B < i + j, h_B \leq i \leq k_B, h_C \leq j \leq k_C} b_i c_j \e^{i + j} + \sum_{h_C \leq j \leq k_C}  c_j \e^j o_B(\e^{k_B}) }{b_{h_B}e^{h_B} + \cdots +  b_{k_B} \e^{k_B} + o_B(\e^{k_B})}  \nonumber \\
& =  - \frac{  \sum_{k_B - h_B < i + j, h_B \leq i \leq k_B, h_C \leq j \leq k_C} b_i c_j \e^{i + j - h_B}}{b_{h_B} + \cdots +  b_{k_B} \e^{k_B - h_B} + o_B(\e^{k_B}) \e^{- h_B}} \nonumber \\
& \quad   - \frac{\sum_{h_C \leq j \leq k_C}  c_j \e^{j - h_B} o_B(\e^{k_B})}{b_{h_B} + \cdots +  b_{k_B} \e^{k_B - h_B} + o_B(\e^{k_B}) \e^{- h_B}}.
\end{align}

The assumptions made in proposition {\bf (iv)} of Lemma 4, imply that $B(\e) \neq 0$ and the following inequality holds for $0 < \e \leq \e_C$, where $\e_C$  is given in relation  (\ref{naswe}),
\begin{equation}\label{berto}
|b_{h_B} + \cdots +  b_{k_B} \e^{k_B - h_B} + o_B(\e^{k_B}) \e^{- h_B}| \geq \frac{|b_{h_B|}}{2} > 0, 
\end{equation} 

The assumptions made in proposition {\bf (iv)} of Lemma 4 and inequality (\ref{berto}) finally  imply that the following inequality holds, for $0 < \e \leq \e_C$, 
\begin{align}
|o_C(\e^{k_C})| & \leq  \e^{k_B - 2h_B + \delta_B}  (\frac{|b_{h_B}|}{2})^{-1}  \nonumber \\ 
& \quad  \times \big( \sum_{k_B - h_B < i + j, h_B \leq i \leq k_B, h_C \leq j \leq  k_C} |b_i| |c_j| \e_C^{i + j - k_B + h_B - \delta_B} \nonumber \\
& \quad \quad  + G_B \sum_{h_C \leq j \leq k_C} |c_j| \e_C^{j + h_B} \big).  
\end{align}

This inequality proofs the proposition {\bf (iv)} of Lemma 4.

The first statement of proposition {\bf (v)} in Lemma 3 states that   function $D(\e)$ can be represented as $(h_D, k_D)$-expansion with parameters $h_D, k_D$ and coefficients $d_{h_d}, \ldots, d_{k_D}$ given in this proposition and relation (\ref{first}). It is  the direct corollary of propositions {\bf (iii)} and   {\bf (iv)} of Lemma 3, which, just,  should be  applied to the product 
$D(\e) = A(\e) \cdot \frac{1}{B(\e)}, \e \in (0, \e'_0]$.

Note that, in this case, parameters $h_D = h_A  + h_C  = h_A - h_B$ and $k_D = (k_A + h_C) \wedge (h_A + k_C) = (k_A - h_B) \wedge (h_A + k_B - 2 h_B)$.  

Now, when it is already proved that $D(\e)$ is $(h_D, k_D)$-expansion,  its coefficients can be also computed  by equalising coefficients for 
for powers $\e^l$ for $l = h_D, \ldots, k_D$ on the left and right hand sides of relation, 
\begin{align}\label{bertoma}
A(\e)  & = B(\e) D(\e) \nonumber \\ 
& = (b_{h_B}e^{h_B} + \cdots +  b_{h_B} \e^{k_B} + o_B(\e^{k_B})) \nonumber \\ 
& \quad \times ( d_{h_D}e^{h_D} + \cdots +  d_{h_D} \e^{k_D} + o_D(\e^{k_D})). 
\end{align} 

This  procedure yields the second statement of  proposition {\bf (v)} in Lemma 3 and the corresponding formulas given in relation (\ref{firsta}). 

The first statement of proposition proposition {\bf (v)} in Lemma 4 and relations (\ref{hoputr}) can be obtained by direct application of 
propositions {\bf (iii)} and  {\bf (iv)} and relations (\ref{profa}) and (\ref{naswe}) given in  Lemma 4, to the product $D(\e) = A(\e) \cdot \frac{1}{B(\e)}$.

Proposition {\bf (iii)} of Lemma 3,  applied to the product on the right hand side in  (\ref{bertoma}), permits to represent this product in the form of $(h, k)$-expansion 
with parameters $h =  h_B + h_D = h_B + h_A - h_B = h_A$ and $k = (h_B + k_D) \wedge (k_B +  h_D) = (h_B + (k_A - h_B)  \wedge (h_A + k_B - 2 h_B)) \wedge
(k_B + h_A - h_B) = k_A \wedge (k_B + h_A - h_B)$.

By canceling coefficient for $\e^l$ on the left and right hand sides in  (\ref{bertoma}), for $l = h_A, \ldots, k_A \wedge (k_B + h_A - h_B)$ and then solving equation  (\ref{bertoma}) with respect to  the remainder $o_D(\e^{k_D})$ yields the following formula for this remainder, 
\begin{align}
o_D(\e^{k_D}) &  =  \frac{  \sum_{k_A \wedge (k_B + h_A - h_B) < l \leq k_A} a_l \e^{l} + o_A(\e^{k_A}) }{b_{h_B}e^{h_B} + \cdots +  b_{k_B} \e^{k_B} + o_B(\e^{k_B})}  \nonumber \\
&  \quad  - \frac{ \sum_{k_A \wedge (k_B + h_A - h_B) < i + j, h_B \leq i \leq k_B, h_D \leq j \leq k_D} b_i d_j \e^{i + j} }{b_{h_B}e^{h_B} + \cdots +  b_{k_B} \e^{k_B} + o_B(\e^{k_B})}  \nonumber \\
&  \quad  - \frac{ \sum_{h_D \leq j \leq k_D}  d_j \e^j o_B(\e^{k_B}) }{b_{h_B}e^{h_B} + \cdots +  b_{k_B} \e^{k_B} + o_B(\e^{k_B})} \nonumber \\
&  =  \frac{  \sum_{k_A \wedge (k_B + h_A - h_B) < l \leq k_A} a_l \e^{l - h_B} + o_A(\e^{k_A}) \e^{- h_B}}{b_{h_B} + \cdots +  b_{k_B} \e^{k_B - h_B} + o_B(\e^{k_B})\e^{- h_B} }  \nonumber \\
&  \quad  - \frac{ \sum_{k_A \wedge (k_B + h_A - h_B) < i + j, h_B \leq i \leq k_B, h_D \leq j \leq k_D} b_i d_j \e^{i + j - h_B} }{b_{h_B} + \cdots +  b_{k_B} \e^{k_B - h_B} + o_B(\e^{k_B})\e^{- h_B} }  \nonumber \\
&  \quad   - \frac{ \sum_{h_D \leq j \leq k_D}  d_j \e^{j - h_B}o_B(\e^{k_B})  }{b_{h_B} + \cdots +  b_{k_B} \e^{k_B - h_B} + o_B(\e^{k_B})\e^{- h_B}}
\end{align}

The assumptions made in proposition {\bf (v)} of Lemma 4 and inequality (\ref{berto}) finally  imply that the following inequality holds, for $0 < \e \leq \e_D$  given in relation (\ref{hoputrk}), 
\begin{align}
|o_D(\e^{k_D})| & \leq  \e^{k_D  + \delta_D}  (\frac{|b_{h_B}|}{2})^{-1}  \nonumber \\ 
& \quad  \times \big( \sum_{k_A \wedge (k_B + h_A - h_B) < l \leq  k_A} |a_l| \e_D^{l - k_D   - h_B - \delta_D} \nonumber \\
& \quad  \quad +  \sum_{k_A \wedge (k_B + h_A - h_B) < i + j, h_B \leq i \leq k_B, h_D \leq j \leq  k_D} |b_i| |d_j| \e_D^{i + j - k_D - h_B - \delta_D} \nonumber \\
& \quad \quad  + G_A e_D^{k_A + \delta_A - h_B - k_D  - \delta_D} \nonumber \\
& \quad \quad  +  G_B \sum_{h_D \leq j \leq k_D} |d_j| \e_D^{j +k_B + \delta_B  - h_B - k_D  - \delta_D} \big).  
\end{align}

{\bf 2.4. Algebraic and related properties of operational rules for Laurent asymptotic expansions}.
Let us also introduce parameter $w_A = k_A - h_A$, which is a length 
of a  Laurent asymptotic expansion $A(\e) = a_{h}\e^{h_A} + \cdots + a_k\e^{k_A} + o_A(\e^{k_A})$.

The following useful  lemma takes place. \vspace{1mm}

{\bf Lemma 7}. {\em The following relations hold for Laurent asymptotic expansions penetrating 
Lemma 3{\rm :}
\begin{itemize}

\item[{\rm \bf (i)}] If $C(\e) = c A(\e)$, then  $w_C = w_A$.

\item[{\rm \bf (ii)}] If $C(\e) = A(\e) + B(\e)$, then  $w_A \wedge w_B \leq w_C \leq w_A \vee w_B$.

\item[{\rm \bf (iii)}]  If $C(\e) = A(\e) \cdot B(\e)$, then  $w_C = w_A \wedge w_B$.

\item[{\rm \bf (iv)}]  If $C(\e) = 1/B(\e)$, then  $w_C = w_B$.

\item[{\rm \bf (v)}]  If $D(\e) = A(\e)/B(\e)$, then  $w_D = w_A \wedge w_B$.

\item[{\rm \bf (vi)}] If $w_A = w_B = w$ then $w_C  = w_D  = w$ for all  Laurent asymptotic expansions penetrating Lemma 3.

\end{itemize}
}

The proof of this simple lemma readily follows from formulas for parameters $h$ and $k$ penetrating  propositions {\bf (i)} -- {\bf (v)} of Lemma 3. 

Let us again consider four Laurent asymptotic  expansions, $A(\e) = a_{h_A}\e^{h_A} +
\cdots + a_{k_A}\e^{k_A} + o_A(\e^{k_A})$, $B(\e) = b_{h_B}\e^{h_B} +
\cdots + b_{k_B}\e^{k_B} + o_B(\e^{k_B})$, $C(\e) = c_{h_C}\e^{h_C} +
\cdots + c_{k_C}\e^{k_C} + o_C(\e^{k_C})$, and $D(\e) =
d_{h_D}\e^{h_D} + \cdots + d_{k_D}\e^{k_D} + o_D(\e^{k_D})$ defined for $0 < \e \leq \e_0$, for some $0 < \e_0 \leq 1$. 

Below,  sums $\sum_{l = h}^k d_l$ are counted as $0$ if $k < h$.

The following lemma is also a corollary of Lemma 3. \vspace{1mm}

{\bf Lemma 8}.  {\em The summation and multiplication operations for the Laurent asymptotic expansions penetrating propositions {\bf (ii)} and {\bf (iii)} in Lemma 3 possess the following algebraic properties, which should be understood as identities for the corresponding  asymptotic expansions{\rm :}
\begin{itemize}
\item[{\rm \bf (i)}] The summation operation is  commutative, i.e., $C(\e) = A(\e) + B(\e) = B(\e) + A(\e)$, where 
$h_C = h_{A + B} = h_{B + A} =  h_A \wedge h_B, \, k_C = k_{A + B} =  k_{B + A} =  k_A \wedge k_B$, and,
{\small 
\begin{equation}
C(\e) = \sum_{l = 0}^{k_C - h_C} (a_{h_C + l} + b_{h_C + l}) \e^{h_C + l} + o_C(\e^{k_c}),
\end{equation}} 
where $a_{h_C + l} = 0$ for $0 \leq  l < h_A - h_C $, $b_{h_C + l} = 0$ for $0 \leq l < h_B - h_C$. 
\item[{\rm \bf (ii)}] The summation operation is  associative, i.e., $D(\e) = (A(\e) + B(\e)) + C(\e) = A(\e) + (B(\e) +  C(\e)) = A(\e) + B(\e) +  C(\e)$, where
$h_D  = h_{(A+ B) + C} = h_{A + (B + C)} = h_{A + B + C} = h_A \wedge h_B \wedge h_C, \, k_D =   k_{(A+ B) + C} = k_{A + (B + C)}  = k_{A + B + C} = k_A \wedge k_B \wedge k_C$, and,
{\small 
\begin{equation}
D(\e) = \sum_{l = 0}^{k_D - h_D} (a_{h_D + l} + b_{h_D + l}  + c_{h_D + l}) \e^{h_D + l} + o_D(\e^{k_D}), 
\end{equation}} 
where $a_{h_D + l}  = 0$  for $0 \leq l < h_{A} - h_D$, $b_{h_D + l}  = 0$ for $0 \leq l < h_{B} - h_D$, $c_{h_D + l}  = 0$ for $0 \leq l < h_{C} - h_D$. 
\item[{\rm \bf (iii)}] The multiplication operation is  commutative, i.e., $C(\e) = A(\e) \cdot B(\e) = B(\e) \cdot A(\e)$, where 
$h_C = h_{A \cdot  B} = h_{B \cdot  A} =  h_A + h_B, \,  k_C = k_{A \cdot  B} =  k_{B \cdot  A} =  (h_A + k_B) \wedge(k_A +  h_B)$, and, 
{\small 
\begin{equation}
C(\e) = \sum_{l = 0}^{k_C - h_C} \big( \sum_{l_1 + l_2 = l, l_1, l_2 \geq 0} a_{h_A + l_1} b_{h_B + l_2} \big) \e^{h_C + l} + o_C(\e^{k_c}). 
\end{equation}} 
\item[{\rm \bf (iv)}] The multiplication operation is  associative, i.e., $D(\e) = (A(\e) \cdot B(\e)) \cdot  C(\e) = A(\e) \cdot  (B(\e) \cdot   C(\e)) = A(\e) \cdot  B(\e) \cdot   C(\e)$, where $h_D  = h_{(A \cdot  B) \cdot  C} = h_{A \cdot  (B \cdot  C)} =  h_{A \cdot  B \cdot  C} =  h_A +  h_B  + h_C, \, k_D =   k_{(A \cdot  B) \cdot  C} = k_{A \cdot  (B \cdot  C)} =  k_{A \cdot  B \cdot  C} = (h_A + h_B + k_C) \wedge (h_A + k_B + h_C) \wedge (k_A + h_B + h_C)$, and,  
{\small
\begin{equation}
D(\e) = \sum_{l = 0}^{k_D - h_D} \big( \sum_{l_1 + l_2 + l_3 = l, l_1, l_2, l_3 \geq 0} a_{h_A + l_1} b_{h_B + l_2} c_{h_C + l_3} \big) \e^{h_D + l} + o_D(\e^{k_D}).
\end{equation}}
\item[{\rm \bf (v)}]  The summation and multiplication operations possess distributive property, i.e.,  $D(\e)  = 
(A(\e) + B(\e)) \cdot C(\e) = A(\e)\cdot C(\e) + B(\e) \cdot C(\e)$, where $h_D = h_{(A + B) \cdot C} = h_{A \cdot C + B \cdot C} = h_A \wedge h_B + h_C = (h_A + h_C) \wedge (h_B + h_C), \, k_D = k_{(A + B) \cdot C} = k_{A \cdot C + B \cdot C} = (h_A \wedge h_B + k_C) \wedge  (k_A \wedge k_B + h_C) = (h_A + k_C) \wedge (k_A + h_C) \wedge (h_B + k_C) \wedge (k_B + h_C)$, and,
{\small
\begin{align}
D(\e) & = \sum_{l = 0}^{k_D - h_D} \big( \sum_{l_1 + l_2 = l, l_1, l_2 \geq 0} (a_{h_A \wedge h_B + l_1} \nonumber \\ 
& \quad  + b_{h_A \wedge h_B + l_1}) c_{h_C + l_2} \big) \e^{h_D + l} + o_D(\e^{k_D}) \nonumber \\  
& = \sum_{l = 0}^{k_D - h_A - h_C}  \big( \sum_{l_1 + l_2 = l, l_1, l_2 \geq 0} a_{h_A + l_1}c_{h_C + l_2} \big) \e^{h_A + h_C +l}  \nonumber \\
& \quad + \sum_{l = 0}^{k_D - h_B - h_C}  \big( \sum_{l_1 + l_2 = l, l_1, l_2 \geq 0} b_{h_B + l_1}c_{h_C + l_2} \big) \e^{h_B + h_C +l} + o_D(\e^{k_D}).
\end{align}}
where $a_{h_A \wedge h_B + l} = 0$. for $0 \leq  l < h_A - h_A \wedge h_B$, $b_{h_A \wedge h_B  + l} = 0$,  for $0 \leq l < h_B - h_A \wedge h_B$. 
\end{itemize}
}

The summation and multiplication rules for computing of upper bounds remainders penetrating propositions  {\bf (ii)} and {\bf (iii)} in Lemma 4 possess the  communicative property. This follows from formulas {\rm (\ref{dvadat})} and 
{\rm (\ref{dvadatasa})} given in Lemma 6. 

However,  the summation and multiplication rules for computing of upper bounds for remainders presented in propositions  {\bf (ii)} and {\bf (iii)} of Lemma 4  do not possess associative and distributional properties. The question about the form of upper bounds for the corresponding remainders, which would possess these properties, remains open. 

As follows from  Lemma 4, operational rules presented in this lemma possess special property that let one give an effective low bounds for parameter $\delta_A$ for any $(h_{A}, k_{A}, \delta_{A}, G_{A}, \e_{A})$-expansion $A(\e)$ obtained as the result of a  finite sequence of operations (multiplication by a constant,  summation, multiplication, and division) performed with  with $(h_{A_i}, k_{A_i},  \delta_{A_i}, G_{A_i}, \e_{A_i})$-expansions  $A_i(\e), i = 1, \ldots, N$ from some finite set of such expansions.

The following lemma takes place. \vspace{1mm}

{\bf Lemma 9}.  {\em The operational rules for computing remainders of asymptotic expansions with explicit upper bounds for remainders presented in 
propositions {\bf (ii)} and {\bf (iii)} of Lemma 4 possess the following properties{\rm :}
\begin{itemize}

\item[{\rm \bf (i)}] If $C(\e) = A(\e) +  B(\e) = B(\e) +  A(\e)$ then   $\delta_C = \delta_{A + B} = \delta_{B + A}, G_C = G_{A + B} = 
G_{B + A}$ and  $\e_C = \e_{A + B}  = \e_{B + A}$, where parameters $\delta_C, G_C$ and $\e_C$ are given by formula {\rm (\ref{profac})} in proposition  {\bf (ii)} of Lemma 4. 

\item[{\rm \bf (ii)}] If $C(\e) = A(\e) \cdot  B(\e) = B(\e) \cdot  A(\e)$ then   $\delta_C = \delta_{A \cdot B} = \delta_{B \cdot A}, G_C = G_{A \cdot B} = G_{B \cdot A}$  and $\e_C = \e_{A \cdot B} = \e_{B \cdot A}$, where parameters $\delta_C, G_C$ and $\e_C$ are given by formula {\rm (\ref{profa})} in proposition  {\bf (iii)} of Lemma 4. 

\item[{\rm \bf (iii)}]  If  $A(\e)$  is $(h_{A}, k_{A}, \delta_{A}, G_{A}, \e_{A})$-expansion obtained as the result of a finite sequence of operations {\rm (}multiplication by a constant,  summation, multiplication, and quotient{\rm )} performed with $(h_{A_i}, k_{A_i}, \delta_{A_i}, G_{A_i}, \e_{A_i})$-expansi- ons  $A_i(\e), i = 1, \ldots, N$ from some finite set of such expansions, then $\delta_A \geq \delta^*_N = \min_{1 \leq i \leq N} \delta_{A_i}$ that makes it possible to rewrite  $A(\e)$ as the 
$(h_{A}, k_{A}, \delta^*_N, G^*_{A, N}, \e_{A})$-expansion, with parameter $G^*_{A, N} = G_{A} \e_A^{\delta_A - \delta^*_N}$.

\end{itemize}

}

\vspace{1mm}

{\bf 3. Perturbed Markov chains and semi-Markov processes} \\

In this section, we  formulate basic perturbation conditions for Markov chains and semi-Markov processes. \vspace{1mm}

{\bf 3.1. Perturbed Markov chains}.  
Let  ${\mathbb X} = \{1, \ldots, N \}$ and $\eta^{(\e)}_n, n = 0, 1, \ldots$ be, for every $\e \in (0, \e_0]$, a homogeneous Markov chain with a phase space ${\mathbb X}$, an initial distribution $p^{(\e)}_i = \PP \{\eta^{(\e)}_0 = i \}, i \in {\mathbb X}$ and transition probabilities defined for $i, j \in \XX$,
\begin{equation}\label{semi}
p_{ij}(\e) = \PP \{ \eta^{(\e)}_{n+1} = j / \eta^{(\e)}_{n} = i \}. 
\end{equation}

Let us assume that the following condition holds:

\ \\

\begin{itemize}
\item [${\bf A}$:] There exist sets $\YY_i \subseteq \XX, i \in \XX$ and $\e_0 \in (0, 1]$ such that: {\bf (a)} probabilities $p_{ij}(\e)  > 0, j \in \YY_i, i \in \XX$, for   $\e \in (0, \e_0]$; {\bf (b)} probabilities $p_{ij}(\e)  = 0, j \in  \overline{\YY}_i, i \in \XX$, for   $\e \in (0, \e_0]$;
{ \bf (c)} there exist $n_{ij} \geq 1$ and a chain of states $i = l_{ij, 0}, l_{ij, 1}, \ldots, l_{ij, n_{ij}} = j$ such that  $l_{ij, 1} \in \YY_{l_{ij, 0}}, \ldots, l_{ij, n_{ij}} \in 
 \YY_{l_{ij, n_{ij}-1}}$,  for every pair of states $i, j \in \XX$. 
\end{itemize}

We refer to sets  $\YY_i, i \in \XX$ as transition sets.

Conditions ${\bf A}$ implies that all sets $\YY_i \neq \emptyset, i \in \XX$, since  matrix $\| p_{ij}(\e) \|$ is stochastic, for every $\e \in (0, \e_0]$.

We now assume that the following perturbation condition holds:
\begin{itemize}
\item [${\bf B}$:] $p_{ij}(\e) =  \sum_{l = l_{ij}^-}^{l_{ij}^+} a_{ij}[l]\e^l + o_{ij}(\e^{l_{ij}^+})$, where $a_{ij}[ l_{ij}^-] > 0$ and $0 \leq l_{ij}^- \leq l_{ij}^+ < \infty$, for $j \in \YY_i, i \in \XX$,  and $o_{ij}(\e^{l_{ij}^+})/\e^{l_{ij}^+} \to 0$ as $\e \to 0$, for $j \in \YY_i, i \in \XX$.  
\end{itemize}

Some additional  conditions should be imposed on parameters $\e_0 \in ( 0, 1]$ and  $l_{ij}^{\pm}, j \in \YY_i, i \in \XX$, and coefficients $a_{ij}[l], l = l_{ij}^-, \ldots, l_{ij}^+,  j \in \YY_i, i \in \XX$, penetrating  the  asymptotic expansions 
condition ${\bf B}$, in order this condition would be consistent with the model assumption that matrix $\| p_{ij}(\e) \|$ is stochastic, for every $\e \in (0, \e_0]$, and with condition ${\bf A}$.

Condition ${\bf B}$ implies that there exits $\e_0 \in (0, 1]$ such that the following relation holds,
\begin{equation}\label{setabi}
p_{ij}(\e)   = \sum_{l = l_{ij}^-}^{l_{ij}^+} a_{ij}[l]\e^l + o_{ij}(\e^{l_{ij}^+})  > 0, \ j \in \YY_i, \  i \in \XX, \ \e \in (0, \e_0].
\end{equation}

Thus, condition ${\bf B}$ is consistent with condition  ${\bf A}$ {\bf (a)}.

The model assumption that matrix $\| p_{ij}(\e) \|$ is stochastic is, under conditions ${\bf A}$ {\bf (a)}  and {\bf (b)}, equivalent to the following relation, 
\begin{equation}\label{stocha}
 \sum_{j \in  \YY_{i}}  p_{ij}(\e) = 1, \ j \in \YY_i, i   \in \XX, \  \e \in (0, \e_0].
\end{equation}

Condition ${\bf B}$  and proposition {\bf (i)} (the multiple summation rule) of Lemma 5 imply that sum $\sum_{j \in \ZZ} p_{ij}(\e)$
can, for every subset $\ZZ \subseteq \YY_i$ and $i \in \XX$, be represented in the form of the following asymptotic 
expansion, 
\begin{equation}\label{expako}
\sum_{j \in \ZZ} p_{ij}(\e)   = \sum_{l = \check{l}_{i, \ZZ}^-}^{l_{i, \ZZ}^+} a_{i, \ZZ}[l]\e^l + o_{i, \ZZ}(\e^{l_{i, \ZZ}^+}),
\end{equation}
where
\begin{equation}\label{parame}
l_{i, \ZZ}^{-} = \min_{j \in \ZZ} l_{ij}^{-}, \ l_{i, \ZZ}^{+} = \min_{j \in \ZZ} l_{ij}^{+},
\end{equation} 
\begin{equation}\label{paramase}
a_{i, \ZZ}[l] = \sum_{j \in \ZZ} a_{ij}[l], \ l =  l_{i, \ZZ}^-, \ldots, l_{i, \ZZ}^+, 
\end{equation} 
where $a_{ij}[l] = 0$, for $0 \leq l < l_{ij}^-, j \in \ZZ$, and 
\begin{equation}
o_{i, \ZZ}(\e^{l_{i, \ZZ}^+}) = \sum_{j \in \ZZ} ( \sum_{l_{i, \ZZ}^+ 
<  l \leq l_{ij}^+} a_{ij}[l]\e^l + o_{ij}(\e^{l_{ij}^+}) ).
\end{equation}

In terms of asymptotic expansions, constant $1$ can be represented, for every $n = 0, 1, \ldots$, in the form of the following pivotal $(0, n)$-expansion, 
\begin{equation}\label{vopit}
1 = 1 + 0 \e + \cdots 0 \e^n + o_{n}(\e^n), 
\end{equation}
where remainders $o_{n}(\e^n) \equiv 0, n = 0, 1, \ldots$. 

Moreover, the above expansion is a $(0, n, 1, G, \e_0)$-expansion for any $0 < G < \infty$ and $n = 0, 1, \ldots$.  

Relation (\ref{stocha}) permits us  apply Lemma 2   to the asymptotic expansions  given in relations 
(\ref{expako}) and ({\ref{vopit}). Not that, in this case,  $l_{i, \YY_i}^- = 0$, otherwise the expression on the right hand side in  
(\ref{expako}) would converge to zero as $\e \to 0$. 
Let us take $n = l_{i,  \YY_i}^+$ in relation (\ref{vopit}). In this case $h_1  = 0$ and $k_1 = l_{i, \YY_i}^+$ in the asymptotic expansion given in relation (\ref{vopit}). 

Lemma 2  and the model stochasticity assumption  (\ref{stocha}) imply that, in this case, the following condition should hold for the coefficients of asymptotic expansions penetrating condition ${\bf B}$: 
\begin{itemize}
\item [${\bf C}$:] {\bf (a)} $a_{i, \YY_i}[l] =  \sum_{j \in \YY_i } a_{ij}[l] = {\rm I}(l = 0), \
0 \leq l \leq l_{i, \YY_i}^+, \ i \in \XX$, where $a_{ij}[l] = 0$, for $0 \leq l < l_{ij}^-, j \in \YY_i, i \in \XX$;  {\bf (b)} $o_{i, \YY_i}(\e^{l_{i, \YY_i}^+}) \equiv 
o_{l_{i,  \YY_i}^+}(\e^{l_{i,  \YY_i}^+}) \equiv 0,  i \in {\mathbb X}$. 
\end{itemize}

{\bf Remark 1}. It is possible to prove that conditions  ${\bf A}$ -- ${\bf C}$ and the model stochasticity assumption (\ref{stocha})  imply that the asymptotic expansion in (\ref{expako}) satisfy, for every $\ZZ \subseteq \YY_i$ and $i \in \XX$,  one of 
the following additional conditions: (a) $l_{i, \ZZ}^- > 0$; (b) $l_{i, \ZZ}^- = 0, a_{i, \ZZ}[0] < 1$; (c)  $l_{i, \ZZ}^- = 0, a_{i, \ZZ}[0] = 1$ and there exists  $0 < l_{i, \ZZ} \leq l^+_{i, \ZZ}$ such that $a_{i, \ZZ}[l]  = 0,  0 <  l < l_{i, \ZZ}$, but $a_{i, \ZZ}[ l_{i, \ZZ}] < 0$; or (d) $l_{i, \ZZ}^- = 0, a_{i, \ZZ}[0] = 1$ and 
$a_{i, \ZZ}[l]  = 0, 0 <   l \leq l^+_{i, \ZZ}$, but the remainder $o_{i, \ZZ}(\e^{l_{i, \ZZ}^+})$ is a nonpositive function of $\e$.

The above proposition implies that there exists $\e_0 \in (0, 1]$ such that  $\sum_{l = \check{l}_{i, \ZZ}^-}^{l_{i, \ZZ}^+} a_{i, \ZZ}[l]\e^l + o_{i, \ZZ}(\e^{l_{i, \ZZ}^+}) \leq 1, \ZZ \subseteq \YY_i, i \in \XX, \e \in (0, \e_0]$. 
Thus, conditions ${\bf A}$ -- ${\bf C}$  are also consistent with the relations $\sum_{j \in \ZZ} p_{ij}(\e) \leq 1, \ZZ\subseteq \YY_i, i \in \XX, \e \in (0, \e_0]$, which follows from the model  stochasticity assumption (\ref{stocha}).

In the case, where the  asymptotic expansions penetrating condition ${\bf B}$ are supposed to be given in the form of asymptotic expansions with explicit upper bounds for remainders, we replace it by the following condition:
\begin{itemize}
\item [${\bf B'}$:] $p_{ij}(\e) =  \sum_{l = l_{ij}^-}^{l_{ij}^+} a_{ij}[l]\e^l + o_{ij}(\e^{l_{ij}^+})$, where $a_{ij}[ l_{ij}^-] > 0$ and $0 \leq l_{ij}^- \leq l_{ij}^+ < \infty$, for $j \in \YY_i, i \in \XX$, and $|o_{ij}(\e^{l_{ij}^+})| \leq G_{ij}\e^{l_{ij}^+ + \delta_{ij}}, 0 < \e \leq \e_{ij}$, for $j \in \YY_i, i \in \XX$, where $ 0 < \delta_{ij} \leq 1, 0 < G_{ij} < \infty$ and  $0 < \e_{ij} \leq \e_0$. 
\end{itemize}

{\bf 3.2. Perturbed semi-Markov processes}. Let  ${\mathbb X} = \{1, \ldots, N \}$ and $(\eta^{(\e)}_n, \kappa^{(\e)}_n), n = 0, 1, \ldots$ be, for every $\e \in (0, 1]$, a Markov renewal process, i.e., a homogeneous Markov chain with the phase space 
${\mathbb X} \times [0, \infty)$, an initial distribution $p^{(\e)}_i = \PP \{\eta^{(\e)}_0 = i, \kappa^{(\e)}_0 = 0 \} = \PP \{\eta^{(\e)}_0 = i \}, i \in {\mathbb X}$ and transition probabilities defined for $(i, s), (j, t) \in  {\mathbb X} \times [0, \infty)$, 
\begin{equation}\label{semika}
Q^{(\e)}_{ij}(t) = \PP \{ \eta^{(\e)}_{n+1} = j, \kappa^{(\e)}_{n+1} \leq t / \eta^{(\e)}_{n} = i, \kappa^{(\e)}_{n} = s \}.
\end{equation}

In this case,  the random sequence $\eta^{(\e)}_n, n = 0, 1, \ldots$ is also a homogeneous (embedded) Markov chain with the  phase space $\XX$ and transition probabilities defined for $i, j \in \XX$,
\begin{equation}\label{embed}
p_{ij}(\e) = \PP \{ \eta^{(\e)}_{n+1} = j / \eta^{(\e)}_{n} = i \} = Q^{(\e)}_{ij}(\infty). 
\end{equation}

We assume that condition ${\bf A}$ holds. This implies that Markov chain $\eta^{(\e)}_n$ has one class of communicative states, for every  $\e \in (0, \e_0]$.

We also assume that the following condition excluding instant transitions holds: 
\begin{itemize}
\item [${\bf D}$:]  $Q_{ij}^{(\e)}(0) = 0, \ i, j \in \XX$, for every $\e \in (0, \e_0]$.
\end{itemize}

Let us now introduce  a semi-Markov process,
\begin{equation}\label{sepr}
\eta^{(\e)}(t) = \eta^{(\e)}_{\nu^{(\e)}(t)}, \ t \geq 0.
\end{equation}
where
\begin{equation}\label{numj}
\nu^{(\e)}(t) = \max(n \geq 0: \zeta^{(\e)}_n \leq t), \ t \geq 0, 
\end{equation}
is a number of jumps in the  time interval  $[0, t]$, and 
\begin{equation}\label{mome}
\zeta^{(\e)}_n = \kappa^{(\e)}_1 + \cdots + \kappa^{(\e)}_n, \ n = 0, 1, \ldots, 
\end{equation}
are sequential moments of jumps for the semi-Markov process $\eta^{(\e)}(t)$.

If $Q^{(\e)}_{ij}(t) = (1 - e^{- \lambda_i(\e) t})p_{ij}(\e), t \geq 0, i, j \in \XX$, then $\eta^{(\e)}(t), t \geq 0$ is a continuous time homogeneous Markov chain. 

If  $Q^{(\e)}_{ij}(t) = I(t \geq 1)p_{ij}(\e), t \geq 0, i, j \in \XX$,  then $\eta^{(\e)}(t) = \eta^{(\e)}_{[t]}, t \geq 0$ is a discrete time homogeneous Markov chain embedded in continuous time.

Let us also introduce expectations of sojourn times,
\begin{equation}\label{expe}
e_{ij}(\e) = \EE_i   \kappa^{(\e)}_{1} I(  \eta^{(\e)}_{1} = j) = \int_0^\infty t Q^{(\e)}_{ij}(dt), \ i, j \in  \XX.
\end{equation}

We also assume that the following condition holds: 
\begin{itemize}
\item [${\bf E}$:] $e_{ij}(\e) < \infty, \ i, j \in \XX$, for   $\e \in (0, \e_0]$.
\end{itemize}

Here and henceforth,  notations $\PP_i $ and $ \EE_i $ are used for conditional probabilities and expectations under 
condition $\eta^{(\e)}_0 = i$.

In the case of continuous time Markov chain,  $e_{ij}(\e) = \frac{1}{\lambda_i(\e)}p_{ij}(\e)$,  $i, j \in \XX$.

In the case  of discrete time Markov chain, $e_{ij}(\e) = p_{ij}(\e)$,  $i, j \in \XX$. 

Conditions ${\bf A}$ and ${\bf D}$ imply that, for every  $\e \in (0, \e_0]$, expectations  $e_{ij}(\e)  > 0$, for $j \in \YY_i, i \in \XX$, and  $e_{ij}(\e)  = 0$, for $j \in  \overline{\YY}_i, i \in \XX$. 

We now assume  that the following perturbation condition holds:
\begin{itemize}
\item [${\bf F}$:]  $e_{ij}(\e) =  \sum_{l = m_{ij}^-}^{m_{ij}^+} b_{ij}[l]\e^l + \dot{o}_{ij}(\e^{m_{ij}^+})$, where $b_{ij}[ m_{ij}^-] > 0$ and $-\infty < m_{ij}^- \leq m_{ij}^+ < \infty$, for $j \in \YY_i, i \in \XX$
and  $ \dot{o}_{ij}(\e^{m_{ij}^+})/\e^{m_{ij}^+} \to 0$ as $\e \to 0$, for $j \in \YY_i, i \in \XX$. \vspace{1mm}
\end{itemize}

In particular, in the case of discrete time Markov chain, condition ${\bf B}$ implies condition ${\bf F}$ to hold, since, in this case, expectations $e_{ij}(\e) = p_{ij}(\e) , j \in \YY_i, i \in \XX$. 

Condition ${\bf F}$ implies that there exits $\e_0 \in (0, 1]$ such that the following relation holds, 
\begin{equation}\label{setabik}
e_{ij}(\e)  > 0, \ j \in \YY_i, \  i \in \XX, \ \e \in (0, \e_0]. 
\end{equation}

This is consistent with condition ${\bf D}$.

In the case, where the  asymptotic expansions penetration condition ${\bf F}$ are given in the form of asymptotic expansions with explicit upper bounds for  remainders, we assumer that the following condition holds:
\begin{itemize}
\item [${\bf F'}$:] $e_{ij}(\e) =  \sum_{l = m_{ij}^-}^{m_{ij}^+} b_{ij}[l]\e^l + \dot{o}_{ij}(\e^{m_{ij}^+})$, where $b_{ij}[ l_{ij}^-] > 0$ and $- \infty < \leq m_{ij}^- \leq m_{ij}^+ < \infty$, for $j \in \YY_i, i \in \XX$,  and $|\dot{o}_{ij}(\e^{l_{ij}^+})| \leq \dot{G}_{ij}\e^{m_{ij}^+ + \dot{\delta}_{ij}}, 0 < \e \leq \dot{\e}_{ij}$, for $j \in \YY_i, i \in \XX$, where $ 0 < \dot{\delta}_{ij} \leq 1, 0 < \dot{G}_{ij} < \infty$ and  $0 < \dot{\e}_{ij} \leq \e_0$. 
\end{itemize} 
\vspace{2mm}

It is also worse to note that the perturbation conditions ${\bf B}$ and ${\bf F}$ are independent.

To see this, let us take arbitrary functions $p_{ij}(\e), j \in {\mathbb Y}_i, i \in \XX$ and $e_{ij}(\e), j \in {\mathbb Y}_i, i \in \XX$ satisfying, respectively, conditions  ${\bf B}$ and ${\bf F}$, and, also,  relations (\ref{setabi}), (\ref{stocha}) and (\ref{setabik}). Then,   there exist semi-Markov transition probabilities $Q_{ij}^{(\e)}(t), t \geq 0, j \in {\mathbb Y}_i, i \in \XX$ such that $Q_{ij}^{(\e)}(\infty) = p_{ij}(\e), j \in {\mathbb Y}_i, i \in \XX$ and $\int_0^\infty t Q_{ij}^{(\e)}(dt) = e_{ij}(\e), j \in {\mathbb Y}_i, i \in \XX$, for every  $\e \in (0, \e_0]$.

It is readily seen that, for example,  semi-Markov transition probabilities $Q_{ij}^{(\e)}(t) = I(t \geq e_{ij}(\e)/p_{ij}(\e))p_{ij}(\e), t \geq 0, j \in {\mathbb Y}_i, i \in \XX$ satisfy the above relations. \\

{\bf 4. Stationary distributions for  semi-Markov processes} \\

In this section, we give basic formulas for stationary distributions  for semi-Markov processes, in particular, formulas connecting  stationary 
distributions with expectations of return times. \vspace{1mm}

{\bf 4.1. Stationary distributions for  semi-Markov processes}. Condition ${\bf A}$ guarantees that the phase space $\XX$ is one class of communicative states for Markov chain $\eta^{(\e)}_n$, for every $\e \in (0, \e_0]$, i.e., the Markov chain $\eta^{(\e)}_n$ is ergodic, and, thus,  for every   $\e \in (0, \e_0]$, there exist the unique stationary distribution $\bar{\rho}(\e) = \langle \rho_1(\e), \ldots$, $\rho_N(\e) \rangle$, which is given by the following ergodic relation,
\begin{equation}\label{stati}
\bar{\mu}^{(\e)}_{i,n} = \frac{1}{n} \sum_{k = 1}^n I(\eta^{(\e)}_k = i) 
\stackrel{a.s.}{\longrightarrow} \rho_i(\e)  \ {\rm as} \ n \to \infty, \ i \in \XX.
 \end{equation}
 
It is useful to note that the ergodic relation (\ref{stati}) holds for any initial distribution $\bar{p}^{(\e)} = \langle p^{(\e)}_1, \ldots p^{(\e)}_N \rangle$ and the stationary distribution $\bar{\rho}(\e)$  does not depend on the initial distribution.

As known, $\rho_i(\e), i \in \XX$ is 
the unique positive solution for the system of linear equations, 
\begin{equation}\label{statiop}
\left\{
\begin{array}{ll}
\rho_{j}(\e) = \sum_{i \in \XX} \rho_i(\e) p_{ij}(\e), j \in \XX, \vspace{2mm} \\
 \sum_{i \in \XX} \rho_i = 1.  
\end{array}
\right.
\end{equation}

Conditions ${\bf A}$, ${\bf D}$ and   ${\bf E}$ imply that, for every  $\e \in (0, \e_0]$,  the semi-Markov process $\eta^{(\e)}(t)$ is also ergodic and its stationary distribution $\bar{\pi}(\e) = \langle \pi_1(\e), \ldots$, $\pi_N(\e) \rangle$ is given by the following ergodic relation, for  $i \in \XX$,
\begin{equation}\label{statik}
\bar{\mu}^{(\e)}_i(t) = \frac{1}{t} \int_0^t I(\eta^{(\e)}(s) = i)ds \stackrel{a.s.}{\longrightarrow} \pi_i(\e)  \ {\rm as} \ t \to \infty.
 \end{equation}
 
As in  (\ref{stati}), the ergodic relation (\ref{statik}) holds for any initial distribution $\bar{p}^{(\e)}$ and the stationary distribution 
$\bar{\pi}(\e)$  does not depend on the initial distribution.

The stationary distributions for the semi-Markov process $\eta^{(\e)}(t)$ and the embedded Markov chain  $\eta^{(\e)}_n$ are connected by the following relation,
\begin{equation}\label{stad}
\pi_i(\e) = \frac{ \rho_i(\e) e_i(\e)}{\sum_{j \in \XX} \rho_j(\e) e_j(\e)}, \ i \in \XX,
\end{equation}
where
\begin{equation}\label{expec}
e_{i}(\e) = \EE_i  \kappa^{(\e)}_{1}  = \sum_{j \in \XX} e_{ij}(\e),  \ i \in   \XX.
\end{equation}

Condition  ${\bf B}$ implies that there exist limits, 
\begin{equation}\label{limp}
p_{i j}(0) = \lim_{\e \to 0} p_{ij}(\e) = \left\{
\begin{array}{cl}
a_{ij}[0] & \ \text{if} \ l_{ij}^-  = 0, j \in \YY_i, i \in \XX, \\
0 & \ \text{if} \ l_{ij}^-  > 0,  j \in \YY_i, i \in \XX,  \\
0 & \ \text{if} \ j \in \overline{\YY}_i, i \in \XX.
\end{array}
\right.
\end{equation}

Matrix $\| p_{ij}(\e) \|$ is stochastic, for every $\e \in (0, \e_0]$ and, thus, matrix  $\| p_{ij}(0) \|$ is also stochastic. 

However, it is possible that matrix $\| p_{ij}(0) \|$ has more zero elements than matrices $\| p_{ij}(\e) \|$. 

Therefore, a Markov chain $\eta^{(0)}_n, n = 0, 1, \ldots$, with the phase space $\XX$ and the matrix of transition probabilities 
$\| p_{ij}(0) \|$ can be not ergodic, and its phase space $\XX$ can consist of one or several closed classes of communicative states plus, possibly, a class of transient states. 

Condition ${\bf F}$ implies that there exist limits,
\begin{equation}\label{limm}
e_{i j}(0) = \lim_{\e \to 0} e_{ij}(\e) = \left \{
\begin{array}{cl}
\infty & \ \text{if} \ m_{ij}^-  < 0,  j \in \YY_i, i \in \XX,  \\
b_{ij}[0] & \ \text{if} \ m_{ij}^-  = 0, j \in \YY_i, i \in \XX, \\
0 & \ \text{if} \ m_{ij}^-  > 0,  j \in \YY_i, i \in \XX,  \\
0 & \ \text{if} \ j \in \overline{\YY}_i, i \in \XX.
\end{array}
\right.
\end{equation}

Out goal is to design an effective  algorithm for constructing asymptotic  expansions for stationary probabilities 
$\pi_i(\e), i \in \XX$,  under assumption that conditions ${\bf A}$ -- ${\bf F}$ hold. 

As we shall se, the proposed algorithm, based on a special techniques of sequential phase space reduction,  can be applied for models with asymptotically coupled and uncoupled phase spaces and all types of asymptotic behavior of  expected sojourn  times.

The models of continuous and  discrete Markov chains are particular cases.

In particular,  asymptotic expansions for stationary probabilities $\rho_i(\e), i \in \XX$ coincide with expansions for stationary probabilities $\pi_i(\e), i \in \XX$,  for the discrete time Markov chain, where expectations $e_{ij}(\e) = p_{ij}(\e), i, j \in \XX$. \vspace{1mm}

{\bf 4.2. Expected hitting times and  stationary probabilities for semi-Markov processes}. 
Let us define hitting times, which are random variables given by the following relation, for $j \in \XX$,
\begin{equation}\label{hitt}
\tau^{(\e)}_j = \sum_{n = 1}^{\nu^{(\e)}_j} \kappa^{(\e)}_n,
\end{equation}
where
\begin{equation}\label{hitta}
\nu^{(\e)}_j = \min(n \geq 1: \eta^{(\e)}_n = j).
\end{equation}

Let us denote,
\begin{equation}\label{hitaex}
E_{ij}(\e) = \EE_i \tau^{(\e)}_j, \ i,j \in \XX.
\end{equation}

As is known, conditions ${\bf A}$, ${\bf D}$ and ${\bf E}$ imply that, for every $\e \in (0, \e_0]$, 
\begin{equation}\label{hitaexa}
0 < E_{ij}(\e) < \infty, \ i,j \in \XX.
\end{equation}

Moreover, under the above conditions, the expectations $E_{ij}(\e), i \in \XX$ are, for every $j \in \XX$, the unique solution for the following system of linear equations,
\begin{equation}\label{systera}
\big\{ E_{ij}(\e) = e_i(\e) + \sum_{r \neq j} p_{ir}(\e) E_{rj}(\e), \ i \in \XX.
\end{equation}

The following relation plays an important role in what follows,
\begin{equation}\label{hittana}
\pi_i(\e) = \frac{e_i(\e)}{E_{ii}(\e)}, \ i \in \XX.
\end{equation}

In fact this formula is an alternative form of relation (\ref{stad}). Indeed, as is known, $E_{ii}(\e) = \sum_{j \in \XX} e_j(\e) f_{ii,j}(\e)$, where $f_{ii,j}(\e)$ is the expected number of visits by the Markov chain $\eta^{(\e)}_n$ the state $j$ between two sequential 
visits of the state $i$. As also known, $f_{ii,j}(\e) = \rho_j(\e)  / \rho_i(\e) , i, j \in \XX$.  

Formula (\ref{hittana}) permits reduce the problem of constructing asymptotic expansions for semi-Markov stationary probabilities 
$\pi_i(\e)$ to the problem of constructing Laurent asymptotic expansions for expectation of hitting times $E_{ii}(\e)$. \\

{\bf 5. Semi-Markov processes with reduced phase spaces} \\

In this section, we present a procedure for one-step procedure of phase space reduction for semi-Markov processes. \vspace{1mm}

{\bf 5.1. Reduction of  a phase space for semi-Markov process}. Let us choose some state $r$ and consider the reduced phase space $_r\XX = \{i \in \XX, i \neq r \}$,  with the state $r$ excluded from the phase space $\XX$.

Let us assume that $p^{(\e)}_r = \PP \{\eta^{(\e)}_0 = r \} = 0$ and  define  the sequential moments of hitting the reduced space $_r\XX$ by the embedded Markov chain $\eta^{(\e)}_n$,
\begin{equation}\label{sequen}
_r\xi^{(\e)}_n = \min(k > \, _r\xi^{(\e)}_{n-1}, \  \eta^{(\e)}_k \in \, _r\XX), \ n = 1, 2, \ldots, \  _r\xi^{(\e)}_0 = 0.
\end{equation}

Now, let us define the random sequence,
\begin{equation}\label{sequena}
(_r\eta^{(\e)}_n, \, _r\kappa^{(\e)}_n) = \left\{
\begin{array}{ll}
(\eta^{(\e)}_{_r\xi^{(\e)}_n} \,, 
\sum_{k = \, _r\xi^{(\e)}_{n-1} +1}^{_r\xi^{(\e)}_n} \kappa^{(\e)}_k) & \ \text{for}  \ n = 1, 2, \ldots, \vspace{2mm} \\
(\eta^{(\e)}_0, 0) & \ \text{for}  \ n = 0. \\
\end{array}
\right.
\end{equation}

This sequence is also a Markov renewal process with a phase space $_r\XX \times [0, \infty)$,   the initial distribution
$p^{(\e)}_i = \PP \{\eta^{(\e)}_0 = i \}, i \in \, _r\XX$ (remind that $p^{(\e)}_r = 0$), and transition probabilities defined for 
$(i, s), (j, t) \in \, _r\XX \times [0, \infty)$, 
\begin{equation}\label{semir}
_rQ^{(\e)}_{ij}(t) = \PP \{ \, _r\eta^{(\e)}_{n+1} = j, \, _r\kappa^{(\e)}_{n+1} \leq t / \, _r\eta^{(\e)}_{n} = i, \, _r\kappa^{(\e)}_{n} = s \}. 
\end{equation}

Respectively, one can define the transformed semi-Markov process with the  reduced phase space  $_r\XX$, 
\begin{equation}\label{sepra}
_r\eta^{(\e)}(t) = \, _r\eta^{(\e)}_{_r\nu^{(\e)}(t)}, \ t \geq 0.
\end{equation}
where
\begin{equation}\label{numja}
_r\nu^{(\e)}(t) = \max(n \geq 0: \, _r\zeta^{(\e)}_n \leq t), \ t \geq 0, 
\end{equation}
is a number of jumps at time interval  $[0, t]$, and 
\begin{equation}\label{momea}
_r\zeta^{(\e)}_n = \, _r\kappa^{(\e)}_1 + \cdots + \, _r\kappa^{(\e)}_n, \ n = 0, 1, \ldots, 
\end{equation}
are sequential moments of jumps for the semi-Markov process $_r\eta^{(\e)}(t)$. 

The transition probabilities $_rQ^{(\e)}_{ij}(t)$ are expressed via the transition probabilities $Q^{(\e)}_{ij}(t)$ by the 
following formula, for $i, j \in \, _r\XX, t \geq 0$, 
\begin{equation}\label{numjabo}
_rQ^{(\e)}_{ij}(t) = Q^{(\e)}_{ij}(t) + \sum_{n = 0}^\infty Q^{(\e)}_{ir}(t) * Q^{(\e) *n}_{rr}(t) * Q^{(\e)}_{rj}(t). 
\end{equation}

Here, symbol $*$ is used to denote a convolution of distribution functions (possibly improper) and 
$Q^{(\e) *n}_{rr}(t)$ is the $n$ times convolution of the distribution function  $Q^{(\e)}_{rr}(t)$ given by the following recurrent formula, for $r \in \XX$,
\begin{equation}\label{recu}
Q^{(\e) *n}_{rr}(t) =  \left \{
\begin{array}{cl}
\int_0^t  Q^{(\e) *(n -1)}_{rr}(t -s)  Q^{(\e)}_{rr}(ds)    & \ \text{for} \  t \geq 0 \ \text{and} \ n \geq 1, \\
I(t \geq 0)  & \ \text{for} \  t \geq 0 \ \text{and} \ n = 0. \\
\end{array}
\right.
\end{equation}

Relation (\ref{numjabo}) directly implies the following formula for transition probabilities of the embedded Markov chain $_r\eta^{(\e)}_n$, for $i, j \in \, _r\XX$,
\begin{align}\label{transit}
_rp_{ij}(\e) & = \, _rQ^{(\e)}_{ij}(\infty) \nonumber \\
& = p_{ij}(\e)  + \sum_{n = 0}^\infty p_{ir}(\e)  p_{rr}(\e)^n p_{rj}(\e) \nonumber \\
& = p_{ij}(\e) + p_{ir}(\e) \frac{p_{rj}(\e)}{1 - p_{rr}(\e)}.
\end{align}

Let us denote,
\begin{align}\label{denoja}
\YY_{i r}^+ =  \{ j \in  \YY_i: j \neq r  \}, \ i, r \in \XX. 
\end{align}
and
\begin{align}\label{deno}
\YY_{ir}^- =  \{ j \in \, _r\XX: r \in \YY_i, \ j \in \YY_r  \}, \ i \in \, _r\XX. 
\end{align}

Condition ${\bf A}$ implies that sets $\YY_{rr}^+ \neq \emptyset, r \in \XX$.

Thus, probabilities $1 - p_{rr}(\e) > 0, r \in \XX$, for every $\e \in (0, \e_0]$.

That is why,
\begin{align}\label{denojat}
_r\YY_{i} & =  \{ j \in \,  _r\XX:  \, _rp_{ij}(\e) > 0, \e \in (0, \e_0] \} 
\nonumber \\
&   = \{ j \in \,  _r\XX: j \in \YY_i \}  \cup \{ j \in \,  _r\XX: r \in \YY_i,  
j \in \YY_r \}, \nonumber \\
& = \YY_{ir}^+ \cup  \YY_{ir}^-, \  i \in \, _r\XX. 
\end{align}

Relation (\ref{transit}) and condition ${\bf A}$,  assumed to hold for the Markov chain  $\eta^{(\e)}_n$,  imply that  condition ${\bf A}$  also holds for the Markov chain $_r\eta^{(\e)}_n$, with the  sets $_r\YY_i , i \in \, _r\XX$.

Indeed, let $i \in \, _r\XX$. If $j \in  \YY_{ir}^+$ then $p_{ij}(\e) > 0$ and, thus, $_rp_{ij}(\e) > 0$. If $j \in  \YY_{ir}^-$ then $p_{ir}(\e), p_{rj}(\e) > 0$ and, again,  $_rp_{ij}(\e) > 0$. If  $j \notin \YY_{ir}^+ \cup  \YY_{ir}^-$ then $p_{ij}(\e) = 0$ and  $p_{ir}(\e), p_{rj}(\e) = 0$. By relation (\ref{transit}), this implies that $_rp_{ij}(\e) = 0$.

Let $i \in \, _r\XX$. If $\YY_{ir}^+ \neq \emptyset$ then $_r\YY_i \neq \emptyset$. If $\YY_{ir}^+ = \emptyset$ then $r \in \YY_i$ and, thus,  $p_{ir}(\e) > 0$. Then, $\YY_{ir}^- = \{ j \in \, _r\XX: \, p_{rj}(\e) > 0 \} = \YY_{rr}^+  \neq \emptyset$. Therefore,
sets $_r\YY_i \neq \emptyset, i \in \, _r\XX$.

Thus, conditions ${\bf A}$ {\bf (a)} and  {\bf (b)} assumed to hold for the Markov chain  $\eta^{(\e)}_n$,  imply  that  these conditions   also hold for the Markov chain $_r\eta^{(\e)}_n$, with sets $_r\YY_i, i \in \, _r\XX$ replacing sets $\YY_i, i \in \XX$.

Also, let $i, j \in  \, _r\XX$ and  $i = l_0, l_1, \ldots, l_{n_{ij}} = j$ be a chain of states such that 
$l_{1} \in \YY_{l_0}, \ldots, l_{n} \in  \YY_{l_{n_{ij}-1}}$. As was remarked above, we can always to assume that states $l_1, \ldots, l_{n_{ij}-1}$ are different and that $l_1, \ldots, l_{n_{ij}-1} \neq i, j$. This implies that either $l_1, \ldots, l_{n_{ij}-1} \neq r$ or there exist at most one $1 \leq k \leq n_{ij} -1$ such that $i_k = r$. In the first case, $l_{1} \in \, _r\YY_{l_0}, \ldots, l_{n_{ij}} \in  \, _r\YY_{l_{n_{ij}-1}}$. In the second case, $l_{1} \in \, _r\YY_{l_0}, \ldots, l_{k-1} \in  \, _r\YY_{l_{k-2}},  l_{k-1} \in  \, _r\YY_{l_{k+1}}, \ldots, l_{n_{ij}} \in  \, _r\YY_{l_{n_{ij}-1}}$. 

Thus, condition ${\bf A}$ {\bf (c)}  assumed to hold for the Markov chain  $\eta^{(\e)}_n$,  imply  that  this condition   also holds for the Markov chain $_r\eta^{(\e)}_n$.

Let us define distribution functions,
\begin{equation}\label{distr}
F^{(\e)}_{i}(t)  = \sum_{j \in \XX} Q^{(\e)}_{ij}(t), t \geq 0, \ i, j \in \XX.
\end{equation}
and 
\begin{equation}\label{distra}
F^{(\e)}_{ij}(t) =  \left \{
\begin{array}{cl}
Q^{(\e)}_{ij}(t)/ p_{ij}(\e)     & \ \text{for} \  t \geq 0 \ \text{if} \ p_{ij}(\e) > 0, \\
 F^{(\e)}_{i}(t) & \ \text{for} \  t \geq 0 \ \text{if} \  p_{ij}(\e) = 0. \\
\end{array}
\right.
\end{equation}

Obviously,
\begin{equation}\label{expa}
\tilde{e}_{ij}(\e) = \int_0^\infty t F^{(\e)}_{ij}(dt)  =  \left \{
\begin{array}{cl}
e_{ij}(\e)/ p_{ij}(\e)     &  \ \text{if} \ p_{ij}(\e) > 0, \\
e_i(\e) &  \ \text{if} \  p_{ij}(\e) = 0, 
\end{array}
\right.
\end{equation}
and 
\begin{equation}\label{dexp}
e_i(\e)  = \int_0^\infty t  F^{(\e)}_{i}(dt), \ i \in \XX.
\end{equation}

Also, let us introduce expectations,
\begin{equation}\label{dexpop}
_re_{ij}(\e)  = \int_0^\infty t  \, _rQ^{(\e)}_{ij}(dt), \ i, j \in \, _r\XX.
\end{equation}

Relation (\ref{numjabo}) directly implies the following formula for expectations  $_re_{ij}(\e)$, $i, j \in \, _r\XX$,
\begin{align}\label{expectaga}
_re_{ij}(\e) & =  \tilde{e}_{ij}(\e) p_{ij}(\e) + \sum_{n = 0}^\infty \big( \tilde{e}_{ir}(\e)  + n \tilde{e}_{rr}(\e)  
+ \tilde{e}_{rj}(\e) \big)  p_{ir}(\e)  p_{rr}(\e)^n p_{rj}(\e) \nonumber \\
& = e_{ij}(\e)  + e_{ir}(\e) \frac{p_{rj}(\e)}{1 - p_{rr}(\e)} \nonumber \\
& \quad  + \, e_{rr}(\e)  \frac{p_{ir}(\e)}{1 - p_{rr}(\e)} \frac{p_{rj}(\e) }{1 - p_{rr}(\e)} +  e_{rj}(\e)  \frac{p_{ir}(\e)}{1 - p_{rr}(\e)}. 
\end{align}

Relation (\ref{expectaga}) implies that conditions ${\bf D}$ and  ${\bf E}$, assumed to hold for the semi-Markov process   
$\eta^{(\e)}(t)$,  imply  that  these conditions   also hold for the semi-Markov process $_r\eta^{(\e)}(t)$. \vspace{1mm}

{\bf 5.2. Hitting times for reduced semi-Markov processes}.
The first hitting times to a state $j \neq r$ are connected for Markov chains  $\eta^{(\e)}_n$  and 
 $_r\eta^{(\e)}_n$ by the following relation,
\begin{align}\label{relan}
\nu^{(\e)}_j & = \min(n \geq 1: \eta^{(\e)}_n = j) \nonumber \\ 
& = \min(_r\xi^{(\e)}_n \geq 1: \, _r\eta^{(\e)}_n  = j) = \, _r\xi^{(\e)}_{_r\nu^{(\e)}_j },
\end{align} 
where
\begin{equation}\label{relana}
_r\nu^{(\e)}_j = \min(n \geq 1: \, _r\eta^{(\e)}_n = j). 
\end{equation}

Relations  (\ref{relan}) and  (\ref{relana}) imply that the following relation hold for the first hitting times to a state $j \neq r$ for the semi-Markov 
processes $\eta^{(\e)}(t)$ and $_r\eta^{(\e)}(t)$,
\begin{align}\label{relanak}
\tau^{(\e)}_j & = \sum_{n = 1}^{\nu^{(\e)}_j} \kappa^{(\e)}_n  = \sum_{n = 1}^{_r\xi^{(\e)}_{_r\nu^{(\e)}_j }} \kappa^{(\e)}_n \nonumber \\
& =  \sum_{n = 1}^{_r\nu^{(\e)}_j } \, _r\kappa^{(\e)}_n = \, _r\tau^{(\e)}_j.  
\end{align}

Let us summarize the above remarks in the following theorem, which play the key role in what follows. \vspace{1mm}

{\bf Theorem 1}. {\em Let conditions ${\bf A}$, ${\bf D}$ and  ${\bf E}$ hold and the initial distribution satisfies  the 
assumption, $p^{(\e)}_r = 0$, for every $\e \in(0, \e_0]$. Then, for any state $j \neq r$,  the first hitting times $\tau^{(\e)}_j$ and $_r\tau^{(\e)}_j$ to the state $j$, respectively, for semi-Markov processes $\eta^{(\e)}(t)$ and $_r\eta^{(\e)}(t)$, coincide.}  \\

{\bf 6. Asymptotic expansions for transition characteristics of perturbed semi-Markov processes with  reduced phase spaces} \\

In this section, we present algorithms for re-calculation of asymptotic expansions for perturbed semi-Markov processes with   reduced phase spaces.

\vspace{1mm}

{\bf 6.1 Asymptotic expansions for non-absorption probabilities}.  As was mentioned above, condition ${\bf A}$ implies that the non-absorption probability $\bar{p}_{ii}(\e) = 1 - p_{ii}(\e)  > 0, i \in \XX, \e \in (0, \e_0]$.

Let us introduce the set,
\begin{equation}
\YY = \{ i \in \XX: i \in \YY_i \}.
\end{equation}

{\bf Algorithm 1.} This is an algorithm for constructing  asymptotic expansions for non-absorption probabilities  $\bar{p}_{ii}(\e), i \in {\mathbb X}$. \vspace{1mm}

{\bf Case 1}: $i \in {\mathbb Y}$. 

Let us use the following relation, which holds, for every $i \in \YY$ and $\e \in  (0, \e_0]$, 
\begin{align}\label{expabas}
\bar{p}_{ii}(\e)  = 1 - p_{ii}(\e) = \sum_{j \in \YY^+_{ii}} p_{ij}(\e)
\end{align}

{\bf 1.1.} To construct the $(h'_i, k'_i)$-expansion for the non-absorption probability $\bar{p}_{ii}(\e) = 1 - p_{ii}(\e)$ by applying  the propositions {\bf (i)} (the multiplication by a constant rule) and  {\bf (ii)} (the summation rule) of Lemma 3 to the $(l^-_{ii}, l^+_{ii})$-expansion for transition probability $p_{ii}(\e)$ given in condition ${\bf B}$ (first, this expansion is multiplied by constant $-1$ and, second, is summated with constant $1$ represented as $(0, l^+_{ii})$-expansion given in relation (\ref{vopit})).  In this case, parameters $h'_i = 0, k'_i = l_{ii}^+$.  

{\bf 1.2.} To construct the $(h''_i, k''_i)$-expansion for the non-absorption probability $\bar{p}_{ii}(\e) = \sum_{j \in \YY^+_{ii}} p_{ij}(\e)$ using the corresponding asymptotic expansions for transition probabilities $p_{ij}(\e), j \in {\mathbb Y}^+_{ii}$ given in condition ${\bf B}$, and the proposition {\bf (i)} (the multiple summation rule) of Lemma 5. In this case, parameters $h''_i = \min_{j \in {\mathbb Y}^+_{ii}} l^-_{ij}, k''_i =  \min_{j \in {\mathbb Y}^+_{ii}} l^+_{ij}$. This asymptotic  expansion is pivotal. 

{\bf 1.3.} To construct the $(\bar{h}_i, \bar{k}_i)$-expansion for the non-absorption probability $\bar{p}_{ii}(\e)$ using relation (\ref{expabas}), and propositions {\bf (i)} -- {\bf (iv)} of  Lemma 1. In this case, parameters $\bar{h}_i = 0 \vee h''_i = h''_i, \bar{k}_i = k'_i \vee k''_i = l_{ii}^+ \vee \min_{j \in {\mathbb Y}^+_{ii}} l^+_{ij}$. This asymptotic expansion is pivotal. \vspace{1mm} 

It should be noted that $(l^-_{ii}, l^+_{ii})$-expansion for the transition  probability $p_{ii}(\e)$ given in condition ${\bf B}$ and $(h''_i, k''_i)$-expansion  for function $\bar{p}_{ii}(\e) = \sum_{j \in \YY^+_{ii}} p_{ij}(\e)$ given in Step 1.2,  satisfy, for every $i \in \XX$, additional conditions given in Remark 1, respectively, for set $\ZZ = \{ i \}$ and  set $\ZZ  = \YY^+_{ii}$.
\vspace{1mm}

{\bf Case 2}: $i \in \overline{{\mathbb Y}}$.  \vspace{1mm}

{\bf 1.4.}  In this case, the non-absorption probability $\bar{p}_{ii}(\e)  \equiv 1$. If necessary, it can be represented in the form of $(0, n)$-expansion given by relation
(\ref{vopit}), for any $n = 0, 1, \ldots$. \vspace{1mm}

The above remarks  can be summarized in the following theorem. \vspace{1mm}

{\bf Lemma  10}. {\em Let conditions ${\bf A}$, ${\bf B}$ and  ${\bf C}$ hold. Then, the asymptotic expansions for the non-absorption probabilities $\bar{p}_{ii}(\e) , i \in \XX$ are given in \mbox{Algorithm 1.}} \vspace{1mm}

{\bf Algorithm 2.} This is an algorithm for computing upper bounds for remainders of asymptotic expansions for 
non-absorption probabilities $\bar{p}_{ii}(\e)$, $i \in \XX$. 

\vspace{1mm}

{\bf Case 1}: $i \in {\mathbb Y}$. \vspace{1mm}

{\bf 2.1.} To construct $(h'_i, k'_i, \delta'_i, G'_i, \e'_i)$-expansion for the non-absorption probability $\bar{p}_{ii}(\e) = 1 - p_{ii}(\e)$ by applying the propositions {\bf (i)} (the multiplication by a constant rule) and  {\bf (ii} (the summation rule) of Lemma 4 to the $(l^-_{ii}, l^+_{ii}, \delta_{ii}, G_{ii}, \e_{ii})$-expansion for the transition probability $p_{ii}(\e)$ given in condition ${\bf B'}$ and  (first, this expansion  is multiplied by constant $-1$ and, second, is summated with constant $1$ represented as $(0, l^+_{ii}, 1, G, \e_0)$-expansion given in relation (\ref{vopit})), third, constant $G$ can be replaced by $0$, since it can be taken an arbitrary small. In this case, parameters $\delta'_i = \delta_{ii}, G'_i = G_{ii}, \e'_{i} = \e_{ii}$. 

{\bf 2.2.} To construct the $(h''_i, k''_i, \delta''_i, G''_i, \e''_i)$-expansion for the non-absorption probability $\bar{p}_{ii}(\e) = \sum_{j \in \YY^+_{ii}} p_{ij}(\e)$ using the  $(l^-_{ij}, l^+_{ij}, \delta_{ij}, G_{ij}, \e_{ij})$-expansions for transition probabilities $p_{ij}(\e), j \in {\mathbb Y}^+_{ii}$ given in condition ${\bf B'}$, and the proposition {\bf (i)} (the multiple summation rule) of Lemma 6. In this case, parameters $\delta''_i, G''_i, \e''_i$ are given by the corresponding variant of relation 
(\ref{dvad}). 

{\bf 2.3.} To construct the $(\bar{h}_i, \bar{k}_i, \bar{\delta}_i, \bar{G}_i, \bar{\e}_i)$-expansion for the non-absorption probability $\bar{p}_{ii}(\e)$ using relation (\ref{expabas}), and proposition {\bf (i)}  of  Lemma 2. In this case, parameters $ \bar{\delta}_i, \bar{G}_i, \bar{\e}_i$ are given by the corresponding variant of relation (\ref{equar}). \vspace{1mm}

{\bf Case 2}: $i \in \overline{{\mathbb Y}}$.  \vspace{1mm}

{\bf 2.4.}  In this case, the non-absorption probability $\bar{p}_{ii}(\e)  \equiv 1$. If necessary, it can be represented in the form of $(0, n, 1, G, \e_0)$-expansion given by relation(\ref{vopit}), for any $0 < G < \infty$ and $n = 0, 1, \ldots$. \vspace{1mm}

The above remarks  can be summarized in the following theorem. \vspace{1mm}

{\bf Lemma  11}. {\em Let conditions ${\bf A}$, ${\bf B'}$ and  ${\bf C}$ hold. Then, the asymptotic expansions for the non-absorption probabilities $\bar{p}_{ii}(\e) , i \in \XX$ with explicit upper bounds for remainders are given in Algorithm 2.} \vspace{1mm}

{\bf 6.2. Asymptotic expansions for transition probabilities  of 
reduced embedded Markov chains}.
Relation (\ref{transit}) can be re-written in the following form more convenient for constructing  asymptotic expansions for probabilities $_rp_{ij}(\e), i, j \in \, _r\XX$, 
\begin{equation}\label{transitak}
_rp_{ij}(\e) =  \left \{
\begin{array}{cl}
p_{ij}(\e) + p_{ir}(\e) \frac{p_{rj}(\e)}{1 - p_{rr}(\e)}    & \ \text{if} \ j \in \YY_{ir}^+  \cap \YY_{ir}^-, \vspace{1mm} \\
p_{ij}(\e)  & \ \text{if} \  j \in  \YY_{ir}^+ \cap \overline{\YY}_{ir}^-,  \vspace{1mm} \\
p_{ir}(\e) \frac{p_{rj}(\e)}{1 - p_{rr}(\e)} & \ \text{if} \  j \in  \overline{\YY}_{ir}^+ \cap \YY_{ir}^-, \vspace{1mm}  \\
0 & \ \text{if} \  j \in \,   \overline{\YY}_{ir}^+ \cap  \overline{\YY}_{ir}^- = \, _r\overline{\YY}_{i}.
\end{array}
\right.
\end{equation}

{\bf Algorithm 3.} This is an algorithm for constructing asymptotic expansions for transition probabilities   
$_rp_{ij}(\e), i, j \in \, _r\XX$. \vspace{1mm}

{\bf Case 1: $r \in \YY$}. \vspace{1mm}

{\bf 3.1.} To construct $(\tilde{h}_{rj}, \, \tilde{k}_{rj})$-expansions for conditional probabilities $\tilde{p}_{rj}(\e)$ $= \frac{p_{rj}(\e)}{1 -  p_{rr}(\e)}, \ j \in \YY_{rr}^+$,  using the $(l^-_{rj}, l^+_{rj})$-expansions for transition probabilities $p_{rj}(\e)$ given in condition ${\bf B}$, the 
$(\bar{h}_r, \bar{k}_r)$-expansion  for the non-absorption probability $\bar{p}_{rr}(\e) = 1 -  p_{rr}(\e)$ given in  Algorithm 1, and the proposition {\bf (v)} (the division rule) of Lemma 3. In this case, parameters $\tilde{h}_{rj}  = l^-_{rj} - \bar{h}_{r},   \tilde{k}_{rj} = (l^+_{rj} -  \bar{h}_{r}) \wedge (l^-_{rj}  +  \bar{k}_{r}  - 2\bar{h}_{r}),  \ j \in \YY_{rr}^+$.  These asymptotic  expansions are pivotal.

{\bf 3.2.} To construct $(\grave{h}_{irj}, \grave{k}_{irj})$-expansions for products $\grave{p}_{irj}(\e) = p_{ir}(\e)\tilde{p}_{rj}(\e)$   $ =   p_{ir}(\e)$  $\cdot \frac{p_{rj}(\e)}{1 -  p_{rr}(\e)}, \ j \in \YY_{ir}^-, \ i \in \, _r\XX$, using the $(l^-_{ir}, l^+_{ir})$-expansions for transition probabilities $p_{ir}(\e)$  given  in condition ${\bf B}$, the $(\tilde{h}_{rj}, \tilde{k}_{rj})$-expansions for conditional probabilities  $\tilde{p}_{rj}(\e)$ given in the above Step 3.1
and the proposition {\bf (iii)} (the multiplication rule) of Lemma 3. In this case, parameters $\grave{h}_{irj} =  l^-_{ir} + \tilde{h}_{rj}, \grave{k}_{irj} = (l^-_{ir} + \tilde{k}_{rj}) \wedge (l^+_{ir} + \tilde{h}_{rj}),  \ j \in \YY_{ir}^-, \ i \in \, _r\XX$.   These asymptotic  expansions are pivotal. 

{\bf 3.3.} To construct $(\acute{h}_{irj}, \acute{k}_{irj})$-expansions for sums $\acute{p}_{irj}(\e) = p_{ij}(\e) + \grave{p}_{irj}(\e)    =$ $p_{ij}(\e) +  p_{ir}(\e) \cdot \frac{p_{rj}(\e)}{1 -  p_{rr}(\e)}, \  j \in \YY_{ir}^+  \cap \YY_{ir}^-, \ i \in \, _r\XX$, using the $(l^-_{ij}, l^+_{ij})$-expansions for transition probabilities  $p_{ij}(\e)$ given in condition ${\bf B}$,  the $(\grave{h}_{irj}$, $\grave{k}_{irj})$-expansions for quantities  $\grave{p}_{irj}(\e)$ given in the above Step 3.2
and the proposition {\bf (ii)} (the summation rule) of Lemma 3. In this case, parameters $\acute{h}_{irj} =  l^-_{ij} \wedge \grave{h}_{irj}, \acute{k}_{irj} = l^+_{ij} \wedge \grave{k}_{irj},  \ j \in \YY_{ir}^-, \ i \in \, _r\XX$.  These asymptotic expansions are pivotal. 

{\bf 3.4.} To construct $(_rl_{ij}^-, \,  _rl_{ij}^+)$-expansions for transition probabilities $_rp_{ij}(\e)$ $= \sum_{l = \, _rl_{ij}^-}^{_rl_{ij}^+} \, _ra_{ij}[l]\e^l + 
o(\e^{_rl_{ij}^+}), i, j \in \, _r\XX$, using the $(l_{ij}^-, l_{ij}^+)$-expansions for transition probabilities $p_{ij}(\e)$ given in condition ${\bf B}$,  the 
$(\grave{h}_{irj}, \grave{k}_{irj})$-expansions for quantities  $\grave{p}_{irj}(\e)$ and  $(\acute{h}_{irj}, \acute{k}_{irj})$-expansions for quantities $\acute{p}_{irj}(\e)$ given, respectively, in the above Steps 3.2 and 3.3,  and the corresponding variants of formulas for  transition probabilities $_rp_{ij}(\e)$ given in relation (\ref{transitak}). In this case, parameters  $_rl_{ij}^- = \acute{h}_{irj}, \, _rl_{ij}^+ = \acute{k}_{irj}$ if $j \in \YY_{ir}^+  \cap \YY_{ir}^-, \ i \in \, _r\XX$, or $_rl_{ij}^- = l_{ij}^-, \, _rl_{ij}^+ = l_{ij}^+$  if $j \in  \YY_{ir}^+ \cap \overline{\YY}_{ir}^-,  \ i \in \, _r\XX$, or $_rl_{ij}^- = \grave{h}_{irj}, \, _rl_{ij}^+ = \grave{k}_{irj}$  if $j \in  \overline{\YY}_{ir}^+ \cap \YY_{ir}^-,   \ i \in \, _r\XX$. These asymptotic expansions are pivotal. \vspace{1mm}

{\bf Case 2: $r  \in \overline{\YY}$}. \vspace{1mm}

{\bf 3.5.} The corresponding algorithm is a particular case of the algorithm given in Steps 3.1 -- 3.4. In this case the non-absorption probability 
$ \bar{p}_{rr}(\e)  = 1 -  p_{rr}(\e) \equiv 1$ and, thus, conditional probabilities  $\tilde{p}_{rj}(\e)= \frac{p_{rj}(\e)}{1 -  p_{rr}(\e)} = p_{rj}(\e)$, $j \in {\mathbb Y}^+_{rr}$. This permits one replace the $(\tilde{h}_{rj}, \, \tilde{k}_{rj})$-expansions for conditional probabilities  $\tilde{p}_{rj}(\e)$ by the $(l^-_{rj}, l^+_{rj})$-expansions for transition probabilities $p_{rj}(\e)$. This is the only change in the algorithm for construction of asymptotic expansions for transition probabilities $_rp_{ij}(\e), i, j \in \, _r\XX$ given in Steps 3.1 -- 3.4, which is required.

The above remarks  can be summarized in the following theorem. \vspace{1mm}

{\bf Theorem 2}. {\em Conditions ${\bf A}$, ${\bf B}$ and  ${\bf C}$ assumed to hold for the Markov chain $\eta^{(\e)}_n$, also hold for the reduced Markov chain $_r\eta^{(\e)}_n$, for every $r \in \XX$. The asymptotic expansions penetrating condition ${\bf B}$ are given  for transition probabilities $_rp_{ij}(\e) , j \in \, _r\YY_{i}, \ i \in \, _r\XX, r \in \XX$ in Algorithm 3.} \vspace{1mm}

{\bf Algorithm 4.} This is an algorithm for computing upper bounds for remainders in  asymptotic expansions for 
transition probabilities  $_rp_{ij}(\e), i, j \in \, _r\XX$.  \vspace{1mm}

{\bf 4.1.} To construct $(\tilde{h}_{rj}, \, \tilde{k}_{rj}, \tilde{\delta}_{rj}, \tilde{G}_{rj}, \tilde{\e}_{rj})$-expansions for conditional probabilities $\tilde{p}_{rj}(\e)$ $= \frac{p_{rj}(\e)}{1 -  p_{rr}(\e)}$,  $j \in \YY_{rr}^+$ using the $(l^-_{rj}, l^+_{rj}, \delta_{rj},  G_{rj},  \e_{rj})$-expansions for transition probabilities $p_{rj}(\e)$ given in condition ${\bf B'}$, the $(\bar{h}_r, \bar{k}_r, \bar{\delta}_{rj}, \bar{G}_{rj}, \bar{\e}_{rj})$-expansion  for the non-absorption probability $\bar{p}_{rr}(\e) = 1 -  p_{rr}(\e)$ given in  Algorithm 2, and the proposition {\bf (v)} (the division rule) of Lemma 4. In this case, parameters $\tilde{\delta}_{rj}, \tilde{G}_{rj}, \tilde{\e}_{rj},  \ j \in \YY_{rr}^+$ are given by the corresponding variants of relation 
(\ref{hoputr}).  

{\bf 4.2.} To construct $(\grave{h}_{irj}, \grave{k}_{irj}, \grave{\delta}_{irj}, \grave{G}_{irj}, \grave{\e}_{irj})$-expansions for products $\grave{p}_{irj}(\e)$ $= p_{ir}(\e)\tilde{p}_{rj}(\e)$   $ =   p_{ir}(\e)$  $\cdot \frac{p_{rj}(\e)}{1 -  p_{rr}(\e)}, \ j \in \YY_{ir}^-, \ i \in \, _r\XX$, using the $(l^-_{ir}, l^+_{ir}, \delta_{ir},  G_{ir},  \e_{ir})$-expansions for transition probabilities $p_{ir}(\e)$  given  in condition ${\bf B'}$, the $(\tilde{h}_{rj}, \tilde{k}_{rj}$, $\tilde{\delta}_{rj}, \tilde{G}_{rj}, \tilde{\e}_{rj})$-expansions for conditional probabilities  $\tilde{p}_{rj}(\e)$ given in the above Step 4.1
and the proposition {\bf (iii)} (the multiplication rule) of Lemma 4. In this case, parameters $\grave{\delta}_{irj}, \grave{G}_{irj}, \grave{\e}_{irj},  \ j \in \YY_{ir}^-, \ i \in \, _r\XX$ are given by the corresponding variants of relation (\ref{profa}).

{\bf 4.3.} To construct $(\acute{h}_{irj}, \acute{k}_{irj}, \acute{\delta}_{irj}, \acute{G}_{irj}, \acute{\e}_{irj})$-expansions for sums $\acute{p}_{irj}(\e) = p_{ij}(\e) + \grave{p}_{irj}(\e)    =$ $p_{ij}(\e) +  p_{ir}(\e) \cdot \frac{p_{rj}(\e)}{1 -  p_{rr}(\e)}, \  j \in \YY_{ir}^+  \cap \YY_{ir}^-, \ i \in \, _r\XX$, using the $(l^-_{ij}, l^+_{ij}, \delta_{ij},  G_{ij},  \e_{ij})$-expansions for transition probabilities  $p_{ij}(\e)$ given in condition ${\bf B'}$,  the $(\grave{h}_{irj}$, $\grave{k}_{irj}, \grave{\delta}_{irj}, \grave{G}_{irj}, \grave{\e}_{irj})$-expansions for quantities  $\grave{p}_{irj}(\e)$ given in the above Step 4.2
and the proposition {\bf (ii)} (the summation rule) of Lemma 4. In this case, parameters $\acute{\delta}_{irj}, \acute{G}_{irj}, \acute{\e}_{irj},  \ j \in \YY_{ir}^-, \ i \in \, _r\XX$  are given by the corresponding variants of relation (\ref{profac}).  

{\bf 4.4.} To construct $(_rl_{ij}^-, \,  _rl_{ij}^+, \, _r\delta_{ij}, \, _rG_{ij}, \, _r\e_{ij})$-expansions for transition probabilities $_rp_{ij}(\e)$ $= \sum_{l = \, _rl_{ij}^-}^{_rl_{ij}^+} \, _ra_{ij}[l]\e^l + o(\e^{_rl_{ij}^+}), i, j \in \, _r\XX$ using the  $(l_{ij}^-, l_{ij}^+, \delta_{ij}, G_{ij}$, $\e_{ij})$-expansions for transition probabilities $p_{ij}(\e)$ given in condition ${\bf B'}$,  the $(\grave{h}_{irj}, \grave{k}_{irj},  \grave{\delta}_{irj},  \grave{G}_{irj},  \grave{\e}_{irj})$-expansions for quantities  $\grave{p}_{irj}(\e)$ and  $(\acute{h}_{irj}, \acute{k}_{irj}, \acute{\delta}_{irj}$, $\acute{G}_{irj}, \acute{\e}_{irj})$-expansions for quantities $\acute{p}_{irj}(\e)$ given, respectively, in the above Steps 4.2 and 4.3,  and the corresponding variants of formulas for  transition probabilities $_rp_{ij}(\e)$ given in relation (\ref{transitak}). In this case, parameters  $_r\delta_{ij} = \acute{\delta}_{irj}, \, _rG_{ij} = \acute{G}_{irj},  \, _r\e_{ij} =  \acute{\e}_{irj}$ if $j \in \YY_{ir}^+  \cap \YY_{ir}^-, \ i \in \, _r\XX$, or 
$_r\delta_{ij}= \delta_{ij}, \, _rG_{ij}  = G_{ij}, \, _r\e_{ij} = \e_{ij}$  if $j \in  \YY_{ir}^+ \cap \overline{\YY}_{ir}^-,  \ i \in \, _r\XX$, or $_r\delta_{ij} = 
\grave{\delta}_{irj}, \, _rG_{ij} = \grave{G}_{irj}, \, _r\e_{ij} =  \grave{\e}_{irj}$  if $j \in  \overline{\YY}_{ir}^+ \cap \YY_{ir}^-,   \ i \in \, _r\XX$.  \vspace{1mm}

{\bf Case 2: $r  \in \overline{\YY}$}. \vspace{1mm}

{\bf 4.5.} The corresponding algorithm is a particular case of the algorithm given in Steps 4.1 -- 4.4. In this case the non-absorption probability 
$ \bar{p}_{rr}(\e)  = 1 -  p_{rr}(\e) \equiv 1$ and, thus, conditional probabilities  $\tilde{p}_{rj}(\e)= \frac{p_{rj}(\e)}{1 -  p_{rr}(\e)} = p_{rj}(\e), i \in {\mathbb Y}^+_{rr}$. This permits one replace the $(\tilde{h}_{rj}, \, \tilde{k}_{rj}, \tilde{\delta}_{rj}, \tilde{G}_{rj}, \tilde{\e}_{rj})$-expansions for conditional probabilities  $\tilde{p}_{rj}(\e)$ by the $(l^-_{rj}, l^+_{rj},  \delta_{rj},  G_{rj},  \e_{rj})$-expansions for transition probabilities $p_{rj}(\e)$. This is the only change in the algorithm for construction of asymptotic expansions for transition probabilities $_rp_{ij}(\e), i, j \in \, _r\XX$ given in Steps 4.1 -- 4.4, which is required.

The above remarks  can be summarized in the following theorem. \vspace{1mm}

{\bf Theorem 3}. {\em Conditions ${\bf A}$, ${\bf B'}$ and  ${\bf C}$  assumed to hold for the Markov chain $\eta^{(\e)}_n$, also hold for the reduced Markov chain $_r\eta^{(\e)}_n$, for every $r \in \XX$. The upper bounds for remainders in asymptotic expansions penetrating condition ${\bf B'}$ are given  for transition probabilities $_rp_{ij}(\e) , j \in \, _r\YY_{i}, \ i \in \, _r\XX, r \in \XX$ in Algorithm 4.} \\

{\bf 6.3. Asymptotic expansions for expectations of sojourn  times for reduced semi-Markov processes.} 
Relation (\ref{expectaga}) can be re-written in the following form more convenient for constructing the corresponding asymptotic expansions probabilities $_re_{ij}(\e)$, $i, j \in \, _r\XX$, 
\begin{equation}\label{transitakne}
_re_{ij}(\e) =  \left \{
\begin{array}{ll}
e_{ij}(\e)  + e_{ir}(\e) \frac{p_{rj}(\e)}{1 - p_{rr}(\e)}   & \vspace{1mm} \\
+ \, e_{rr}(\e) \frac{p_{ir}(\e)}{1 - p_{rr}(\e)} \frac{p_{rj}(\e) }{1 - p_{rr}(\e)} & \vspace{1mm} \\
+  e_{rj}(\e)  \frac{p_{ir}(\e)}{1 - p_{rr}(\e)}   & \ \text{if} \ j \in \YY_{ir}^+  \cap \YY_{ir}^-,  \vspace{2mm} \\
e_{ij}(\e)  & \ \text{if} \  j \in  \YY_{ir}^+ \cap \overline{\YY}_{ir}^-,  \vspace{2mm} \\
e_{ir}(\e) \frac{p_{rj}(\e)}{1 - p_{rr}(\e)}  \vspace{1mm} \\
+ \, e_{rr}(\e)\frac{p_{ir}(\e)}{1 - p_{rr}(\e)} \frac{p_{rj}(\e) }{1 - p_{rr}(\e)} & \vspace{1mm} \\
+  e_{rj}(\e)  \frac{p_{ir}(\e)}{1 - p_{rr}(\e)}  & \ \text{if} \  j \in  \overline{\YY}_{ir}^+ \cap \YY_{ir}^-, \vspace{2mm} \\
0 & \ \text{if} \  j \in \,   \overline{\YY}_{ir}^+ \cap  \overline{\YY}_{ir}^-.
\end{array}
\right.
\end{equation}

{\bf Algorithm 5.} This is an algorithm for computing asymptotic expansions for 
expectations $_re_{ij}(\e), i, j \in \, _r\XX$.  \vspace{1mm}

{\bf Case 1: $r \in \YY$}. \vspace{1mm}

{\bf 5.1.} To construct $(\hat{h}_{ir}, \, \hat{k}_{ir})$-expansions  for quantities $\hat{p}_{ir}(\e)  = \frac{p_{ir}(\e)}{\bar{p}_{rr}(\e)} 
= \frac{p_{ir}(\e) }{1 - p_{rr}(\e)}, i \in \, _r\XX$, using the  $(l^-_{ir}, \, l^+_{ir})$-expansions for transition probabilities  $p_{ir}(\e)$ given in condition 
${\bf B}$, the $(\bar{h}_{r}, \, \bar{k}_{r})$-expansion for the non-absorption probability  $\bar{p}_r(\e) = 1 - p_{rr}(\e)$ given in Algorithm 1, and  proposition {\bf (v)} (the division rule) of Lemma 3. In this case, parameters $\hat{h}_{ir} = h_{ir} - \bar{h}_{r}, \hat{k}_{ir} = 
(h_{ir} + \bar{k}_{r} - 2 \bar{h}_{r}) \wedge (k_{ir} - \bar{h}_{r}), \ i \in \, _r\XX$.  These asymptotic  expansions are pivotal. 

{\bf 5.2.} To construct $(\check{h}_{irj}, \, \check{k}_{irj})$-expansions for  products $ \check{p}_{irj}(\e)  = \hat{p}_{ir}(\e) \tilde{p}_{rj}(\e)$ 
$= \frac{p_{ir}(\e)}{1 - p_{rr}(\e)}  \frac{p_{rj}(\e) }{1 - p_{rr}(\e)}, \ j \in  \YY_{ir}^-, \ i \in \, _r\XX$, using   $(\hat{h}_{ir}, \, \hat{k}_{ir})$-expansions  given in the above Step 5.1, the $(\tilde{h}_{rj}, \, \tilde{k}_{rj})$-expansions for conditional probabilities  $\tilde{p}_{rj}(\e) = \frac{p_{rj}(\e)}{1 - p_{rr}(\e)}$ given in Step 3.1 of Algorithm 3, and the proposition {\bf (iii)} (the multiplication rule) of Lemma 3. In this case, parameters $\check{h}_{ir} = \hat{h}_{ir} + \tilde{h}_{rj}, \check{k}_{ir} = (\hat{h}_{ir} + \tilde{k}_{rj}) \wedge (\hat{k}_{ir} + \tilde{h}_{rj}),  \ j \in  \YY_{ir}^-, \ i \in \, _r\XX$. 
These asymptotic  expansions are pivotal. 

{\bf 5.3.} To construct $(\dot{\tilde{h}}_{irj}, \, \dot{\tilde{k}}_{irj})$-expansions  for products $\tilde{e}_{irj}(\e) = e_{ir}(\e) \tilde{p}_{rj}(\e)$ 
$=  e_{ir}(\e)\frac{p_{rj}(\e)}{1 - p_{rr}(\e)},  \ j \in  \YY_{ir}^-, \ i \in \, _r\XX, $ using the  $(m^-_{ir}, \, m^+_{ir})$-expansions for expectations  
$e_{ir}(\e)$ given in condition ${\bf F}$, the $(\tilde{h}_{rj}, \, \tilde{k}_{rj})$-expansions for conditional probabilities  $\tilde{p}_{rj}(\e) = \frac{p_{rj}(\e)}{1 - p_{rr}(\e)}$ given in Step 3.1 of Algorithm 3,  and the proposition {\bf (iii)} (the multiplication rule) of Lemma 3. In this case, parameters $\dot{\tilde{h}}_{irj} = 
m^-_{ir} + \tilde{h}_{rj}, \dot{\tilde{k}}_{irj} = (m^-_{ir} + \tilde{k}_{rj}) \wedge (m^+_{ir}  + \tilde{h}_{rj}),  \ j \in  \YY_{ir}^-, \ i \in \, _r\XX$. 
These asymptotic  expansions are pivotal. 

{\bf 5.4.} To construct $(\dot{\check{h}}_{irj}, \, \dot{\check{k}}_{irj})$-expansions  for products $\check{e}_{irj}(\e) = e_{rr}(\e)$ $\cdot \check{p}_{irj}(\e)$ 
$=  e_{rr}(\e)\frac{p_{ir}(\e)}{1 - p_{rr}(\e)}  \frac{p_{rj}(\e)}{1 - p_{rr}(\e)},  \ j \in  \YY_{ir}^-, \ i \in \, _r\XX, $ using the  $(m^-_{rr}, \, m^+_{rr})$-expansions for expectations  $e_{rr}(\e)$ given in condition ${\bf F}$, the $(\check{h}_{irj}, \, \check{k}_{irj})$-expan\-sions for  quantities  
$\check{p}_{irj}(\e)$ given in the above Step 5.2,  and the proposition {\bf (iii)} (the multiplication rule) of Lemma 3. In this case, parameters 
$\dot{\check{h}}_{irj} = m^-_{rr} + \check{h}_{irj}, \dot{\check{k}}_{irj} = (m^-_{rr} + \check{k}_{irj}) \wedge (m^+_{rr}  + \check{h}_{irj}),  \ j \in  \YY_{ir}^-, \ i \in \, _r\XX$. These asymptotic  expansions are pivotal. 

{\bf 5.5.} To construct $(\dot{\hat{h}}_{irj}, \, \dot{\hat{k}}_{irj})$-expansions  for products $\hat{e}_{irj}(\e) = e_{rj}(\e)$ $\cdot \hat{p}_{ir}(\e)$ 
$=  e_{rj}(\e) \frac{p_{ir}(\e)}{1 - p_{rr}(\e)},  \ j \in  \YY_{ir}^-, \ i \in \, _r\XX, $ using the  $(m^-_{rj}, \, m^+_{rj})$-expansions for expectations  $e_{rj}(\e)$ given in condition ${\bf F}$, the $(\hat{h}_{ir}, \, \hat{k}_{ir})$-expansions for quantities  
$\hat{p}_{ir}(\e)$ given in the above Step 5.1,  and the proposition {\bf (iii)} (the multiplication rule) of Lemma 3. In this case, parameters 
$\dot{\hat{h}}_{irj} = m^-_{rj} + \hat{h}_{ir}, \dot{\hat{k}}_{irj} = (m^-_{rj} + \hat{k}_{ir}) \wedge (m^+_{rj}  + \hat{h}_{ir}),  \ j \in  \YY_{ir}^-, \ i \in \, _r\XX$. These asymptotic  expansions are pivotal. 

{\bf 5.6.} To construct $(\ddot{h}_{irj}, \, \ddot{k}_{irj})$-expansions for sums $\ddot{e}_{irj}(\e) = \tilde{e}_{irj}(\e) + \check{e}_{irj}(\e) + 
\hat{e}_{irj}(\e) = e_{ir}(\e) \frac{p_{rj}(\e)}{1 - p_{rr}(\e)}  +  e_{rr}(\e)\frac{p_{ir}(\e)}{1 - p_{rr}(\e)} \frac{p_{rj}(\e) }{1 - p_{rr}(\e)}
+  e_{rj}(\e)  \frac{p_{ir}(\e)}{1 - p_{rr}(\e)}$, using the $(\dot{\tilde{h}}_{irj}, \, \dot{\tilde{k}}_{irj})$-expansions for quantities $\tilde{e}_{irj}(\e)$, the 
$(\dot{\check{h}}_{irj}, \, \dot{\check{k}}_{irj})$-expansions  for quantities $\check{e}_{irj}(\e)$ and the $(\dot{\hat{h}}_{irj}, \, \dot{\hat{k}}_{irj})$-expansions  for  quantities $\hat{e}_{irj}(\e)$ given, respectively, in the above Steps 5.3, 5.4 and 5.5, and the proposition {\bf (i)} (the summation rule) of Lemma 5. In this case, parameters $\ddot{h}_{irj} = \dot{\tilde{h}}_{irj} \wedge \dot{\check{h}}_{irj} \wedge \dot{\hat{h}}_{irj}, \ddot{k}_{irj} = \dot{\tilde{k}}_{irj} \wedge \dot{\check{k}}_{irj} \wedge \dot{\hat{k}}_{irj},  \ j \in  \YY_{ir}^-, \ i \in \, _r\XX$. These asymptotic expansions are pivotal.

{\bf 5.7.} To construct $(\dddot{h}_{irj}, \, \dddot{k}_{irj})$-expansions for sums $\dddot{e}_{irj}(\e) =  e_{ij}(\e)  +  \ddot{e}_{irj}(\e)  =  e_{ij}(\e)  + e_{ir}(\e) \frac{p_{rj}(\e)}{1 - p_{rr}(\e)}  +  e_{rr}(\e)\frac{p_{ir}(\e)}{1 - p_{rr}(\e)} \frac{p_{rj}(\e) }{1 - p_{rr}(\e)}
+  e_{rj}(\e)  \frac{p_{ir}(\e)}{1 - p_{rr}(\e)}$, using   $(m^-_{ij}, \, m^+_{ij})$-expansions for expectations  $e_{ij}(\e)$ given in condition ${\bf F}$, 
the $(\ddot{h}_{irj}, \, \ddot{k}_{irj})$-expansions-expansions  for  quantities $\ddot{e}_{irj}(\e)$ given in the above Step 5.6, and the proposition {\bf (ii)} (the summation rule) of Lemma 3. In this case, parameters $\dddot{h}_{irj} = m^-_{ij} \wedge \ddot{\tilde{h}}_{irj}, \dddot{k}_{irj} = m^+_{ij} \wedge \ddot{k}_{irj},  \ j \in  \YY_{ir}^-, \ i \in \, _r\XX$. These asymptotic expansions are pivotal.

{\bf 5.8.} To construct $(_rm_{ij}^-, \, _rm_{ij}^+)$-expansions for probabilities $_re_{ij}(\e)  =$ $\sum_{l = \, _rm_{ij}^-}^{_rm_{ij}^+} \, _rb_{ij}[l]\e^l + o(\e^{_rm_{ij}^+}),  i, j \in \, _r\XX$ using the asymptotic expansions  for expectations $e_{ij}(\e)$ given in condition ${\bf F}$, the 
$(\ddot{h}_{irj}, \, \ddot{k}_{irj})$-expansions for quantities  $\ddot{e}_{irj}(\e)$ and  $(\dddot{h}_{irj}, \, \dddot{k}_{irj})$-expansions for quantities 
$\dddot{e}_{irj}(\e)$ given, respectively,  in Steps 5.6  and 5.7,  and the corresponding variants of formulas for  expectations $_re_{ij}(\e)$ given in relation (\ref{transitakne}). In this case, parameters  $_rm_{ij}^- = \dddot{h}_{irj}, \, _rm_{ij}^+ = \dddot{k}_{irj}$ if $j \in \YY_{ir}^+  \cap \YY_{ir}^-, \ i \in \, _r\XX$, or $_rm_{ij}^- = m_{ij}^-, \, _rm_{ij}^+ = m_{ij}^+$  if $j \in  \YY_{ir}^+ \cap \overline{\YY}_{ir}^-,  \ i \in \, _r\XX$, or $_rm_{ij}^- = \ddot{h}_{irj}, \, _rm_{ij}^+ = \ddot{k}_{irj}$  if $j \in  \overline{\YY}_{ir}^+ \cap \YY_{ir}^-,   \ i \in \, _r\XX$. These asymptotic expansions are pivotal. \vspace{1mm}

{\bf Case 2: $r \in \overline{\YY}$}. \vspace{1mm}

{\bf 5.9.} The corresponding algorithm is a particular case of the algorithm given in Steps 5.1 -- 5.8. In this case the non-absorption probability 
$ \bar{p}_{rr}(\e)  = 1 -  p_{rr}(\e) \equiv 1$ and, thus, conditional probabilities  $\tilde{p}_{rj}(\e)= \frac{p_{rj}(\e)}{1 -  p_{rr}(\e)} = p_{rj}(\e)$, $j \in {\mathbb Y}^+_{rr}$ and quantities  $\hat{p}_{ir}(\e)  = \frac{p_{ir}(\e)}{1 -  p_{rr}(\e)} =  p_{ir}(\e), \, i \in \, _r\XX$. This permits one replace the $(\tilde{h}_{rj}, \, \tilde{k}_{rj})$-expansions for conditional probabilities  $\tilde{p}_{rj}(\e)$ by the $(l^-_{rj}, l^+_{rj})$-expansions for transition probabilities $p_{rj}(\e)$ and 
the $(\hat{h}_{rj}, \, \hat{k}_{rj})$-expansions for quantities  $\hat{p}_{rj}(\e)$ by the $(l^-_{ir}, l^+_{ir})$-expansions for transition probabilities $p_{ir}(\e)$. These are the only changes in the algorithm for construction of asymptotic expansions for expectations $_re_{ij}(\e), i, j \in \, _r\XX$ given in Steps 5.1 -- 5.8, which are required.

The above remarks  can be summarized in the following theorem. \vspace{1mm}

{\bf Theorem 4}. {\em Conditions ${\bf A}$  -- ${\bf F}$ assumed to hold for the semi-Markov process $\eta^{(\e)}(t)$, also hold for the reduced semi-Markov process $_r\eta^{(\e)}(t)$, for every $r \in \XX$. The asymptotic expansions penetrating conditions  ${\bf B}$ and  ${\bf F}$ are given  for transition probabilities $_rp_{ij}(\e) , j \in \, _r\YY_{i}, \ i \in \, _r\XX, r \in \XX$ and expectations $_re_{ij}(\e) , j \in \, _r\YY_{i}, \ i \in \, _r\XX, r \in \XX$ in Algorithms 3 and 5.} \vspace{1mm}

{\bf Algorithm 6.} This is an algorithm for computing  upper bounds for remainders in asymptotic expansions for 
 expectations $_re_{ij}(\e), i, j \in \, _r\XX$.  \vspace{1mm}

{\bf 6.1.} To construct $(\hat{h}_{ir}, \, \hat{k}_{ir},  \hat{\delta}_{ir},  \hat{G}_{ir},  \hat{\e}_{ir})$-expansions  for quantities $\hat{p}_{ir}(\e)  = \frac{p_{ir}(\e)}{\bar{p}_{rr}(\e)} = \frac{p_{ir}(\e) }{1 - p_{rr}(\e)}, i \in \, _r\XX$, using the  $(l^-_{ir}, \, l^+_{ir}, \delta_{ij},  G_{ij},  \e_{ij})$-expansions for transition probabilities  $p_{ir}(\e)$ given in condition ${\bf B'}$, the $(\bar{h}_{r}, \, \bar{k}_{r}, \bar{\delta}_{ir},  \bar{G}_{ir},  \bar{\e}_{ir})$-expansi\-on for the non-absorption probability  $\bar{p}_r(\e) = 1 - p_{rr}(\e)$ given in Algorithm 2, and  proposition {\bf (v)} (the division rule) of Lemma 4. In this case, parameters $\hat{\delta}_{ir},  \hat{G}_{ir},  \hat{\e}_{ir}, \ i \in \, _r\XX$ are given by the corresponding variants of relation 
(\ref{hoputr}).   

{\bf 6.2.} To construct $(\check{h}_{irj}, \, \check{k}_{irj}, \check{\delta}_{irj},  \check{G}_{irj},  \check{\e}_{irj})$-expansions for  products $ \check{p}_{irj}(\e)$  $= \hat{p}_{ir}(\e) \tilde{p}_{rj}(\e) = \frac{p_{ir}(\e)}{1 - p_{rr}(\e)}  \frac{p_{rj}(\e) }{1 - p_{rr}(\e)}, \ j \in  \YY_{ir}^-, \ i \in \, _r\XX$, using   $(\hat{h}_{ir}, \, \hat{k}_{ir}, \hat{\delta}_{ir},  \hat{G}_{ir},  \hat{\e}_{ir})$-expansions  given in the above Step 6.1, the $(\tilde{h}_{rj}, \, \tilde{k}_{rj},  \tilde{\delta}_{rj},  \tilde{G}_{rj},  \tilde{\e}_{rj})$-expansions for conditional probabilities  $\tilde{p}_{rj}(\e) = \frac{p_{rj}(\e)}{1 - p_{rr}(\e)}$ given in Step 4.1 of Algorithm 4, and the proposition {\bf (iii)} (the multiplication rule) of Lemma 4. In this case, parameters $\check{\delta}_{irj},  \check{G}_{irj},  \check{\e}_{irj},  \ j \in  \YY_{ir}^-, \ i \in \, _r\XX$ are given by the corresponding variants of relation (\ref{profa}).   

{\bf 6.3.} To construct $(\dot{\tilde{h}}_{irj}, \, \dot{\tilde{k}}_{irj},  \dot{\tilde{\delta}}_{irj},  \dot{\tilde{G}}_{irj},  \dot{\tilde{\e}}_{irj})$-expansions  for products $\tilde{e}_{irj}(\e)$ $= e_{ir}(\e) \tilde{p}_{rj}(\e) =  e_{ir}(\e)\frac{p_{rj}(\e)}{1 - p_{rr}(\e)},  \ j \in  \YY_{ir}^-, \ i \in \, _r\XX, $ using the  $(m^-_{ir}, \, m^+_{ir}, \delta_{ir}, G_{ir}$, $\e_{ir})$-expansions for expectations  
$e_{ir}(\e)$ given in condition ${\bf F'}$, the $(\tilde{h}_{rj}, \, \tilde{k}_{rj}, \tilde{\delta}_{rj}$, $\tilde{G}_{rj}, \tilde{\e}_{rj})$-expansions for conditional probabilities  $\tilde{p}_{rj}(\e) = \frac{p_{rj}(\e)}{1 - p_{rr}(\e)}$ given in Step 4.1 of Algorithm 4,  and the proposition {\bf (iii)} (the multiplication  rule) of Lemma 4. In this case, parameters $ \dot{\tilde{\delta}}_{irj},  \dot{\tilde{G}}_{irj},  \dot{\tilde{\e}}_{irj},  \ j \in  \YY_{ir}^-, \ i \in \, _r\XX$ are given by the corresponding variants of relation (\ref{profa}).    

{\bf 6.4.} To construct $(\dot{\check{h}}_{irj}, \, \dot{\check{k}}_{irj}, \dot{\check{\delta}}_{irj}, \dot{\check{G}}_{irj}, \dot{\check{\e}}_{irj})$-expansions  for products $\check{e}_{irj}(\e)$ $= e_{rr}(\e)$ $\cdot \check{p}_{irj}(\e) =  e_{rr}(\e)\frac{p_{ir}(\e)}{1 - p_{rr}(\e)}  \frac{p_{rj}(\e)}{1 - p_{rr}(\e)},  \ j \in  \YY_{ir}^-, \ i \in \, _r\XX, $ using the  $(m^-_{rr}, \, m^+_{rr}$, $\delta_{rr}, G_{rr}, \e_{rr})$-expansions for expectations  $e_{rr}(\e)$ given in condition ${\bf F'}$, the $(\check{h}_{irj}, \, \check{k}_{irj}, \check{\delta}_{irj}, \check{G}_{irj}, \check{\e}_{irj})$-expansions for  quantities  
$\check{p}_{irj}(\e)$ given in the above Step 6.2,  and the proposition {\bf (iii)} (the multiplication  rule) of Lemma 4. In this case, parameters 
$\dot{\check{\delta}}_{irj}, \dot{\check{G}}_{irj}, \dot{\check{\e}}_{irj}),  \ j \in  \YY_{ir}^-, \ i \in \, _r\XX$ are given by the corresponding variants of relation (\ref{profa}).     

{\bf 6.5.} To construct $(\dot{\hat{h}}_{irj}, \, \dot{\hat{k}}_{irj}, \dot{\hat{\delta}}_{irj}, \dot{\hat{G}}_{irj}, \dot{\hat{\e}}_{irj})$-expansions  for products $\hat{e}_{irj}(\e)$ $= e_{rj}(\e)$ $\cdot \hat{p}_{ir}(\e) =  e_{rj}(\e) \frac{p_{ir}(\e)}{1 - p_{rr}(\e)},  \ j \in  \YY_{ir}^-, \ i \in \, _r\XX, $ using the  $(m^-_{rj}, \, m^+_{rj}, \delta_{rj}, G_{rj}$, $\e_{rj})$-expansions for expectations  $e_{rj}(\e)$ given in condition ${\bf F'}$, the $(\hat{h}_{ir}, \, \hat{k}_{ir}, \hat{\delta}_{ir}$, 
$\hat{G}_{ir}, \hat{\e}_{ir})$-expansions for quantities  
$\hat{p}_{ir}(\e)$ given in the above Step 6.1,  and the proposition {\bf (iii)} (the multiplication  rule) of Lemma 4. In this case, parameters 
$\dot{\hat{\delta}}_{irj}, \dot{\hat{G}}_{irj}, \dot{\hat{\e}}_{irj}),  \ j \in  \YY_{ir}^-, \ i \in \, _r\XX$ are given by the corresponding variants of 
relation (\ref{profa}).     

{\bf 6.6.} To construct $(\ddot{h}_{irj}, \, \ddot{k}_{irj}, \ddot{\delta}_{irj}, \ddot{G}_{irj}, \ddot{\e}_{irj})$-expansions for 
sums $\ddot{e}_{irj}(\e) = \tilde{e}_{irj}(\e) + \check{e}_{irj}(\e) + 
\hat{e}_{irj}(\e) = e_{ir}(\e) \frac{p_{rj}(\e)}{1 - p_{rr}(\e)}  +  e_{rr}(\e)\frac{p_{ir}(\e)}{1 - p_{rr}(\e)} \frac{p_{rj}(\e) }{1 - p_{rr}(\e)}
+  e_{rj}(\e)  \frac{p_{ir}(\e)}{1 - p_{rr}(\e)}$, using the $(\dot{\tilde{h}}_{irj}, \dot{\tilde{k}}_{irj}, \dot{\tilde{\delta}}_{irj}, \dot{\tilde{G}}_{irj}, 
\dot{\tilde{\e}}_{irj})$-expansions for quantities $\tilde{e}_{irj}(\e)$, the 
$(\dot{\check{h}}_{irj}$,  $\dot{\check{k}}_{irj}, \dot{\check{\delta}}_{irj}, \dot{\check{G}}_{irj}, \dot{\check{\e}}_{irj})$-expansions  for quantities $\check{e}_{irj}(\e)$ and the $(\dot{\hat{h}}_{irj}$,  $\dot{\hat{k}}_{irj}, \dot{\hat{\delta}}_{irj}$, $\dot{\hat{G}}_{irj}, \dot{\hat{\e}}_{irj})$-expansions  for  quantities $\hat{e}_{irj}(\e)$ given, respectively, in the above Steps 6.3, 6.4 and 6.5, and the proposition {\bf (i)} (the summation rule) of Lemma 6. In this case, parameters $\ddot{\delta}_{irj}, \ddot{G}_{irj}, \ddot{\e}_{irj},  \ j \in  \YY_{ir}^-, \ i \in \, _r\XX$ are given by the corresponding variants of 
relation (\ref{dvad}).

{\bf 6.7.} To construct $(\dddot{h}_{irj}, \, \dddot{k}_{irj},  \dddot{\delta}_{irj},  \dddot{G}_{irj},  \dddot{\e}_{irj})$-expansions for sums $\dddot{e}_{irj}(\e)$ 
$=  e_{ij}(\e)  +  \ddot{e}_{irj}(\e)  =  e_{ij}(\e)  + e_{ir}(\e) \frac{p_{rj}(\e)}{1 - p_{rr}(\e)}  +  e_{rr}(\e)\frac{p_{ir}(\e)}{1 - p_{rr}(\e)} \frac{p_{rj}(\e) }{1 - p_{rr}(\e)}
+  e_{rj}(\e)  \frac{p_{ir}(\e)}{1 - p_{rr}(\e)}$, using   $(m^-_{ij},  m^+_{ij}, \dot{\delta}_{ij}, \dot{G}_{ij}, \dot{\e}_{ij})$-expansions for expectations  $e_{ij}(\e)$ given in condition ${\bf F'}$, the $(\ddot{h}_{irj}, \ddot{k}_{irj}, \ddot{\delta}_{irj}, \ddot{G}_{irj}, \ddot{\e}_{irj})$-expansions  for  quantities $\ddot{e}_{irj}(\e)$ given in the above Step 6.6, and the proposition {\bf (ii)} (the summation rule) of Lemma 4. In this case, parameters $ \dddot{\delta}_{irj},  \dddot{G}_{irj},  \dddot{\e}_{irj},  \ j \in  \YY_{ir}^-, \ i \in \, _r\XX$ are given by the corresponding variants of 
relation (\ref{profac}).

{\bf 6.8.} To construct $(_rm_{ij}^-, \,  _rm_{ij}^+, \, _r\dot{\delta}_{ij}, \, _r\dot{G}_{ij}, \, _r\dot{\e}_{ij})$-expansions expansions for expectations  $_re_{ij}(\e)  = \sum_{l = \, _rm_{ij}^-}^{_rm_{ij}^+} \, _rb_{ij}[l]\e^l + o(\e^{_rm_{ij}^+}),  i, j \in \, _r\XX$, using the $(m_{ij}^-$, $m_{ij}^+, \dot{\delta}_{ij}, 
\dot{G}_{ij}, \dot{\e}_{ij})$-expansions  for expectations $e_{ij}(\e)$ given in condition ${\bf F'}$, the  $(\ddot{h}_{irj}$, $\ddot{k}_{irj}, \ddot{\delta}_{irj}, \ddot{G}_{irj}, \ddot{\e}_{irj})$-expansions for quantities $\ddot{e}_{irj}(\e)$ and the $(\dddot{h}_{irj}$, $\dddot{k}_{irj}$,  $\dddot{\delta}_{irj}$,  $\dddot{G}_{irj},  \dddot{\e}_{irj})$-expansions for quantities 
$\dddot{e}_{irj}(\e)$ given, respectively,  in Steps 6.6 and 6.7,  and the corresponding variants of formulas for  expectations $_re_{ij}(\e)$ given in relation (\ref{transitakne}). In this case, parameters  $_r\dot{\delta}_{ij} = \dddot{\delta}_{irj}, \, _r\dot{G}_{ij} = \dddot{G}_{irj}, \, _r\dot{\e}_{ij} =  \dddot{\e}_{irj}$ if $j \in \YY_{ir}^+  \cap \YY_{ir}^-, \ i \in \, _r\XX$, or $_r\dot{\delta}_{ij} = \delta_{ij}, \, _r\dot{G}_{ij} = G_{ij},  \, _r\dot{\e}_{ij} = \e_{ij}$  if $j \in  \YY_{ir}^+ \cap \overline{\YY}_{ir}^-,  \ i \in \, _r\XX$, or $_r\dot{\delta}_{ij} = \ddot{\delta}_{irj}, \,  _r\dot{G}_{ij} = \ddot{G}_{irj}, \, _r\dot{\e}_{ij} =  \ddot{\e}_{irj}$  if $j \in  \overline{\YY}_{ir}^+ \cap \YY_{ir}^-,   \ i \in \, _r\XX$. \vspace{1mm}

{\bf Case 2: $r \in \overline{\YY}$}. \vspace{1mm}

{\bf 6.9.} The corresponding algorithm is a particular case of the algorithm given in Steps 5.1 -- 5.8. In this case the non-absorption probability 
$ \bar{p}_{rr}(\e)  = 1 -  p_{rr}(\e) \equiv 1$ and, thus, conditional probabilities  $\tilde{p}_{rj}(\e)= \frac{p_{rj}(\e)}{1 -  p_{rr}(\e)} = p_{rj}(\e)$, $j \in {\mathbb Y}^+_{rr}$ and quantities  $\hat{p}_{ir}(\e)  = \frac{p_{ir}(\e)}{1 -  p_{rr}(\e)} =  p_{ir}(\e), \, i \in \, _r\XX$. This permits one replace the $(\tilde{h}_{rj}, \, \tilde{k}_{rj}, \tilde{\delta}_{rj}, \tilde{G}_{rj}, \tilde{\e}_{rj})$-expansions for conditional probabilities  $\tilde{p}_{rj}(\e)$ by the $(l^-_{rj}, l^+_{rj},  \delta_{rj},  G_{rj},  \e_{rj})$-expansions for transition probabilities $p_{rj}(\e)$ and 
the $(\hat{h}_{rj},  \hat{k}_{rj}, \hat{\delta}_{rj}, \hat{G}_{rj}, \hat{\e}_{rj})$-expansions for quantities  $\hat{p}_{rj}(\e)$ by 
the $(l^-_{ir}, l^+_{ir}, \delta_{ir}, G_{ir}, \e_{ir})$-expansions for transition probabilities $p_{ir}(\e)$. These are the only changes in the algorithm for construction of asymptotic expansions for expectations $_re_{ij}(\e), i, j \in \, _r\XX$ given in Steps 6.1 -- 6.8, which are required.

The above remarks  can be summarized in the following theorem. \vspace{1mm}

{\bf Theorem 5}. {\em Conditions ${\bf A}$, ${\bf B'}$, ${\bf C}$ -- ${\bf E}$, ${\bf F'}$ assumed to hold for the semi-Markov process $\eta^{(\e)}(t)$, also hold for the reduced semi-Markov process $_r\eta^{(\e)}(t)$, for every $r \in \XX$. The upper bounds for remainders in expansions penetrating conditions  ${\bf B'}$ and  ${\bf F'}$ are given  for transition probabilities $_rp_{ij}(\e) , j \in \, _r\YY_{i}, \ i \in \, _r\XX, r \in \XX$ and expectations $_re_{ij}(\e) , j \in \, _r\YY_{i}, \ i \in \, _r\XX, r \in \XX$ in Algorithms 4 and 6.} \\

{\bf 7. Sequential reduction of phase space and asymptotic expansions for stationary distributions of perturbed semi-Markov processes} \\

In this section, we present algorithms of sequential reduction of phase spaces for semi-Markov processes and construction of 
asymptotic expansions for their stationary distributions. \vspace{1mm}

{\bf 7.1. Algorithms of sequential reduction of phase spaces for semi-Markov processes}.  Let $\eta^{(\e)}(t)$ be a semi-Markov process with the phase space $\XX = \{1, \ldots, N \}$, which satisfy conditions ${\bf A}$ -- ${\bf F}$. 

Let $\langle r_1, \ldots, r_N \rangle$ is a permutation of the sequence  $\langle 1, \ldots, N \rangle$, and $\bar{r}_n = (r_1, \ldots, r_n), n = 1, \ldots, N$ is the corresponding sequence of growing chains of states from space $\XX$. 

Let us choose  state $i \in \XX$, and a  permutation $\langle r_1, \ldots, r_N \rangle$ such that $r _N = i$. 

Let us also assume that initial distribution $\bar{p}{(\e)}$ is concentrated in the state $i$, i.e., $p_i^{(\e)} = 1$.
\vspace{1mm}

{\bf Algorithm 7.} This is an algorithm for sequential reduction of the phase space for the semi-Markov process $\eta^{(\e)}(t)$ and constructing  asymptotic expansions for transition probabilities and expectation of sojourn times for semi-Markov processes with reduced phase spaces. \vspace{1mm}

{\bf 7.1.} Let $_{\bar{r}_1}\eta^{(\e)}(t) = \, _{r_1}\eta^{(\e)}(t)$ be the reduced  semi-Markov process which is the result of reduction of state $r_1$ for the semi-Markov process  $\eta^{(\e)}(t)$. This semi-Markov process has the phase space $_{\bar{r}_1}\XX = \XX \setminus \{ r_1 \}$, transition probabilities of the embedded Markov chain $_{\bar{r}_1}p_{i'j'}(\e), i', j' \in \, _{r_1}\XX$ and expectations of transition times $_{\bar{r}_1}e_{i'j'}(\e), i', j' \in \, _{r_1}\XX$, which are determined by the transition probabilities and the expectations of transition times for the process $\eta^{(\e)}(t)$ via relations (\ref{transit}) and (\ref{expectaga}). According Theorem 1, the expectations of hitting times $E_{i'j'}(\e), i', j' \in \, _{\bar{r}_1}\XX$ coincide for the semi-Markov processes $\eta^{(\e)}(t)$ and $_{\bar{r}_1}\eta^{(\e)}(t)$. According Theorems 2 and 4, the  semi-Markov process $_{\bar{r}_1}\eta^{(\e)}(t)$ satisfy conditions ${\bf A}$ -- ${\bf F}$. The transition sets 
$_{\bar{r}_1}\YY_{i'}   = \, _{r_1}\YY_{i'}, i' \in \, _{\bar{r}_1}\XX$ are determined for the process $_{\bar{r}_1}\eta^{(\e)}(t)$ by condition ${\bf A}$  and relation (\ref{denojat}). Therefore,  the $(_{\bar{r}_1}l_{i'j'}^-, \, _{\bar{r}_1}l_{i'j'}^+)$-expansions for transition probabilities $_{\bar{r}_1}p_{i'j'}(\e), j' \in \,  _{\bar{r}_1}\YY_{i'}, i' \in \,  _{\bar{r}_1}\XX$ and $(_{\bar{r}_1}m_{i'j'}^-, \, _{\bar{r}_1}m_{i'j'}^+)$-expansions for expectations $_{\bar{r}_1}e_{i'j'}(\e), j' \in \,  _{\bar{r}_1}\YY_{i'}, i' \in \,  _{\bar{r}_1}\XX$ can be constructed by applying Algorithms 1, 3 and 5 to the $(l_{i'j'}^-, l_{i'j'}^+)$-expansions for transition probabilities $p_{i'j'}(\e), j' \in \YY_{i'}, i' \in \XX$ and $(m_{i'j'}^-, m_{i'j'}^+)$-expansions for expectations $e_{i'j'}(\e), j' \in \YY_{i'}, i' \in \XX$. These expansions ara pivotal.

{\bf 7.2.} Let $_{\bar{r}_2}\eta^{(\e)}(t)$  be the reduced  semi-Markov process which is the result of reduction of state  $r_2$ for the  semi-Markov process $_{\bar{r}_1}\eta^{(\e)}(t)$. This semi-Markov process has the phase space $_{\bar{r}_2}\XX = \XX \setminus \{ r_1, r_2 \}$, the transition probabilities of the embedded Markov chain $_{\bar{r}_2}p_{i'j'}(\e), i', j' \in _{r_1}\XX$ and the expectations of transition times $_{\bar{r}_2}e_{i'j'}(\e), i', j' \in \, 
_{r_2}\XX$, which are determined by the transition probabilities and the expectations of transition times for the process $_{\bar{r}_1}\eta^{(\e)}(t)$ via relations (\ref{transit}) and (\ref{expectaga}). According Theorem 1, the expectations of hitting times $E_{i'j'}(\e), i', j' \in \, _{\bar{r}_2}\XX$  coincide for the semi-Markov processes $\eta^{(\e)}(t)$,  $_{\bar{r}_1}\eta^{(\e)}(t)$ and $_{\bar{r}_2}\eta^{(\e)}(t)$. According Theorems 2 and 4, the  transition probabilities of the embedded Markov chain $_{\bar{r}_2}p_{i'j'}(\e), i', j' \in 
\, _{\bar{r}_2}\XX$ and the expectations of transition times $_{\bar{r}_2}e_{i'j'}(\e), i', j' \in \, _{\bar{r}_2}\XX$ satisfy conditions ${\bf A}$ -- ${\bf F}$.  The transition sets $_{\bar{r}_2}\YY_{i'},  i' \in \, _{\bar{r}_2}\XX$ are determined for the process $_{\bar{r}_2}\eta^{(\e)}(t)$ by condition ${\bf A}$ and relation (\ref{denojat}) in the same way as the transition sets $_{\bar{r}_1}\YY_{i'}, i' \in \, _{r_1}\XX$ are determined by condition ${\bf A}$ and relation (\ref{denojat})  for the process $_{\bar{r}_1}\eta^{(\e)}(t)$. Therefore,  the $(_{\bar{r}_2}l_{i'j'}^-, \, _{\bar{r}_2}l_{i'j'}^+)$-expansions for transition probabilities $_{\bar{r}_2}p_{i'j'}(\e), j' \in \,  _{\bar{r}_2}\YY_{i'}, i' \in \,  _{\bar{r}_2}\XX$ and $(_{\bar{r}_2}m_{i'j'}^-, \, _{\bar{r}_2}m_{i'j'}^+)$-expansions for expectations $_{\bar{r}_2}e_{i'j'}(\e), j' \in \,  _{\bar{r}_2}\YY_{i'}, i' \in \,  _{\bar{r}_2}\XX$ can be constructed by applying Algorithms 1, 3 and 5 to the 
 $(_{\bar{r}_1}l_{i'j'}^-, \, _{\bar{r}_1}l_{i'j'}^+)$-expansions for transition probabilities $_{\bar{r}_1}p_{i'j'}(\e), j' \in \,  _{\bar{r}_1}\YY_{i'}, i' \in \,  _{\bar{r}_1}\XX$ and $(_{\bar{r}_1}m_{i'j'}^-, \, _{\bar{r}_1}m_{i'j'}^+)$-expansions for expectations $_{\bar{r}_1}e_{i'j'}(\e), j' \in \,  _{\bar{r}_1}\YY_{i'}, i' \in \,  _{\bar{r}_1}\XX$. These expansions ara pivotal.

{\bf 7.3.} By continuing the above procedure of  phase space reduction for states $r_3, \ldots, r_{N-1}$,  we construct the semi-Markov process $_{\bar{r}_{N-1}}\eta^{(\e)}(t)$ with the phase space $_{\bar{r}_{N-1}}\XX = \XX \setminus \{ r_1, r_2, 
\ldots, r_{N-1} \} = \{ i \}$ (which is a one-point set), the transition probabilities of the embedded Markov chain $_{\bar{r}_{N-1}}p_{ii}(\e) = 1$, and the expectations of transition times $_{\bar{r}_{N -1}}e_{ii}(\e)$, which are determined by the transition probabilities and the expectations of transition times of the process $_{\bar{r}_{N-2}}\eta^{(\e)}(t)$ via relations (\ref{transit}) and (\ref{expectaga}). According Theorem 1, the expectations of hitting times $E_{ii}(\e)$ for the semi-Markov processes $\eta^{(\e)}(t)$,  $_{\bar{r}_1}\eta^{(\e)}(t), \ldots, \, _{\bar{r}_{N-1}}\eta^{(\e)}(t)$ coincide. According Theorems 2 and 4, the transition probabilities of the embedded Markov chain $_{\bar{r}_{N-1}}p_{ii}(\e) = 1$ and the expectations of transition times $_{\bar{r}_{N-1}}e_{ii}(\e)$ satisfy conditions ${\bf A}$ -- ${\bf F}$.  In this case, the transition set  $_{\bar{r}_{N-1}}\YY_i = \{ i \}$, for every $i \in \XX$. Therefore,  the $(_{\bar{r}_{N-1}}l_{i'j'}^-, \, _{\bar{r}_{N-1}}l_{i'j'}^+)$-expansions for transition probabilities $_{\bar{r}_{N-1}}p_{i'j'}(\e) = 1, j' \in \,  _{\bar{r}_{N-1}}\YY_{i'}, i' \in \,  _{\bar{r}_{N-1}}\XX$ (which take the form of relation (\ref{vopit})) and $(_{\bar{r}_{N-1}}m_{i'j'}^-$, $_{\bar{r}_{N-1}}m_{i'j'}^+)$-expansions for expectations $_{\bar{r}_{N-1}}e_{i'j'}(\e), j' \in \,  _{\bar{r}_{N-1}}\YY_{i'}$, $i' \in \,  _{\bar{r}_{N-1}}\XX$ can be constructed by applying Algorithms 1, 3 and 5 to the 
 $(_{\bar{r}_{N-2}}l_{i'j'}^-, \, _{\bar{r}_{N-2}}l_{i'j'}^+)$ -expansions for transition probabilities \, $_{\bar{r}_{N-2}}p_{i'j'}(\e), \, j' \in$ $_{\bar{r}_{N-2}}\YY_{i'}, \, i' \in \,  _{\bar{r}_{N-2}}\XX$ and $(_{\bar{r}_{N-2}}m_{i'j'}^-, \, _{\bar{r}_{N-2}}m_{i'j'}^+)$-expansions for expectations $_{\bar{r}_{N-2}}e_{i'j'}(\e), j' \in \,  _{\bar{r}_{N-2}}\YY_{i'}$, $i' \in \,  _{\bar{r}_{N-2}}\XX$.  These expansions ara pivotal.

{\bf 7.4.}  The semi-Markov process $_{\bar{r}_{N-1}}\eta^{(\e)}(t)$ has the one-point phase space $_{\bar{r}_{N-1}}\XX = \{ i \}$ and, thus, the transition probability $_{\bar{r}_{N-1}}p_{ii}(\e) \equiv 1$,  while the expectation of transition time   $_{\bar{r}_{N -1}}e_{ii}(\e) 
= E_{ii}(\e)$. The above algorithm of sequential reduction of phase space should be repeated for every $i \in \XX$. In this way, the Laurent asymptotic expansions for quantities  $E_{ii}(\e), i \in \XX$ can be written down. These asymptotic expansions have the following form, 
\begin{equation}\label{fina}
E_{ii}(\e) = \sum_{l = M_{ii}^-}^{M_{ii}^+} B_{ii}[l]\e^l + \hat{o}_i(\e^{M_{ii}^+}), i \in \XX,  
\end{equation}
where parameters $M_{ii}^{\pm} = \, _{\bar{r}_{N-1}}m_{ii}^{\pm}, i \in \XX$ and the coefficients $B_{ii}[l] = \, _{\bar{r}_{N-1}}b_{ii}[l]$, $l = M_{ii}^{-}, \ldots, M_{ii}^{+}, i \in \XX$, where  $_{\bar{r}_{N-1}}b_{ii}[l]$ are coefficients in the corresponding $(_{\bar{r}_{N-1}}m_{i'j'}^-, \, _{\bar{r}_{N-1}}m_{i'j'}^+)$-expansions for expectations $_{\bar{r}_{N-1}}e_{i'j'}(\e), j' \in \,  _{\bar{r}_{N-1}}\YY_{i'}$, $i' \in \,  _{\bar{r}_{N-1}}\XX$.  These expansions are pivotal.

It should be noted that, for every $n = 1, \ldots, N -1$, the reduced semi-Markov process $_{\bar{r}_{n}}\eta^{(\e)}(t)$ is invariant with respect to any permutation 
$\bar{r}'_n = (r'_1, \ldots, r'_n)$ of the the sequence $\bar{r}_n = (r_1, \ldots, r_n)$. 

Indeed, for every such permutation $\bar{r}'_n = (r'_1, \ldots, r'_n)$, the corresponding reduced semi-Markov process $_{\bar{r}'_{n}}\eta^{(\e)}(t)$ is constructed from the initial semi-Markov process $\eta^{(\e)}(t)$,  as the sequence of its states  at sequential moment of hitting into the same reduced phase space $_{\bar{r}'_n}\XX = \XX \setminus \{ r'_1, \ldots, r'_n \} = \, _{\bar{r}_n}\XX  = \XX \setminus \{ r_1, \ldots, r_n \}$ and times between sequential jumps  of the  reduced semi-Markov process $_{\bar{r}'_{n}}\eta^{(\e)}(t)$ which are  times between sequential hitting of the above  reduced space by the initial  semi-Markov process $\eta^{(\e)}(t)$. 

This implies that the expectation of transition time $_{\bar{r}_{n}}e_{i'j'}(\e)$ is, for every $i', j' \in \, _{\bar{r}_n}\XX$ and $n = 1, \ldots, N - 1$,  invariant with respect to any permutation $\bar{r}'_n = (r'_1, \ldots, r'_n)$ of the sequence $\bar{r}_n = (r_1, \ldots, r_n)$. 

Moreover, as follows from the Algorithms 1 -- 7, the expectation of transition time $_{\bar{r}_{n}}e_{i'j'}(\e)$  is a rational function of initial transition probabilities  $p_{ij}(\e), j \in \YY_i, i \in \XX$ and expectations $e_{ij}(\e), j \in \YY_i, i \in \XX$ (a quotient of two sums of products of some of these probabilities and expectations),  which, according the above remarks, is invariant with respect to any permutation $\bar{r}'_n = (r'_1, \ldots, r'_n)$ of the sequence $\bar{r}_n = (r_1, \ldots, r_n)$. 

By using identical arithmetical transformations (disclosure of brackets, imposition of a common factor out of the brackets,  bringing a fractional expression to a common denominator, permutation of summands or multipliers,  elimination of expression with equal absolute values and opposite signs in the sums and  elimination of equal expressions in the quotients, etc.) the rational function  $_{\bar{r}'_{n}}e_{i'j'}(\e)$ given by Algorithm 7 can be transformed in the rational function  $_{\bar{r}_{n}}e_{i'j'}(\e)$ given by Algorithm 7 and wise versa.

By Lemma 8, these transformations do not affect the corresponding asymptotic expansions for  expectation $_{\bar{r}_{n}}e_{i'j'}(\e)$ given by Algorithm 7, and, thus, these asymptotic expansions  are  invariant with respect to any permutation $\bar{r}'_n = (r'_1, \ldots, r'_n)$ of the sequence $\bar{r}_n = (r_1, \ldots, r_n)$.

The above remarks can be summarized in the following theorem.  \vspace{1mm}

{\bf Theorem 6.} {\em Let conditions ${\bf A}$ -- ${\bf F}$ hold for semi-Markov processes  $\eta^{(\e)}(t)$. Then, for every $i \in \XX$, the Laurent asymptotic expansion {\rm (\ref{fina})} for the expectation of hitting times $E_{ii}(\e)$ given by Algorithm 7 can be written down. This expansion is invariant with respect to the choice of 
permutation  $\langle r_1, \ldots, r_{N-1}, i \rangle$ of sequence $\langle 1, \ldots, N \rangle$, in the above algorithm.}  \vspace{1mm}

Let us now assume that conditions ${\bf A}$, ${\bf B'}$, ${\bf C}$ -- ${\bf E}$, ${\bf F'}$ hold for the semi-Markov process $\eta^{(\e)}(t)$.  \vspace{1mm}

{\bf Algorithm 8.} This is an algorithm for  computing upper bounds for remainders in asymptotic expansions for transition probabilities and expectation of sojourn times for semi-Markov processes with reduced phase spaces. 

\vspace{1mm}

{\bf 8.1.} Let $_{\bar{r}_1}\eta^{(\e)}(t) = \, _{r_1}\eta^{(\e)}(t)$ be be the reduced  semi-Markov process,  which is constructed as this is described in Step 7.1 of Algorithm 7. According Theorems 3 and 5, the  semi-Markov process $_{\bar{r}_1}\eta^{(\e)}(t)$ satisfies conditions 
${\bf A}$, ${\bf B'}$, ${\bf C}$ -- ${\bf E}$, ${\bf F'}$. Therefore,  $(_{\bar{r}_1}l_{i'j'}^-, \, _{\bar{r}_1}l_{i'j'}^+,  \, _{\bar{r}_1}\delta_{i'j'}, \, _{\bar{r}_1}G_{i'j'}, \, _{\bar{r}_1}\e_{i'j'})$-expansions for transition probabilities $_{\bar{r}_1}p_{i'j'}(\e), j' \in \,  _{\bar{r}_1}\YY_{i'}, i' \in \,  _{\bar{r}_1}\XX$ and $(_{\bar{r}_1}m_{i'j'}^-, \, _{\bar{r}_1}m_{i'j'}^+,  \, _{\bar{r}_1}\dot{\delta}_{i'j'}, \, _{\bar{r}_1}\dot{G}_{i'j'}, \, _{\bar{r}_1}\dot{\e}_{i'j'})$-expansions for expectations 
$_{\bar{r}_1}e_{i'j'}(\e), j'$ $\in \,  _{\bar{r}_1}\YY_{i'}, i' \in \,  _{\bar{r}_1}\XX$ can be constructed by applying Algorithms 1 -- 5 to the $(l_{i'j'}^-, l_{i'j'}^+,  \delta_{i'j'},  G_{i'j'},  \e_{i'j'})$-expansions for transition probabilities $p_{i'j'}(\e), j' \in \YY_{i'}, i' \in \XX$ and $(m_{i'j'}^-, m_{i'j'}^+, \dot{\delta}_{i'j'},  \dot{G}_{i'j'},  \dot{\e}_{i'j'})$-expansions for expectations $e_{i'j'}(\e)$, $j' \in \YY_{i'}, i' \in \XX$. 

{\bf 8.2.} Let $_{\bar{r}_2}\eta^{(\e)}(t)$ be the reduced semi-Markov process,  which is constructed as this is described in Step 7.2 of Algorithm 7.   According Theorems 3 and 5, the  semi-Markov process $_{\bar{r}_2}\eta^{(\e)}(t)$ satisfies conditions 
${\bf A}$, ${\bf B'}$, ${\bf C}$ -- ${\bf E}$, ${\bf F'}$. Therefore,  $(_{\bar{r}_2}l_{i'j'}^-, \, _{\bar{r}_2}l_{i'j'}^+,  \, _{\bar{r}_2}\delta_{i'j'}, \, _{\bar{r}_2}G_{i'j'}, \, _{\bar{r}_2}\e_{i'j'})$-expansions for transition probabilities $_{\bar{r}_2}p_{i'j'}(\e), j' \in \,  _{\bar{r}_2}\YY_{i'}, i' \in \,  _{\bar{r}_2}\XX$ and $(_{\bar{r}_2}m_{i'j'}^-, \, _{\bar{r}_2}m_{i'j'}^+,  \, _{\bar{r}_2}\dot{\delta}_{i'j'}$, $_{\bar{r}_2}\dot{G}_{i'j'}, \, _{\bar{r}_2}\dot{\e}_{i'j'})$-expansions for expectations 
$_{\bar{r}_2}e_{i'j'}(\e), j'$ $\in \,  _{\bar{r}_2}\YY_{i'}, i' \in \,  _{\bar{r}_2}\XX$ can be constructed by applying Algorithms 1 -- 5 to the \, $(_{\bar{r}_1}l_{i'j'}^-, \, _{\bar{r}_1}l_{i'j'}^+, \, _{\bar{r}_1}\delta_{i'j'}, \,  _{\bar{r}_1}G_{i'j'}$,  $_{\bar{r}_1}\e_{i'j'})$-expansions \, for \, transition \, probabilities $p_{i'j'}(\e), j' \in \YY_{i'}, i' \in  \, _{\bar{r}_1}\XX$ and $(_{\bar{r}_1}m_{i'j'}^-, \,  _{\bar{r}_1}m_{i'j'}^+, \, _{\bar{r}_1}\dot{\delta}_{i'j'},  \,  _{\bar{r}_1}\dot{G}_{i'j'}$,  $_{\bar{r}_1}\dot{\e}_{i'j'})$-expansions for expectations $e_{i'j'}(\e)$, $j' \in \, _{\bar{r}_1}\YY_{i'}, i' \in \, _{\bar{r}_1}\XX$. 

{\bf 8.3.} Finally, let $_{\bar{r}_{N-1}}\eta^{(\e)}(t)$ be the reduced semi-Markov process, which is constructed as this  is described in Step 7.2 of Algorithm 7.   According Theorems 3 and 5, the  semi-Markov process $_{\bar{r}_{N-1}}\eta^{(\e)}(t)$ 
satisfies conditions ${\bf A}$, ${\bf B'}$, ${\bf C}$ -- ${\bf E}$, ${\bf F'}$. Therefore,  $(_{\bar{r}_{N-1}}l_{i'j'}^-, \, _{\bar{r}_{N-1}}l_{i'j'}^+$,  
$_{\bar{r}_{N-1}}\delta_{i'j'}, \, _{\bar{r}_{N-1}}G_{i'j'}, \, _{\bar{r}_{N-1}}\e_{i'j'})$-expansions for transition probabilities $_{\bar{r}_{N-1}}p_{i'j'}(\e)$, $j' \in \,  
_{\bar{r}_{N-1}}\YY_{i'}, i' \in \,  _{\bar{r}_{N-1}}\XX$ and $(_{\bar{r}_{N-1}}m_{i'j'}^-, \, _{\bar{r}_{N-1}}m_{i'j'}^+,  \, _{\bar{r}_{N-1}}\dot{\delta}_{i'j'}, \, 
_{\bar{r}_{N-1}}\dot{G}_{i'j'}, \, _{\bar{r}_{N-1}}\dot{\e}_{i'j'})$-expansions for expectations 
$_{\bar{r}_{N-1}}e_{i'j'}(\e), j'$ $\in \,  _{\bar{r}_{N-1}}\YY_{i'}, i' \in \,  _{\bar{r}_{N-1}}\XX$ can be constructed by applying Algorithms 1 -- 5 to the  
$(_{\bar{r}_{N-2}}l_{i'j'}^-, \, _{\bar{r}_{N-2}}l_{i'j'}^+, \, _{\bar{r}_{N-2}}\delta_{i'j'}$,  $_{\bar{r}_{N-2}}G_{i'j'},  \, _{\bar{r}_{N-2}}\e_{i'j'})$-expansions for  transition \, probabilities $p_{i'j'}(\e), j' \in \YY_{i'}, i'$ $\in  \, _{\bar{r}_{N-2}}\XX$ and $(_{\bar{r}_{N-2}}m_{i'j'}^-, \,  _{\bar{r}_{N-2}}m_{i'j'}^+, \, 
_{\bar{r}_{N-2}}\dot{\delta}_{i'j'},  \,  _{\bar{r}_{N-2}}\dot{G}_{i'j'},  \,  _{\bar{r}_{N-2}}\dot{\e}_{i'j'})$-expansions for expectations $e_{i'j'}(\e)$, $j' \in \, 
_{\bar{r}_{N-2}}\YY_{i'}, i' \in \, _{\bar{r}_{N-2}}\XX$. 

{\bf 8.4.}  Finally, due to equalities  $_{\bar{r}_{N -1}}e_{ii}(\e) = E_{ii}(\e), i \in \XX$, we get that the asymptotic expansion (\ref{fina}) for expectations 
$E_{ii}(\e), i \in \XX$, given in the Step 7.4 of Algorithm 7, is a $(M_{ii}^{-}, M_{ii}^{+}, \delta^\circ_{ii}, G^\circ_{ii}, \e^\circ_{ii})$-expansion with parameters $M_{ii}^{-} = \, _{\bar{r}_{N-1}}m_{ii}^{-}, M_{ii}^{+} = \, _{\bar{r}_{N-1}}m_{ii}^{+},  \delta^\circ_{ii} = \, _{\bar{r}_{N-1}}\dot{\delta}_{ii}, 
G^\circ_{ii} = \, _{\bar{r}_{N-1}}\dot{G}_{ii}, \e^\circ_{ii} = \, _{\bar{r}_{N-1}}\dot{\e}_{ii}$. 

In this case, the invariance  of explicit upper bounds for  remainders given by Algorithm 8,  with respect to the choice of any permutation  
$\langle r_1, \ldots, r_{N-1}, i \rangle$ of sequence $\langle 1, \ldots, N \rangle$,  can not be guaranteed. 

However, Lemma 9 guarantees that the following inequalities hold for the parameters  $\delta^\circ_{ii}, i \in \XX$,
\begin{equation}\label{bound}
\delta^\circ_{ii} \geq \delta^\circ = \min_{j \in \YY_i, i \in \XX} ( \delta_{ij} \wedge \dot{\delta}_{ij}). 
\end{equation}

The following theorem takes place. \vspace{1mm}

{\bf Theorem 7.} {\em Let conditions ${\bf A}$, ${\bf B'}$, ${\bf C}$ -- ${\bf E}$, ${\bf F'}$ hold for semi-Markov processes  $\eta^{(\e)}(t)$. Then, for every $i \in \XX$, the $(M_{ii}^{-}, M_{ii}^{+})$-expansion {\rm (\ref{fina})} for the expectations of hitting times $E_{ii}(\e)$, given by Algorithm 7, is a $(M_{ii}^{-}, M_{ii}^{+}, \delta^\circ_{ii}$, $G^\circ_{ii}, \e^\circ_{ii})$-expansion, with parameters $\delta^\circ_{ii}, G^\circ_{ii}, \e^\circ_{ii}$ given in Algorithm 8. The inequality {\rm (\ref{bound})} holds for  parameters  $\delta^\circ_{ii}, i \in \XX$.}  \vspace{1mm}

{\bf 7.2. Asymptotic expansions  for stationary probabilities of perturbed semi-Markov processes}.
Let us recall relation (\ref{hittana}) for stationary probabilities of the semi-Markov process $\eta^{(\e)}(t)$,
\begin{equation}\label{hittana1}
\pi_i(\e) = \frac{e_i(\e)}{E_{ii}(\e)}, \ i \in \XX.
\end{equation}

{\bf Algorithm 9.}  This is an algorithm for constructing asymptotic expansions  for stationary probabilities of perturbed 
semi-Markov processes. \vspace{1mm}

{\bf 9.1.} Conditions ${\bf A}$ - ${\bf F}$ and proposition {\bf (i)} (the multiple summation rule) of Lemma 5, permits one can construct $(m_{i}^-, m_{i}^+)$-expansions for expectations $e_{i}(\e), i \in \XX$, which take the following forms, 
\begin{align}\label{expaaba}
e_{i}(\e) & = \sum_{j \in \YY_i} e_{ij}(\e) \nonumber \\ 
& = \sum_{j \in \YY_i} \big( \sum_{l = m_{ij}^-}^{m_{i}^+} b_{ij}[l]\e^l + \dot{o}_{ij}(\e^{m_{ij}^+}) \big) \nonumber \\
& = \sum_{l = m_{i}^-}^{m_{i}^+} b_{i}[l]\e^l + \dot{o}_i(\e^{m_{i}^+}),  \ i \in \XX,
\end{align}
where
\begin{equation}\label{paramek}
m_{i}^- = \min_{j \in \, \YY_i} m_{ij}^-, \ m_{i}^+ = \min_{j \in \, \YY_i} m_{ij}^+, \ i \in \XX,
\end{equation}
and
\begin{equation}\label{paramasek}
b_{i}[m_{i}^- + l] = \sum_{j \in \, Y_i} b_{ij}[m_{i}^- + l], \ l = 0, \ldots, m_{i}^+ - m_{i}^-,  \ i \in \XX,
\end{equation}
where $b_{ij}[m_{i}^- + l] = 0$, for $0 \leq l < m_{ij}^- - m_{i}^- , j \in \, Y_i,  \ i \in \XX$. 

The above asymptotic expansions are pivotal for all $i \in \XX$. 

{\bf 9.2.}  Conditions ${\bf A}$ -- ${\bf F}$, relation (\ref{hittana1}) and proposition {\bf (v)} (the division rule) of Lemma 3, 
permits us  construct  
$(n_{i}^-, n_{i}^+)$-expansions for stationary probabilities $\pi_{i}(\e), i \in \XX$, which  take the following forms, 
\begin{equation}\label{expaabanaba}
\pi_{i}(\e) =  \sum_{l = n_{i}^-}^{n_{i}^+} c_{i}[l]\e^l + o_i(\e^{n_{i}^+}), \ i \in \XX,
\end{equation}
where 
\begin{equation}\label{final}
n_i^- = m_{i}^- - M_{ii}^- , \   n^+_i = (m_{i}^+ - M^-_{ii}) \wedge (m_{i}^- + M^+_{ii} - 2 M^-_{ii}), \ i \in \XX,
\end{equation}
and
\begin{align}\label{finala}
c_{i}[n_i^- + l] & =  \frac{ b_{i}[m_{i}^- + l] - \sum_{1 \leq l' \leq l} B_{ii}[M_{ii}^- + l] \, c_{i}[n_i^- + l - l']}{B_{ii}[M_{ii}^-]},  
\nonumber \\
l & = 0, \ldots, n_i^+  - n_i^-,  \ i \in \XX. 
\end{align}

Since $\pi_{i}(\e)  > 0, i \in \XX, \e \in (0, \e_0]$, the asymptotic expansions (\ref{expaaba}) are pivotal, i.e., coefficients, 
\begin{equation}\label{parewtop}
c_{i}[n_i^-] =  b_{i}[m_{i}^-] / B_{ii}[M_{ii}^-] > 0,  \ i \in \XX. 
\end{equation} 

By the definition,  $e_i(\e) \leq E_{ii}(\e), \ i \in \XX, \e \in (0, \e_0]$. This implies that parameters $M_{i}^- \leq  m_{i}^-, i \in \XX$ and thus, parameters \begin{equation}\label{parewt}
n_i^- \geq 0, \ i \in \XX.
\end{equation}

Moreover, since $\sum_{i \in \XX} \pi_i(\e) = 1$, for every $\e \in (0,  \e_0]$, the parameters $n_i^{\pm}, i \in \XX$ and coefficients $c_{i}[l], l = n_{i}^-, \ldots, n_{i}^+, i \in \XX$ satisfies the following relations,
\begin{equation}\label{parew}
n^- = \min_{i \in \XX} n_i^- = 0,
\end{equation}
and
\begin{equation}\label{limpolew}
c[l] = \sum_{i \in \XX} c_{i}[l] = \left\{
\begin{array}{cl}
1 & \ \text{for} \ l  = 0, \\
0 & \ \text{for} \ 0 < l \leq  n^+ =  \min_{i \in \XX} n_{i}^+.  \\
\end{array}
\right.
\end{equation}

Let us introduce sets,
\begin{equation*}\label{limpolewas}
\XX_0 = \{ i \in \XX: n_i^-  = 0 \}.
\end{equation*}

By the above remarks, the following relation takes place,
\begin{equation}\label{limmu}
\pi_{i}(0) = \lim_{\e \to 0} \pi_{i}(\e) = \left \{
\begin{array}{cl}
c_{i}[0] > 0 & \ \text{if} \  i \in \XX_0, \\
0 & \ \text{if} \  i \notin  \XX_0.  \\
\end{array}
\right.
\end{equation}

{\bf Theorem 8.} {\em Let conditions ${\bf A}$ -- ${\bf F}$ hold for semi-Markov processes  $\eta^{(\e)}(t)$. Then, the $(n_{i}^-, n_{i}^+)$-expansions  
{\rm (\ref{expaabanaba})}, for the stationary probabilities $\pi_{i}(\e), i \in \XX$ given by Algorithm 9, can be written down. This expansion is invariant with respect to the choice of permutation  $\langle r_1, \ldots, r_{N-1}, i \rangle$ of sequence $\langle 1, \ldots, N \rangle$, in the above algorithm. Relations {\rm (\ref{parewtop}) -- 
(\ref{limmu})} hold for these expansions.} \vspace{1mm}

{\bf Algorithm 10.}  This is an algorithm for computing  upper bounds for remainders in asymptotic expansions  for stationary probabilities of perturbed semi-Markov processes. \vspace{1mm}

{\bf 10.1.} Conditions ${\bf A}$, ${\bf B'}$, ${\bf C}$ -- ${\bf E}$, ${\bf F'}$ and the proposition {\bf (i)} (the multiple summation  rule) of Lemma 6 imply that the $(m_{i}^-, m_{i}^+)$-expansions for expectations $e_{i}(\e), i \in \XX$ are  $(m_{i}^-, m_{i}^+, \dot{\delta}_i,  \dot{G}_i,  \dot{\e}_i)$-expansions, with parameters $ \dot{\delta}_i,  \dot{G}_i,  \dot{\e}_i, i \in \XX$ given by the following formulas,
\begin{equation*}
\dot{\delta}_{i} =    \min_{j \in \YY_i, \,  m^+_{ij} = m^+_i } \dot{\delta}_{ij}, \makebox[58mm]{} 
\end{equation*}
\begin{align*}
\dot{G}_{i} = &  \sum_{j \in \YY_i } \big( \dot{G}_{ij} \dot{\e}_{i}^{m^+_{ij} + \dot{\delta}_{ij} - m^+_{i}  - \dot{\delta}_{i} }  
+  \sum_{m^+_{i}  < j \leq m^+_{ij}}  |b_{ij}| \dot{\e}_{i}^{j - m^+_{i} - \dot{\delta}_{i}} \big), 
\end{align*}
\begin{equation}\label{dvadre}
\dot{\e}_{i} = \min_{j \in \YY_i} \dot{\e}_{ij}. \makebox[72mm]{}  
\end{equation} 

{\bf 10.2.}  Conditions ${\bf A}$, ${\bf B'}$, ${\bf C}$ -- ${\bf E}$, ${\bf F'}$ and the propositions  {\bf (iv)} (the reciprocal rule) and  {\bf (v)} (the division rule)  of Lemma 6 imply that the $(n_{i}^-, n_{i}^+)$-expansions for expectations $\pi_{i}(\e), i \in \XX$ are  $(n_{i}^-, n_{i}^+, \delta^*_i,  
G^*_i,  \e^*_i)$-expansions, with parameters $\delta^*_i,  G^*_i,  \e^*_i, i \in \XX$ given by the following formulas,
\begin{equation*}
\delta^*_i = \left\{
\begin{array}{ll}
\dot{\delta}_i & \ \text{if} \ n^+_i = m^+_i - M^-_{ii} < n^-_i + M^+_{ii}  - 2 M^-_{ii}, \\
\dot{\delta}_i \wedge \delta^\circ_{ii} & \ \text{if} \ n^+_i = m^+_i - M^-_{ii} < n^-_i + M^+_{ii}  - 2 M^-_{ii},  \makebox[28mm]{}   \\
\delta_B & \ \text{if} \ n^+_i  =  n^-_i + M^+_{ii}  - 2 M^-_{ii} < m^+_i - M^-_{ii},
\end{array}
\right. 
\end{equation*} 
\begin{align*}
G^*_i = & \, (\frac{B_{ii}[M^-_{ii}]}{2})^{-1} \big( \sum_{m^+_i \wedge (m^-_i + M^+_{ii} - M^-_{ii}) < l \leq m^+_i} |b_i[l]| (\e^*_i)^{l - n^+_i  - M^-_{ii} - \delta^*_i} \nonumber \\
& +  \sum_{m^+_i \wedge (m^-_i + M^+_{ii} - M^-_{ii}) < l + k, m^-_i \leq l \leq m^+_i, M^-_{ii} \leq k \leq  M^+_{ii}} |b_i[l]| |c_i[k]|  
(\e^*_i)^{l + k  - n^+_i - M^-_{ii}  - \delta^*_i}   \nonumber \\
& + \dot{G}_i (\e^*_i)^{m^+_i + \dot{\delta}_i - n^+_i - M^-_{ii} - \delta^*_i}  \nonumber \\
& + G^\circ_{ii} \sum_{n^-_i \leq k \leq n^+_i} |c_i[k]| (\e^*_i)^{k + M^+_{ii} + \delta^\circ_i  - n^+_i - M^-_{ii}  - \delta^*_i} \big),
\end{align*}
{\small
\begin{equation}\label{hoputrkop}
\e^*_i = \dot{\e}_i \wedge  \e^\circ_i \wedge \left\{
\begin{array}{ll}
 \frac{B_{ii}[M^-_{ii}]}{2} ( \sum_{M^-_{ii} < l \leq M^+_{ii}} |B_{ii}[l]| (\e^\circ_i)^{l - M^-_{ii} - 1} \vspace{2mm} \\
+  \ G^\circ_{ii} (\e^\circ_i)^{M^+_{ii} + \delta^\circ_i - M^-_{ii}  - 1})^{-1}  & \ \text{if} \ M^-_{ii} < M^+_{ii},   \vspace{2mm} \\
( \frac{B_{ii}[M^-_{ii}]}{2 G^\circ_i})^{\frac{1}{\delta^\circ_i}}  & \ \text{if} \ M^-_{ii} = M^+_{ii}.
\end{array}
\right.
\end{equation}
}

{\bf Theorem 9.} {\em Let conditions ${\bf A}$, ${\bf B'}$, ${\bf C}$ -- ${\bf E}$, ${\bf F'}$ hold for semi-Markov processes  $\eta^{(\e)}(t)$. Then, the $(n_{i}^-, n_{i}^+, \delta^*_i,  G^*_i,  \e^*_i)$-expansions  {\rm (\ref{expaabanaba})} for the stationary probabilities $\pi_{i}(\e), i \in \XX$ given by Algorithms 9 and 10 can be 
written down. The  inequalities $\delta^*_i \geq \delta^\circ, i \in \XX$ hold, where parameter $\delta^\circ$ is given in relation {\rm (\ref{bound})}.} \vspace{1mm}

{\bf 7.3. Laurent asymptotic  expansions for expectations of hitting times}. Algorithms presented above yields the Laurent asymptotic  expansions for expectations of hitting times $E_{ij}(\e), i, j \in \XX$. Indeed, let choose two states $i, j \in \XX$  and a chain of states   $\bar{r}_{N-2} = (r_1, \ldots, r_{N-2}), r_1, \ldots, r_{N-2} \neq i, j$. 

According the Theorem 1 and Algorithm 
8 the expectations $E_{ij}(\e)$ coincides for the initial semi-Markov process $\eta^{(\e)}(t)$ and the semi-Markov process 
$_{\bar{r}_{N-2}}\eta^{(\e)}(t)$.  The semi-Markov process $_{\bar{r}_{N-2}}\eta^{(\e)}(t)$ has a two-points phase space $_{\bar{r}_{N-2}}\XX = \{ i, j \}$. The expectations of hitting times  $E_{i' j'}(\e), i' \in \{ i, j \}$ can be found by solving, for every $j' \in \{ i, j \}$, the system of (two, in this case) linear equations (\ref{systera}) that yields the following formulas, for every $j' \in \{ i, j \}$,
\begin{equation}\label{nopla}
\left \{
\begin{array}{ll}
E_{i' j'}(\e) & = \, _{\bar{r}_{N-2}}e_{i'}(\e) \cdot \frac{1}{_{\bar{r}_{N-2}}p_{i' j'}(\e)},  \vspace{2mm} \\
E_{j' j'}(\e) & = \, _{\bar{r}_{N-2}}e_{j'}(\e)  + \, _{\bar{r}_{N-2}}e_{i'}(\e) \cdot \frac{_{\bar{r}_{N-2}}p_{j' i'}(\e) }{_{\bar{r}_{N-2}}p_{i' j'}(\e)},
\end{array}
\right.
\end{equation}
where $i' \neq j'$ in both equations in ({\ref{nopla}) and, 
\begin{equation}\label{nytr}
_{\bar{r}_{N-2}}e_{i'}(\e) =  \, _{\bar{r}_{N-2}}e_{i' i}(\e) + 
\, _{\bar{r}_{N-2}}e_{i' j}(\e), \ i' \in \{i, j \}.
\end{equation}

The  corresponding asymptotic expansions for $E_{ij}(\e)$  can be constructing by using the  asymptotic expansions for transition probabilities  $p_{j' i'}(\e)$ and expectations $_{\bar{r}_{N-2}}e_{i'}(\e)$ given in Algorithms 7 and 8 and the operational rules  for Laurent asymptotic expansions presented in Lemmas 1 -- 9.  \\

{\bf 8. Future studies and bibliographical remarks} \\

 {\small In this section, we present some directions for future studies and short bibliographical remarks concerned works in the area. \vspace{1mm}

{\bf 8.1. Directions for future studies}. The method of sequential reduction of a phase space presented in the paper  can also be applied  for getting asymptotic expansions for high order power  and exponential moments of hitting times,  for nonlinearly perturbed semi-Markov processes.

In the present paper, we consider the model, where the pre-limiting perturbed  semi-Markov processes have a phase space which is one class of communicative states, while the limiting unperturbed semi-Markov process has a phase space which consists of one or several classes of communicative states and possibly a 
class of transient states.  However, the method of sequential reduction of the phase space can also be applied to nonlinearly perturbed semi-Markov processes with absorption and, therefore, to the model, where the pre-limiting semi-Markov processes also have a phase space, which consists of several classes of communicative states and a class of transient states.  

We are quite sure that combination of results in the above two directions with the methods of asymptotic analysis for nonlinearly perturbed regenerative processes  developed in Silvestrov (1995, 2007, 2010) and  Gyllenberg and Silvestrov (1998, 1999a, 2000a, 2008) will make it possible to expand results concerned  asymptotic expansions for quasi-stationary distributions and other characteristics for nonlinearly perturbed semi-Markov processes with absorption, where the limiting semi-Markov process has a phase space which consists of one class of communicative states and a class of transient states, to a general case,  where the limiting semi-Markov process has a  phase space, which consists of several classes of communicative states and a class of transient states.

The problems of aggregation of steps in the time-space screening procedures for semi-Markov processes, tracing  pivotal orders for different groups of states as well as  getting explicit matrix formulas, for coefficients and parameters of upper bounds for remainders in the corresponding asymptotic expansions for stationary distributions and moments of hitting times,  do require additional studies.  It can be expected that such formulas can be obtained, for example, for birth-death type  semi-Markov processes, for which the proposed algorithms of reduction of a phase space preserve the birth-death structure for reduced semi-Markov processes. 

We are going to present results concerned  Laurent asymptotic expansions for power and exponential moments of hitting times, quasi-stationary distributions and  explicit formulas for coefficients and parameters of upper bounds for remainders for some specific classes of semi-Markov models, as well as applications to some models of population genetics, information networks and queuing systems, in future publications.  \vspace{1mm}

{\bf 8.2. Bibliographical remarks}.
Note first of all that the model of perturbed discrete time Markov chains, at least, in the most difficult case of so-called singularly perturbed  Markov chains and semi-Markov processes  with absorption and asymptotically uncoupled phase spaces, attracted attention of researchers in the mid of the 20th century.  

{\bf (1)} The first  works related to asymptotical problems for the above models are Meshalkin (1958),  Simon, and Ando (1961), Hanen (1963a, b, c, 1964), Seneta (1967, 1968a, b, 1973), Schweitzer (1968), Korolyuk (1969),   Silvestrov (1969, 1970, 1971a, b, 1972a, b, 1974),  Anisimov  (1970, 1971a, b,  1973a, b), Korolyuk and Turbin (1970, 1972), Gusak and  Korolyuk (1971),   Turbin (1971, 1972), Korolyuk, Penev and  Turbin (1972), Kovalenko (1973, 1975),  Poli\v s\v cuk and  Turbin (1973),  Korolyuk, Brodi and Turbin (1974), Pervozvanski\u\i \, and  Smirnov (1974),   Courtois (1975) and Ga\u\i tsgori and  Pervozvanski\u\i \, (1975).  

{\bf (2)}  Convergence results, for distributions and moments of hitting times, eigenvalues, eigenvectors, stationary and quasi-stationary distributions, Perron roots, coefficients of ergodicity, etc. have been studied in works by
Meshalkin (1958), Hanen, (1963a, b, c, 1964), Seneta (1967, 1968a, b, 1973, 2006), Schweitzer (1968, 1984), Korolyuk (1969, 1989),  Silvestrov (1969, 1970, 1971a, 1972a, 1974, 1976, 1978, 1979a, b, 1981, 2000),  Anisimov (1970, 1971a, b, 1973a, b, 1980, 1986, 1988, 2008), Korolyuk and Turbin (1970, 1972, 1976, 1978),    Gusak and  Korolyuk (1971), Turbin (1971),  Korolyuk, Penev and  Turbin (1972),  Kovalenko (1973, 1975),  Korolyuk, Brodi and Turbin (1974), Ga\u\i tsgori and  Pervozvanskiy (1975, 1983),  Allen,  Anderssen and Seneta (1977), Kaplan (1979, 1980),  Korolyuk,  Turbin and  Tomusjak (1979),  Shurenkov (1980a, b), Anisimov and  Chernyak (1982), Anisimov, Vo\u\i na and  Lebedev (1983), Coderch, Willsky, Sastry and Casta\~{n}on (1983),  Korolyuk, D. and Silvestrov (1983, 1984),  Stewart (1983, 1984), Koury,  McAllister and  Stewart (1984),   McAllister,  Stewart and  Stewart, W. (1984),  Cao and  Stewart  (1985),  Kartashov (1986, 1987,  1996b), Burnley (1987), Gibson and  Seneta (1987), Haviv (1987),   Haviv, Ritov and  Rothblum (1987), Rohlichek (1987), Rohlicek and  Willsky (1988a, b), Silvestrov  and Velikii (1988), Hunter (1991a),  Latouche (1991), Pollett and  Stewart (1994),  Hoppensteadt, Salehi and Skorokhod (1996a), Kalashnikov  (1997), Korolyuk and Limnios (1998, 1999, 2000, 2002, 2004a, b, 2005), Marek and Mayer (1998), Yin and  Zhang (1998, 2005, 2013), Craven (2003), Hern\'{a}ndez-Lerma and Lasserre (2003), Yin, Zhang and Badowski (2003),  Silvestrov and Drozdenko (2005, 2006), Drozdenko (2007a, b, 2009), Barbour and Pollett (2010, 2012) and Meyer (2015). 

{\bf (3)} Rates of convergence, errors of approximation, sensitivity  and related stability theorems for Markov chains and related models of stochastic processes have been studied in works by Schweitzer (1968, 1986, 1987, 1991), Silvestrov (1969, 1971b, 1972b),  Seneta (1973, 1984, 1988a, b, 1991, 1993, 2006), Courtois (1975, 1982),  Ga\u\i tsgori and  Pervozvanskiy (1975, 1983), Kovalenko (1975), Kalashnikov (1978, 1997a, 1997b), Louchard and  Latouche (1978, 1982 1990), Berman and Plemmons (1979, 1994), Meyer (1980, 1994),  Kalashnikov and Anichkin (1981), Bobrova (1983), Stewart (1983, 1984a, 1990, 1991, 1993a, b, 1998, 2001, 2003),  Courtois and Semal (1984a, b, c, 1991), Haviv and  Rothblum (1984),  Haviv and  Van der Heyden (1984), Koury, McAllister and  Stewart (1984), McAllister, Stewart, G. and Stewart, W. (1984), Funderlic and Meyer (1985), Kartashov (1985a, b, c, 1986, 1987, 1988, 1996a, 1996b, 2000, 2005),  Vantilborgh (1985), Haviv (1986, 1992, 2006), Haviv and  Ritov (1986, 1994),  Rohlichek (1987),  Rohlicek and  Willsky (1988a, b),  Stewart and  Sun (1990), Hunter  (1991a, b, 2005, 2014), Stewart and Zhang  (1991), Hassin and Haviv  (1992, 1999),  Barlow (1993), Meyn and Tweedie (1993, 2009), Lasserre (1994),  Pollett and Stewart (1994), Stewart, G., Stewart, W. and  McAllister (1994),  Borovkov (1998),  Yin and Zhang (1998, 2003, 2005, 2013), Li,  Yin, G., Yin, K. and Zhang (1999),  Craven (2003), Kontoyiannis and Meyn (2003),  Zhang and  Yin (2004),  Guo  (2006) and Sirl,  Zhang and Pollett (2007).  
     
{\bf (4)} Asymptotic expansions for distributions of hitting times, moments of hitting times, resolvents,  eigenvalues, eigenvectors, stationary and quasi-stationary distributions, Perron roots, etc., have been studied in works by Turbin (1972),  Poli\v s\v cuk and  Turbin (1973), Koroljuk, Brodi and Turbin (1974), Pervozvanski\u\i \, and  Smirnov (1974), Courtois  and  Louchard (1976),  Korolyuk and Turbin  (1976, 1978), Courtois  (1977), Latouche and  Louchard  (1978), Kokotovi\'{c},  Phillips and  Javid (1980), Korolyuk, Penev and  Turbin (1981),  Phillips and  Kokotovi\'{c} (1981), Delebecque (1983),  Abadov (1984),  Kartashov (1985d, 1996b), Haviv (1986),  Korolyuk (1989),   Stewart and  Sun (1990), Silvestrov and  Abadov (1991, 1993), Haviv, Ritov and  Rothblum (1992), Haviv and  Ritov (1993),  Schweitzer and Stewart (1993), Silvestrov (1995, 2007, 2010), Englund and  Silvestrov (1997),  Gyllenberg and Silvestrov (1998, 1999a, 2000a, 2008), Korolyuk and  Limnios (1998, 1999, 2000, 2002, 2004a, b, 2005), Stewart (1998, 2001), Yin and  Zhang (1998, 2003, 2005, 2013), Avrachenkov (1999, 2000), Avrachenkov and  Lasserre (1999), Korolyuk, V.S. and Korolyuk, V.V. (1999), Englund (2000, 2001), Yin, G., Zhang,  Yang and Yin, K. (2001), Avrachenkov and  Haviv (2003, 2004), Craven (2003) and  Avrachenkov,  Filar and Howlett  (2013).

{\bf (5)} Asymptotic expansions for  other characteristics of Markov type processes are presented in works by
Nagaev (1957, 1961), Leadbetter (1963),  Poli\v s\v cuk and  Turbin(1973),  Quadrat  (1983), Abadov (1984), Silvestrov and  Abadov (1984, 1991, 1993),  Stewart and Sun (1990),  Kartashov (1996b, 2013),  Khasminskii, Yin and  Zhang (1996a, b), Wentzell (1996, 1999), Yin and  Zhang  (1996b), Cao (1998),  Gyllenberg and  Silvestrov, (1998, 1999a, 2000a, 2000c, 2008), Fuh and  Lai (2001),  Kontoyiannis and Meyn (2003),  Fuh (2004, 2007), Samo\u\i lenko (2006a, b), Silvestrov  (1995, 2007, 2010), Ni  (2010a, b, 2011, 2012, 2014). Ni, Silvestrov and Malyarenko (2008),  Albeverio,  Koroliuk and  Samoilenko   (2009),  Petersson (2013a, b, 2014),  Avrachenkov,  Filar and Howlett  (2013)  and Silvestrov and  Petersson (2013), Silvestrov D. and Silvestrov, S. (2015). 

{\bf (6)} We would like esspesially to mention books including problems on perturbed Markov chains, semi-Markov processes and related problems. These are
 Seneta (1973, 2006), Silvestrov (1974), Korolyuk and Turbin  (1976, 1978), Courtois (1977), Kalashnikov  (1978, 1997b), Anisimov (1988, 2008), Stewart and  Sun (1990), Korolyuk and Swishchuk (1992), Meyn and  Tweedie (1993, 2009), Kartashov (1996b),   Borovkov (1998), Stewart (1998, 2001),Yin and Zhang (1998, 2005, 2013), Korolyuk, V.S. and Korolyuk, V.V. (1999),  Bini, Latouche and Meini (2005),  Koroliuk and  Limnios (2005), Gyllenberg and  Silvestrov (2008), and   Avrachenkov, Filar and Howlett (2013).   

{\bf (7)} General results of perturbation theory of matrices and linear operators are presented in works by
Vishik and Lyusternik (1960), Stewart (1969, 1973, 1979, 1984b, 1998, 2001),  Plotkin and Turbin (1971, 1975), Korolyuk and Turbin (1976, 1978),  Berman and Plemmons (1979, 1994), Wentzell and Freidlin (1979),   Haviv (1988), Meyer and Stewart (1988), Bielecki and  Stettner (1989), Delebecque (1990),  Stewart and Sun (1990), Hunter (1991b, 2014), Haviv and  Ritov (1994), Lasserre (1994),  Hoppensteadt, Salehi and  Skorokhod (1997),  Avrachenkov (1999),  Korolyuk, V.S. and Korolyuk, V.V. (1999), Li and Stewart (2000), Avrachenkov,  Haviv and  Howlett (2001), Howlett and  Avrachenkov (2001), Howlett,  Pearce and  Torokhti (2003), Torokhti, Howlett and  Pearce (2003), Verhulst ( 2005), Kartashov (1996b), Avrachenkov, Pearce and  Ejov (2009),  Howlett,  Albrecht and Pearce (2010, 2014), Albrecht, Howlett and  Pearce  (2011), and Avrachenkov and  Lasserre (2013).  In particular, we would like to mention some books, which contains materials on general perturbation matrix and operator theory. These are Erd\'{e}lyi (1956), Kato (1966, 2013),  Cole (1968), Korolyuk and Turbin (1976, 1978), Wentzell and  Freidlin (1979), Kevorkian and  Cole (1981, 1985, 1996, 2011), Baumg{\"a}rtel (1985),  Stewart (1998, 2001),  Korolyuk, V.S. and Korolyuk, V.V. (1999), Konstantinov,  Gu, Mehrmann and Petkov  (2003),  Verhulst (2005), Gyllenberg and Silvestrov (2008), and Avrachenkov,  Filar and Howlett  (2013).

{\bf (8)} Applications of results on perturbed Markov type processes  to the control theory, decision processes, Internet, queuing theory, mathematical genetics, population dynamics and epidemic models, insurance and financial mathematics are presented in works by
Simon and Ando (1961), Courtois  (1977), Kalashnikov  (1978, 1997b), Delebecque and  Quadrat (1981),  Quadrat (1983), Schweitzer (1984),  Anisimov, Zakusilo and Donchenko (1987), Latouche (1988), Pervozvanskii and  Gaitsgori (1988), Meyer (1989), Ho and  Cao (1991), Gyllenberg and  Silvestrov (1994, 1999b, 2000b, 2008, 2014), Pollett and Stewart (1994), Abbad and  Filar (1995),  Hoppensteadt, Salehi and Skorokhod (1996b),  Kovalenko,  Kuznetsov, and Pegg (1997), Borovkov (1998), Yin and Zhang (1998, 2005, 2013), Englund (1999a, b),   Yin, G., Zhang, Yang and Yin, K. (2001), Avrachenkov,  Filar and  Haviv (2002), Altman, Avrachenkov and  N\'{u}\~{n}ez-Queija (2004),  Langville and Meyer (2006), Avrachenkov,  Litvak, and Son Pham (2007, 2008),  Andersson and  Silvestrov, S. (2008), Anisimov (2008), Konstantinov and Petkov (2008),  Avrachenkov,  Borkar and  Nemirovsky (2010),   Barbour and Pollett (2010, 2012), H\"{o}ssjer (2011, 2014), Engstr\"{o}m and  Silvestrov, S. (2014), H\"{o}ssjer and  Ryman (2014), Ni (2014), Petersson  (2014), Silvestrov (2014, 2015).  In particular, we would like to mention some books in this area that are Kalashnikov  (1978, 1997b), Anisimov, Zakusilo and Donchenko (1987), Pervozvanski\u\i \, and  Gaitsgori (1988), Kovalenko,  Kuznetsov, and Pegg (1997),  (1998), Anisimov (2008), Gyllenberg and  Silvestrov (2008), Avrachenkov,  Filar and Howlett  (2013),  and  Silvestrov  (2014, 2015). 

{\bf (9)}  Exact and related approximative computational methods for stationary and quasi-stationary distributions of Markov chains and semi-Markov processes and related problems are presented in works by Romanovski\u\i \, (1949), Feller (1950, 1968),  Kemeny and Snell (1960), Golub and  Seneta (1973), Seneta (1973, 2006),  Paige,  Styan and  Wachter (1975),  Silvestrov  (1980a, b, 1996), Chatelin and  Miranker (1984), Harrod and  Plemmons (1984),  Schweitzer (1984, 1991), Grassman, Taksar and  Heyman (1985), Schweitzer,  Puterman and Kindle (1985),   Sheskin (1985), Hunter (1986, 1991a),  Schweitzer and  Kindle (1986), Feinberg and Chiu (1987), Haviv (1987, 1992),  Haviv, Ritov and  Rothblum (1987), Latouche and Ramaswami (1987, 1999), Sumita and Reiders (1988), Mattingly and  Meyer (1991), Stewart, W. (1994), Kim and  Smith (1995), Stewart (1998, 2001),  Kartashov (2000),  Meyer (2000),  H\"{a}ggstr\"{o}m (2002), Bini, Latouche and  Meini (2005), Golub and Van Loan (2013), Silvestrov, Manca and Silvestrova (2014). In particular, we would like to mention some related books that are  Romanovski\u\i \, (1949), Feller (1950, 1968),  Kemeny and Snell (1960), Golub and  Seneta (1973), Seneta (1973, 2006), Berman and  Plemmons (1979, 1994), 
Silvestrov (1980a), Meyer (2000), H\"{a}ggstr\"{o}m (2002), Bini, Latouche and  Meini (2005), Meyn and  Tweedie (1993, 2009), Hern\'{a}ndez-Lerma and Lasserre,  (2003), Gyllenberg and  Silvestrov (2008), and Collet,  Mart\'{\i}nez and San Mart\'{\i}n (2013).
}

\end{document}